\def\ifundefined#1#2{\expandafter\ifx\csname#1\endcsname\relax\input #2\fi}
\input amssym.def 
\def\Ma{{\cal\char'101}}
\def\Mb{{\cal\char'102}}
\def\Md{{\cal\char'104}}
\def\Me{{\cal\char'105}}
\def\Mf{{\cal\char'106}}
\def\Mg{{\cal\char'107}}
\def\Mh{{\cal\char'110}}
\def\Mi{{\cal\char'111}}
\def\Mj{{\cal\char'112}}
\def\Mm{{\cal\char'115}}
\def\Mn{{\cal\char'116}} 
\def\Mo{{\cal\char'117}}
\def\Mr{{\cal\char'122}}
\def\Ms{{\cal\char'123}}
\def\Mt{{\cal\char'124}}
\def\Mv{{\cal\char'126}}

\def\pmb#1{\setbox0=\hbox{$#1$}      
     \kern-.025em\copy0\kern-\wd0
     \kern.05em\copy0\kern-\wd0
     \kern-.025em\box0}

\def\endproof{$\hfill \square$}

\def\abs#1{\left\vert #1\right\vert}    

\def\Cross{\bigm| \kern-5.5pt \not \ \, }
\def\cross{\mid \kern-5.0pt \not \ \, }             
\def\notto{\hbox{$~\rightarrow~\kern-1.5em\hbox{/}\ \ $}}

\def\al{\alpha}
\def\be{\beta}
\def\ga{\gamma}

\def\dl{\delta}
\def\vph{\varphi}
\def\vep{\varepsilon}

\def\Th{\Theta}
\def\vth{\vartheta}
\def\sg{\sigma}

\hyphenation{math-ema-ticians}
\hyphenation{pa-ra-meters}
\hyphenation{pa-ra-meter}
\hyphenation{lem-ma}
\hyphenation{lem-mas}
\hyphenation{to-po-logy}
\hyphenation{to-po-logies}
\hyphenation{homo-logy}
\hyphenation{homo-mor-phy}

\def\nSigma{\Sigma \kern-8.3pt \bigm|\,}

\def\got#1{\hbox{\teneuler #1}}

\font\teneufm=eufm10
\font\eighteufm=eufm8
\font\fiveeufm=eufm5

\newfam\eufam
\textfont\eufam=\teneufm
\scriptfont\eufam=\eighteufm
\scriptscriptfont\eufam=\fiveeufm

\def\got{\fam=\eufam\teneufm}

\def\boxit#1{\vbox{\hrule\hbox{\vrule\kern2.0pt
       \vbox{\kern2.0pt#1\kern2.0pt}\kern2.0pt\vrule}\hrule}}

\def\vlra#1{\hbox{\kern-1pt
       \hbox{\raise2.38pt\hbox{\vbox{\hrule width#1 height0.26pt}}}
       \kern-4.0pt$\rightarrow$}}

\def\vlla#1{\hbox{$\leftarrow$\kern-1.0pt
       \hbox{\raise2.38pt\hbox{\vbox{\hrule width#1 height0.26pt}}}}}

\def\vlda#1{\hbox{$\leftarrow$\kern-1.0pt
       \hbox{\raise2.38pt\hbox{\vbox{\hrule width#1 height0.26pt}}}
       \kern-4.0pt$\rightarrow$}}

\def\longra#1#2#3{\,\lower3pt\hbox{${\buildrel\mathop{#2}
\over{{\vlra{#1}}\atop{#3}}}$}\,}

\def\longla#1#2#3{\,\lower3pt\hbox{${\buildrel\mathop{#2}
\over{{\vlla{#1}}\atop{#3}}}$}\,}

\def\longda#1#2#3{\,\lower3pt\hbox{${\buildrel\mathop{#2}
\over{{\vlda{#1}}\atop{#3}}}$}\,}

\def\overrightharpoonup#1{\vbox{\ialign{##\crcr
	$\rightharpoonup$\crcr\noalign{\kern-1pt\nointerlineskip}
	$\hfil\displaystyle{#1}\hfil$\crcr}}}

\catcode`@=11
\def\@@dalembert#1#2{\setbox0\hbox{$#1\rm I$}
  \vrule height.90\ht0 depth.1\ht0 width.04\ht0
  \rlap{\vrule height.90\ht0 depth-.86\ht0 width.8\ht0}
  \vrule height0\ht0 depth.1\ht0 width.8\ht0
  \vrule height.9\ht0 depth.1\ht0 width.1\ht0 }
\def\dalembert{\mathord{\mkern2mu\mathpalette\@@dalembert{}\mkern2mu}}

\def\@@varcirc#1#2{\mathord{\lower#1ex\hbox{\m@th${#2\mathchar\hex0017 }$}}}
\def\varcirc{\mathchoice
  {\@@varcirc{.91}\displaystyle}{\@@varcirc{.91}\textstyle}
{\@@varcirc{.45}\scriptscriptstyle}}
\catcode`@=12

\font\tensf=cmss10 \font\sevensf=cmss8 at 7pt
\newfam\sffam
\textfont\sffam=\tensf\scriptfont\sffam=\sevensf

\input amssym.def
\input amssym
\magnification=1200

\font\bigsll=cmsl10 scaled\magstep3
\tolerance=500
\overfullrule=0pt
\hsize=6.40 true in
\hoffset=.10 true in
\voffset=0.10 true in
\vsize=8.70 true in
\null
\centerline{\bigsll Crystalline Boundedness Principle}
\bigskip\smallskip
\centerline{\sl Adrian Vasiu, U of Arizona}
\bigskip\smallskip
\centerline{04/17/02; refinements 5/27/03; revisions 5/25/04 and 5/27/05; final version 12/15/05}
\bigskip\noindent
{\bf Abstract}. We prove that an $F$-crystal $(M,\vph)$ over an algebraically closed field $k$ of characteristic $p>0$ is determined by $(M,\vph)$ mod $p^n$, where $n\ge 1$ depends only on the rank of $M$ and on the greatest Hodge slope of $(M,\vph)$. We also extend this result to triples $(M,\vph,G)$, where $G$ is a flat, closed subgroup scheme of ${\bf GL}_M$ whose generic fibre is connected and has a Lie algebra normalized by $\vph$. We get two purity results. If ${\got C}$ is an $F$-crystal over a reduced ${\bf F}_p$-scheme $S$, then each stratum of the Newton polygon stratification of $S$ defined by ${\got C}$, is an affine $S$-scheme (a weaker result was known before for $S$ noetherian). The locally closed subscheme of the Mumford scheme ${\Ma_{d,1,N}}_k$ defined by the isomorphism class of a principally quasi-polarized $p$-divisible group over $k$ of height $2d$, is an affine ${\Ma_{d,1,N}}_k$-scheme. 
\bigskip\noindent
{\bf R\'esum\'e}. Nous prouvons qu'un $F$-cristal $(M,\vph)$ d\'efini sur un corps $k$ alg\'ebriquement clos de caract\'eristique $p>0$ est d\'etermin\'e par $(M,\vph)$ mod $p^n$,
o\`u $n\ge 1$ d\'epend seulement du rang de $M$ et de la plus grand
pente de Hodge de $(M,\vph)$. On \'etend ce r\'esultat aux triplets $(M,\vph,G)$, o\`u $G$ est un sous-groupe 
ferm\'e et plat de ${\bf GL}_M$ dont la fibre g\'en\'erique est connexe et a une alg\`ebre de Lie normalis\'ee par $\vph$. Nous obtenons deux r\'esultats de puret\'e. Si ${\got C}$ est un $F$-cristal sur un ${\bf F}_p$-sch\'ema r\'eduit $S$, alors chaque strate de la stratification du polygone de Newton de $S$ d\'efini par ${\got C}$ est un $S$-sch\'ema affine (un r\'esultat moins g\'en\'eral \'etait d\'ej\`a connu pour $S$ noeth\'erien). Le sous-sch\'ema localement ferm\'e du sch\'ema de Mumford
${\Ma_{d,1,N}}_k$ d\'efini par la classe d'isomorphisme d'un
groupe $p$-divisible principalement quasi-polaris\'e sur
$k$ de hauteur $2d$ est un ${\Ma_{d,1,N}}_k$-sch\'ema affine.
\bigskip\smallskip
\centerline{\bigsll {\bf \S 1 Introduction}}
\bigskip\smallskip
Let $p\in{\bf N}$ be a prime. Let $k$ be a perfect field of characteristic $p$. Let $\bar k$ be an algebraic closure of $k$. Let $W(k)$ be the Witt ring of $k$. Let $B(k):=W(k)[{1\over p}]$ be the field of fractions of $W(k)$. Let $\sg:=\sg_k$ be the Frobenius automorphism of $k$, $W(k)$, and $B(k)$. A group scheme $H$ over ${\rm Spec}(W(k))$ is called integral if $H$ is flat over ${\rm Spec}(W(k))$ and $H_{B(k)}$ is connected (i.e. if the scheme $H$ is integral). Let ${\rm Lie}(H_{B(k)})$ be the Lie algebra over $B(k)$ of $H_{B(k)}$. If $H$ is smooth over ${\rm Spec}(W(k))$, let  ${\rm Lie}(H)$ be the Lie algebra over $W(k)$ of $H$. If $O$ is a free module of finite rank over some commutative ${\bf Z}$-algebra $R$, let ${\bf GL}_O$ be the group scheme over ${\rm Spec}(R)$ of linear automorphisms of $O$. 
\smallskip
Let $(r,d)\in{\bf N}\times ({\bf N}\cup\{0\})$, with $r\ge d$. Let $D$ be a $p$-divisible group over ${\rm Spec}(\bar k)$ of height $r$ and dimension $d$. It is well known that if $d\in\{0,1,r-1,r\}$, then:
\medskip\noindent
{\bf (*)} {\it $D$ is uniquely determined up to isomorphism by its $p$-torsion subgroup scheme $D[p]$}. 
\medskip
But (*) does not hold if $2\le d\le r-2$. In 1963 Manin published an analogue of (*) for $2\le d\le r-2$ but unfortunately he separated it into three parts (see [28, p. 44, 3.6, and 3.8] and below). Only recently, this paper and [36] contain explicit analogues of (*) for $2\le d\le r-2$. The two main reasons for this delay in the literature are: (i) the widely spread opinion, which goes back more than 40 years, that $p$-divisible groups involve an infinite process, and (ii) the classification results of [28, p. 44] were rarely used. Our point of view is that $F$-crystals in locally free sheaves of finite rank over many ${\rm Spec}(k)$-schemes $Y$ involve a {\it bounded infinite process}. In this paper we give meaning to this point of view for the case $Y={\rm Spec}(\bar k)$. We start with few definitions. 
\medskip\smallskip
{\bf 1.1. Definitions.} {\bf (a)} By a {\it latticed $F$-isocrystal with a group} over $k$ we mean a triple $(M,\vph,G)$, where $M$ is a free $W(k)$-module of finite rank, where $\vph$ is a $\sg$-linear automorphism of $M[{1\over p}]$, and where $G$ is an integral, closed subgroup scheme of ${\bf GL}_M$, such that the Lie subalgebra ${\rm Lie}(G_{B(k)})$ of ${\rm End}(M[{1\over p}])$ is normalized by $\vph$. Here we denote also by $\vph$ the $\sg$-linear (algebra) automorphism of ${\rm End}(M[{1\over p}])$ that takes $e\in {\rm End}(M[{1\over p}])$ into $\vph\circ e\circ\vph^{-1}\in {\rm End}(M[{1\over p}])$. If $G={\bf GL}_M$, then often we do not mention $G$ and we omit ``with a group". 
\smallskip
{\bf (b)} By an isomorphism between two latticed $F$-isocrystals with a group $(M_1,\vph_1,G_1)$ and $(M_2,\vph_2,G_2)$ over $k$ we mean a $W(k)$-linear isomorphism $f:M_1\tilde\to M_2$ such that $\vph_2\circ f=f\circ\vph_1$ and  the isomorphism ${\bf GL}_{M_1}\tilde\to {\bf GL}_{M_2}$ induced by $f$, takes $G_1$ onto $G_2$. 
\medskip
The pair $(M[{1\over p}],\vph)$ is called an $F$-isocrystal over $k$. If we have $pM\subseteq\vph(M)\subseteq M$, then the pair $(M,\vph)$ is called a Dieudonn\'e module over $k$. For $g^\prime\in G(B(k))$ let $g^\prime\vph$ be the $\sg$-linear automorphism of $M[{1\over p}]$ that takes $x\in M[{1\over p}]$ into $g^\prime(\vph(x))\in M[{1\over p}]$. The triple $(M,g^\prime\vph,G)$ is also a latticed $F$-isocrystal with a group over $k$. 
\smallskip
Often there exists a ``good" class ${\bf M}$ of motives over $k$ that has the following property. The crystalline realization of any motive $\Mm$ in ${\bf M}$ is naturally identified with $(M,g_{\Mm}\vph)$ for some $g_{\Mm}\in G(W(k))$ and moreover $G_{B(k)}$ is the identity component of the subgroup of ${\bf GL}_{M[{1\over p}]}$ that fixes some tensors of the tensor algebra of $M[{1\over p}]\oplus {\rm Hom}(M[{1\over p}],B(k))$ which do not depend on $\Mm$ and which are (expected to be) crystalline realizations of motives over $k$ that are intrinsically associated to $\Mm$. For instance, see [40, \S5 and \S6] for contexts that pertain to classes of $H^1$ motives of abelian varieties over ${\rm Spec}(k)$ which are associated to $k$-valued points of a (fixed) good integral model of a Shimura variety of Hodge type. The paper [40] and many previous ones (like [25]) deal with particular cases of such triples $(M,\vph,G)$'s: the pair $(M,\vph)$ is a Dieudonn\'e module over $k$, the group scheme $G$ is reductive, and there exists a semisimple element $s_{\vph}\in G(B(k))$ whose eigenvalues are $1$ and $p$ and such that $\vph s_{\vph}^{-1}$ is a $\sg$-linear automorphism of $M$. Any good classification of the triples $(M,g_{\Mm}\vph,G)$ up to isomorphisms defined by elements of $G(W(k))$, is often an important tool toward the classification of motives in ${\bf M}$. 
\smallskip
Classically, one approaches the classification of all triples $(M,g\vph,G)$ with $g\in G(W(k))$, up to isomorphisms defined by elements of $G(W(k))$, in two steps. The first step aims to classify $(M[{1\over p}],g\vph,G_{B(k)})$'s up to isomorphisms defined by elements of $G(B(k))$. The second steps aims to use the first step in order to study $(M,g\vph,G)$'s. 
\smallskip
A systematic and general approach to the first step was started in [24], which works in the context in which the group $G_{B(k)}$ is reductive, $k=\bar k$, and the pair $(M[{1\over p}],G_{B(k)})$ has a ${\bf Q}_p$ structure $(M_{{\bf Q}_p},G_{{\bf Q}_p})$ with respect to which $\vph$ becomes $g_{\vph}(1_{M_{{\bf Q}_p}}\otimes\sg)$ for some $g_{\vph}\in G(B(k))$; thus, in order to classify $(M[{1\over p}],g\vph,G_{B(k)})$'s up to isomorphisms defined by elements of $G(B(k))$, one only has to describe the image $G_{\vph}$ of the set $\{gg_{\vph}|g\in G(W(k))\}$ in the set $B(G_{{\bf Q}_p})$ of $\sg$-conjugacy classes of elements of $G_{{\bf Q}_p}(B(k))=G(B(k))$. Even if $k=\bar k$, in general such ${\bf Q}_p$ structures do not exist (for instance, they do not exist if the group $G_{B(k)}$ is commutative and $({\rm Lie}(G_{B(k)}),\vph)$ has non-zero slopes). 
\smallskip
One can define two natural equivalence relations $I_{\vph}$ and $R_{\vph}$ on the set underlying the group $G(W(k))$ as follows. A pair $(g_1,g_2)\in G(W(k))^2$ belongs to $I_{\vph}$ (resp. to $R_{\vph}$) if and only if there exists $g_{12}\in G(W(k))$ (resp. $g_{12}\in G(B(k))$) such that $g_{12}g_1\vph=g_2\vph g_{12}$. The set of isomorphism classes of $(M,g\vph,G)$'s (up to isomorphisms defined by elements of $G(W(k))$) is in natural bijection to the quotient set $G(W(k))/I_{\vph}$. The quotient set $G(W(k))/R_{\vph}$ is a more general version of the above type of sets $G_{\vph}$. In general, the natural surjective map $G(W(k))/I_{\vph}\twoheadrightarrow G(W(k))/R_{\vph}$ is not an injection and some of its fibres have the same cardinality as $k$. In general, one can not ``recover" $(M,g\vph,G)$ and its reductions modulo powers of $p$ from the equivalence class $[g]\in G(W(k))/I_{\vph}$ and from the triple $(M[{1\over p}],g\vph,G_{B(k)})$. The last two sentences explain why in this paper, for the study of the quotient set $G(W(k))/I_{\vph}$ and of (reductions modulo powers of $p$ of) $(M,g\vph,G)$'s, we can not appeal to the results of [24], [37], etc. In addition, the language of latticed $F$-isocrystals is more general and more suited for reductions modulo powers of $p$, for endomorphisms, for deformations, and for functorial purposes than the language of either $\sg$-conjugacy classes or equivalence classes of $I_{\vph}$.
\smallskip
If $g_1$, $g_2$, $g_{12}\in G(W(k))$ satisfy $g_{12}g_1\vph=g_2\vph g_{12}$, it is of interest to keep track of the greatest number $n_{12}\in{\bf N}\cup\{0\}$ such that $g_{12}$ and $1_M$ are congruent mod $p^{n_{12}}$. As the relation $I_{\vph}$ is not suitable for this purpose, it will not be used outside this introduction.
\smallskip
The set $\{(M,g\vph,G)|g\in G(W(k))\}$ is in natural bijection to $G(W(k))$. Any set of the form $\{(M,g\vph,G)|g\in G(W(k))\}$ will be called a {\it family} of latticed $F$-isocrystals with a group over $k$. This paper is a starting point for general classifications of families of latticed $F$-isocrystals with a group over $\bar k$. The fact that such classifications are achievable is supported by the following universal principle. 
\medskip\smallskip
{\bf 1.2. Main Theorem A (Crystalline Boundedness Principle).} {\it Suppose $k=\bar k$. Let $(M,\vph,G)$ be a latticed $F$-isocrystal with a group over $k$. Then there exists a number $n_{\rm fam}\in {\bf N}\cup\{0\}$ that is effectively bounded from above and that has the property that for any pair $(g,g_{n_{\rm fam}})\in G(W(k))^2$ such that $g_{n_{\rm fam}}$ is congruent mod $p^{n_{\rm fam}}$ to $1_M$, there exist isomorphisms between $(M,g\vph,G)$ and $(M,g_{n_{\rm fam}}g\vph,G)$ which are elements of $G(W(k))$.}
\medskip
Thus the equivalence class $[g]\in G(W(k))/I_{\vph}$ depends only on $g$ mod $p^{n_{\rm fam}}$; this supports our {\it bounded infinite process} point of view. If $G={\bf GL}_M$ and $(M,\vph)$ is a Dieudonn\'e module over $k$, then Main Theorem A is a direct consequence of [28, p. 44, 3.6, and 3.8]. By a classical theorem of Dieudonn\'e (see [7, Thms. 3 and 5], [28, \S2], [5, Ch. IV, \S4], or [14, Ch. III, \S6]), the category of $p$-divisible groups over ${\rm Spec}(k)$ is antiequivalent to the category of Dieudonn\'e modules over $k$. Thus we get a new proof of the following result which in essence is due to Manin and which is also contained in [36].
\medskip\smallskip
{\bf 1.3. Corollary.} {\it There exists a smallest number $T(r,d)\in{\bf N}\cup\{0\}$ such that any $p$-divisible group $D$ over ${\rm Spec}(\bar k)$ of height $r$ and dimension $d$, is uniquely determined up to isomorphism by its $p^{T(r,d)}$-torsion subgroup scheme $D[p^{T(r,d)}]$. Upper bounds of $T(r,d)$ are effectively computable in terms of $r$.}
\medskip\smallskip
{\bf 1.4. On the proof of Main Theorem A.} 
The proof of Main Theorem A (see 3.1) relies on what we call the {\it stairs method}. The method is rooted on the simple fact that for any $t\in{\bf N}$ and every $y,z\in {\rm End}(M)$, the two automorphisms $1_M+p^ty$ and $1_M+p^tz$ of $M$ commute mod $p^{2t}$. To outline the method, we assume in this paragraph that $G$ is smooth over ${\rm Spec}(W(k))$. Let $m\in{\bf N}\cup\{0\}$ be the smallest number for which there exists a $W(k)$-submodule $E$ of ${\rm Lie}(G)$ that contains $p^m({\rm Lie}(G))$ and that has a $W(k)$-basis $\{e_1,e_2,\ldots,e_v\}$ such that for $l\in\{1,\ldots,v\}$ we have $\vph(e_l)=p^{n_l}e_{\pi(l)}$, where $\pi$ is a permutation of the set $\{1,\ldots,v\}$ and where $n_l$'s are integers  that have the following {\it stairs property}. For any cycle $(l_1,\ldots,l_q)$ of $\pi$, the integers $n_{l_1},\ldots,n_{l_q}$ are either all non-negative or all non-positive. The existence of $m$ is implied by Dieudonn\'e's classification of $F$-isocrystals over $k$ (see [28, \S2]). In general, the $W(k)$-submodule $E$ is not a Lie subalgebra of ${\rm Lie}(G)$. For any $\tilde g\in G(W(k))$ congruent mod $p^{2m+t}$ to $1_M$, there exists $\tilde e\in E$ such that the elements $\tilde g$ and $1_M+p^{m+t}\tilde e$ of ${\bf GL}_M(W(k))$ are congruent mod $p^{2m+1+t}$. Due to this and the stairs property, for $p\ge 3$ there exists an isomorphism between $(M,\tilde g\vph,G)$ and $(M,\vph,G)$ which is an element $\tilde g_0\in G(W(k))$ congruent mod $p^{m+t}$ to $1_M$ (see 3.1.1). If $p=2$, then a slight variant of this holds. Exponential maps (see 2.6) substitute from many points of view the classical Verschiebung maps of Dieudonn\'e modules; for instance, one can choose $\tilde g_0$ to be an infinite product of exponential elements of the form $\sum_{i=0}^{\infty} {{p^{i(m+1)}}\over {i!}}e^i$, where $e\in E$. See 2.2 to 2.4 for the $\sg$-linear preliminaries that are necessary for the estimates which give us the effectiveness part of Main Theorem A. These estimates provide inductively upper bounds of $m$ in terms of $\dim(G_{B(k)})$ and of the $s$-number and the $h$-number of the latticed $F$-isocrystal $({\rm Lie}(G),\vph)$ over $k$ (see 2.2.1 (e) for these two non-negative integers which do not change if $\vph$ is replaced by $g\vph$). 
\medskip\smallskip
{\bf 1.5. Complements, examples, and applications.} See 3.2 for interpretations and variants of Main Theorem A in terms of reductions modulo powers of $p$; in particular, see 3.2.4 for the passage from Main Theorem A to Corollary 1.3. In 3.3 we improve (in many cases of interest) the upper bounds (of $n_{\rm fam}$, etc.) we obtain in 3.1.1 to 3.1.5. 
\smallskip
In \S4 we include four examples. It is well known that if the $p$-divisible group $D$ is ordinary, then $D$ is uniquely determined up to isomorphism by $D[p]$ and moreover $D$ has a unique lift to ${\rm Spec}(W(\bar k))$ (called the canonical lift) that has the property that any endomorphism of $D$ lifts to it. Example 2 identifies the type of latticed $F$-isocrystals with a group over $\bar k$ to which the last two facts generalize naturally  (see 4.3.1 and 4.3.2). Example 4 shows that if $r=2d$, $d\ge 3$, and the slopes of the Newton polygon of $D$ are ${1\over d}$ and ${{d-1}\over d}$, then $D$ is uniquely determined up to isomorphism by $D[p^3]$ (see 4.5). 
\smallskip
In \S5 we list four direct applications of Main Theorem A and of 3.2. First we present the homomorphism form of Main Theorem A (see 5.1.1). Second we define {\it transcendental degrees of definition} for many classes of latticed $F$-isocrystals with a group over $\bar k$ (see 5.2). When the transcendental degrees of definition are $0$, we also define (finite) {\it fields of definition}. In particular, Theorem 5.2.3 (when combined with Lemma 3.2.2) implies that it is possible to build up an {\it atlas} and a {\it list of tables} of isomorphism classes of $p$-divisible groups (endowed with certain extra structures) over ${\rm Spec}(\bar k)$ that are definable over the spectrum of a fixed finite field ${\bf F}_{p^q}$, which are similar in nature to the atlas of finite groups (see [3]) and to the list of tables of elliptic curves over ${\rm Spec}({\bf Q})$ (see [4]). 
\smallskip
Let $N\in{\bf N}\setminus\{1,2\}$ be relatively prime to $p$. Let $\Ma_{d,1,N}$ be the smooth, quasi-projective Mumford moduli scheme over ${\rm Spec}({\bf F}_p)$ that parametrizes isomorphism classes of principally polarized abelian schemes with level-$N$ structure and of relative dimension $d$ over ${\rm Spec}({\bf F}_p)$-schemes (see [33, Thms. 7.9 and 7.10]). Third we apply the principally quasi-polarized version of Corollary 1.3 (see 3.2.5) to get a new type of stratification of $\Ma_{d,1,N}$. Here the word stratification is used in a wide sense (see 2.1.1) which allows the number of strata to be infinite. The strata we get are defined by isomorphism classes of principally quasi-polarized $p$-divisible groups of height $2d$ over spectra of algebraically closed fields of characteristic $p$; they are regular and equidimensional (see 5.3.1 and 5.3.2). Moreover, this new type of stratification of $\Ma_{d,1,N}$ satisfies the {\it purity property} we define in 2.1.1, i.e. its strata are affine $\Ma_{d,1,N}$-schemes (see 5.3.1 and 5.3.2). Variants of 1.3, 3.2.5, 3.2.6, and 5.3.2 but without its purity property part, are also contained in [36]. 
\smallskip
Fourth we get a new proof (see 5.4) of the ``Katz open part" of the Grothendieck--Katz specialization theorem for Newton polygons (see [22, 2.3.1 and 2.3.2]). 
\smallskip
The main goal of \S6 is to prove the following result (see 6.1 and 6.2).
\medskip\smallskip
{\bf 1.6. Main Theorem B.} {\it Let ${\got C}$ be an $F$-crystal in locally free sheaves of finite rank over a reduced ${\rm Spec}({\bf F}_p)$-scheme $S$. Then the Newton polygon stratification of $S$ defined by  ${\got C}$ satisfies the purity property (i.e. each stratum of it is an affine $S$-scheme).}
\medskip
A variant of Main Theorem B was obtained first in [10, 4.1], for the particular case when $S$ is locally noetherian. The fact that the variant is a weaker form of Main Theorem B is explained in 6.3 (a). The main new idea of \S6 is: Newton polygons are encoded in the existence of suitable morphisms between different evaluations of $F$-crystals (viewed without connections) at Witt schemes of (effectively computable) finite lengths. The proof of Main Theorem B combines this new idea with the results of Katz (see [22, 2.6 and 2.7]) on isogenies between $F$-crystals of constant Newton polygons over spectra of (perfections of) complete, discrete valuation rings that are of the form $\bar k[[x]]$.
\medskip\smallskip
{\bf Acknowledgements.} We would like to thank U. of Utah and U. of Arizona for good conditions in which to write this paper and D. Ulmer, G. Faltings, and the referee for valuable comments.
\bigskip\smallskip
\centerline{\bigsll {\bf \S 2 Preliminaries}}
\bigskip\smallskip
See 2.1 for our main notations and conventions. See 2.2 for few definitions and simple properties that pertain to latticed $F$-isocrystals with a group over $k$. In particular, in 2.2.2 we define Dieudonn\'e--Fontaine torsions and volumes of latticed $F$-isocrystals. Inequalities and estimates on such torsions are gathered in 2.3 and 2.4 (respectively); they are essential for examples and for the effectiveness part of 1.2. In 2.5 we apply [42] to get ${\bf Z}_p$ structures for many classes of latticed $F$-isocrystals with a group over $\bar k$. In 2.6 and 2.7 we include group scheme theoretical properties that are needed in \S3 and \S4. In 2.8 we present complements on the categories $\Mm(W_q(S))$ we will introduce in 2.1. In 2.9 we recall two results of commutative algebra. Subsections 2.8 and 2.9 are not used before 5.4. For Newton polygons of $F$-isocrystals over $k$ we refer to [22, 1.3].
\medskip\smallskip
{\bf 2.1. Notations and conventions.} 
By $w$ we denote an arbitrary variable. If $q\in{\bf N}$, let ${\bf F}_{p^q}$ be the field with $p^q$ elements. If $R$ is a commutative ${\bf F}_p$-algebra, let $W(R)$ be the Witt ring of $R$ and let $W_q(R)$ be the ring of Witt vectors of length $q$ with coefficients in $R$. We identify $R=W_1(R)$. Let $\Phi_R$ be the canonical Frobenius endomorphism of either $W(R)$ or $W_q(R)$; we have $\Phi_k=\sg_k=\sg$. Let $R^{(p^q)}$ be $R$ but viewed as an $R$-algebra via the $q$-th power Frobenius endomorphism $\Phi_R^q:R\to R$. If $R$ is reduced, let $R^{\rm perf}:={\rm ind.}\,{\rm lim.}\,_{q\in {\bf N}} R^{(p^q)}$ be the perfection of $R$. 
\smallskip
Let $\Mm(W_q(R))$ be the abelian category whose objects are $W_q(R)$-modules endowed with $\Phi_{R}$-linear endomorphisms and whose morphisms are $W_q(R)$-linear maps that respect the $\Phi_{R}$-linear endomorphisms. We identify $\Mm(W_q(R))$ with a full subcategory of $\Mm(W_{q+1}(R))$ and thus we can define $\Mm(W(R)):=\cup_{q\in{\bf N}}\Mm(W_q(R))$. 
\smallskip
If $S$ is a ${\rm Spec}({\bf F}_p)$-scheme, in a similar way we define $W_q(S)$, $\Phi_S$, $\Mm(W_q(S))$, and $\Mm(W(S))$. We view $W_q(S)$ as a scheme and by a $W_q(S)$-module we mean a quasi-coherent module over the structure ring sheaf $\Mo_{W_q(S)}$ of $W_q(S)$. The formal scheme $W(S)$ is used only as a notation. If $S={\rm Spec}(R)$, then we identify canonically $\Mm(W_q(R))=\Mm(W_q(S))$ and $\Mm(W(R))=\Mm(W(S))$. If $t\in \{1,\ldots,q\}$ and $*(q)$ is a morphism of $\Mm(W_q(S))$, let $*(t)$ be the morphism of $\Mm(W_t(S))$ that is the tensorization of $*(q)$ with $W_t(S)$. Let $S^{\rm top}$ be the topological space underlying $S$. All crystals over $S$ (i.e. all crystals on Berthelot's crystalline site ${\rm CRIS}(S/{\rm Spec}({\bf Z}_p))$) are in locally free sheaves of finite rank. An $F$-crystal ${\got C}$ over $S$ comprises from a crystal ${\got M}$ over $S$ and an isogeny $\Phi_S^*({\got M})\to {\got M}$ of crystals over $S$; let $h_{\got C}\in {\bf N}\cup\{0\}$ be the smallest number such that $p^{h_{\got C}}$ annihilates the cokernel of this isogeny. We identify an $F$-crystal (resp. an $F$-isocrystal) over ${\rm Spec}(k)$ with a latticed $F$-isocrystal $(M,\vph)$ over $k$ that has the property that $\vph(M)\subseteq M$ (resp. with an $F$-isocrystal over $k$ as defined in \S1). The  pulls back of $F$-crystals ${\got C}$ and ${\got C}_*$ over $S$ to an $S$-scheme $S_1$ (resp. to an affine $S$-scheme ${\rm Spec}(R_1)$) are denoted by ${\got C}_{S_1}$ and ${\got C}_{*S_1}$ (resp. by  ${\got C}_{R_1}$ and ${\got C}_{*R_1}$). 
\smallskip
Let $(M,\vph,G)$ be a latticed $F$-isocrystal with a group over $k$. We refer to $M$ as its $W(k)$-module. Let $r_M\in {\bf N}\cup\{0\}$ be the rank of $M$. If $f_1$ and $f_2$ are two ${\bf Z}$-endomorphisms of either $M$ or $M[{1\over p}]$, let $f_1f_2:=f_1\circ f_2$. Two ${\bf Z}$-endomorphisms of $M$ are said to be congruent mod $p^q$ if their reductions mod $p^q$ coincide. Let $M^*:={\rm Hom}(M,W(k))$. Let 
$$\Mt(M):=\oplus_{t,u\in{\bf N}\cup\{0\}} M^{\otimes t}\otimes_{W(k)} M^{*\otimes u}.$$ 
We denote also by $\vph$ the $\sg$-linear automorphism of $\Mt(M)[{1\over p}]$ that takes $f\in M^*[{1\over p}]$ into $\sg f\vph^{-1}\in M^*[{1\over p}]$ and that acts on $\Mt(M)[{1\over p}]$ in the natural tensor product way. The canonical identification ${\rm End}(M[{1\over p}])=M[{1\over p}]\otimes_{B(k)} M^*[{1\over p}]$ is compatible with the $\vph$ actions (see 1.1 (a) for the action of $\vph$ on ${\rm End}(M[{1\over p}])$). If $O$ is either a free $W(k)$-submodule or a $B(k)$-vector subspace of $\Mt(M)[{1\over p}]$ such that $\vph(M)\subseteq M$, then we denote also by $\vph$ the $\sg$-linear endomorphism of $O$ induced by $\vph$. The $W(k)$-span of tensors $v_1,\ldots,v_n\in \Mt(M)[{1\over p}]$ is denoted by $<v_1,\ldots,v_n>$. The latticed $F$-isocrystal $(M^*,\vph)$ over $k$ is called the dual of $(M,\vph)$. We emphasize that the pair $(M^*,\vph)$ involves no Tate twist. A bilinear form on $M$ is called perfect if it defines naturally a $W(k)$-linear isomorphism $M\tilde\to M^*$.
\smallskip
Let $\tilde G_{B(k)}$ be a connected subgroup of ${\bf GL}_{M[{1\over p}]}$. As $\vph$ is a $\sg$-linear automorphism of $M[{1\over p}]$, the group $\{\vph\tilde g\vph^{-1}|\tilde g\in \tilde G_{B(k)}(B(k))\}$ is the group of $B(k)$-valued points of the unique connected subgroup of ${\bf GL}_{M[{1\over p}]}$ that has $\vph({\rm Lie}(\tilde G_{B(k)}))$ as its Lie algebra (see [1, Ch. II, 7.1] for the uniqueness part). So as $\vph$ normalizes ${\rm Lie}(G_{B(k)})$, for $g\in G(B(k))$ we have $\vph g\vph^{-1}\in G(B(k))$; in what follows this fact is used without any extra comment. 
\smallskip
In this paragraph we assume $\vph(M)\subseteq M$. We also refer to $(M,\vph,G)$ as an $F$-crystal with a group over $k$. The Hodge slopes of $(M,\vph)$  (see [22, 1.2]) are the non-negative integers $h_1,\ldots,h_{r_M}$ such that the torsion $W(k)$-module $M/\vph(M)$ is isomorphic to $\oplus_{i=1}^{r_M} W(k)/(p^{h_i})$. If $O$ is a $W(k)$-submodule of $M$ such that $\vph(O)\subseteq O$, we denote also by $\vph$ the $\sg$-linear endomorphism of $M/O$ induced by $\vph$. We refer to the triple $(M/p^qM,\vph,G_{W_q(k)})$ as the reduction mod $p^q$ of $(M,\vph,G)$. If $G={\bf GL}_M$, then often we do not mention $G$ and $G_{W_q(k)}$ and we omit ``with a group". The reduction $(M/p^qM,\vph)$ mod $p^q$ of $(M,\vph)$ is an object of $\Mm(W_q(k))$. 
\smallskip
If $a$, $b\in{\bf Z}$ with $b\ge a$, let $S(a,b):=\{a,a+1,\ldots,b\}$. If $l\in {\bf N}$, if $*$ is a small letter, and if $(*_1,\ldots,*_l)$ is an $l$-tuple which is either an element of ${\bf Z}^l$ or an ordered $W(k)$-basis of some $W(k)$-module, then we define $*_t$ for any $t\in{\bf Z}$ via the rule: $*_t:=*_u$, where $u\in\{1,\ldots,l\}\cap (t+l{\bf Z})$. If $x\in{\bf R}$, let $[x]$ be the greatest integer of the interval $(-\infty,x]$. 
\medskip
{\bf 2.1.1. Conventions on stratifications.} Let $K$ be a field. By a {\it stratification} $\Ms$ of a reduced ${\rm Spec}(K)$-scheme $X$ ({\it in potentially an infinite number of strata}), we mean that:
\medskip
{\bf (i)} for any field $L$ that is either $K$ or an algebraically closed field that contains $K$, a set $\Ms_L$ of disjoint reduced, locally closed subschemes of $X_L$ is given such that each point of $X_L$ with values in an algebraic closure of $L$ factors through some element of $\Ms_L$;
\smallskip
{\bf (ii)} if $i_{12}:L_1\hookrightarrow L_2$ is an embedding between two fields as in (a), then the reduced scheme of the pull back to $L_2$ of any member of $\Ms_{L_1}$, is an element of $\Ms_{L_2}$; so we have a natural pull back injective map $\Ms(i_{12}):\Ms_{L_1}\hookrightarrow\Ms_{L_2}$. 
\medskip
If the inductive limit of all maps $\Ms(i_{12})$ exists (resp. does not exist) in the category of sets, then we say that the stratification $\Ms$ has a {\it class} which is (resp. is not) a set. Each element of some set $\Ms_L$ is referred as a {\it stratum} of $\Ms$. We say $\Ms$ satisfies the {\it purity property} if for any field $L$ as in (a), every element of $\Ms_L$ is an affine $X_L$-scheme.${}^1$ $\vfootnote{1} {This is a more practical, refined, and general definition than any other one that relies on codimension $1$ statements on complements. See Remark 6.3 (a) below.}$ Thus $\Ms$ satisfies the purity property if and only if each stratum of it is an affine $X$-scheme. If all maps $\Ms(i_{12})$'s are bijections, then we identify $\Ms$ with $\Ms_K$ and we say $\Ms$ is {\it of finite type}. 
\medskip\smallskip
{\bf 2.2. Definitions and simple properties.} In this Subsection we introduce few notions and simple properties that pertain naturally to latticed $F$-isocrystals.
\medskip
{\bf 2.2.1. Complements to 1.1.} {\bf (a)} A morphism (resp. an isogeny) between two latticed $F$-isocrystals $(M_1,\vph_1)$ and $(M_2,\vph_2)$ over $k$ is a $W(k)$-linear map (resp. isomorphism) $f:M_1[{1\over p}]\to M_2[{1\over p}]$ such that $f\vph_1=\vph_2 f$ and $f(M_1)\subseteq M_2$. If $f$ is an isogeny, then by its degree we mean $p^l$, where $l$ is the length of the artinian $W(k)$-module $M_2/f(M_1)$. 
\smallskip
{\bf (b)} By a {\it latticed $F$-isocrystal with a group and an emphasized family of tensors} over $k$ we mean a quadruple 
$$(M,\vph,G,(t_{\al})_{\al\in\Mj}),$$ 
where $(M,\vph,G)$ is a latticed $F$-isocrystal with a group over $k$, where $\Mj$ is a set of indices, and where $t_{\al}\in\Mt(M)$ is a tensor that is fixed by both $\vph$ and $G$, such that $G_{B(k)}$ is the subgroup of ${\bf GL}_{M[{1\over p}]}$ that fixes $t_{\al}$ for all $\al\in\Mj$. If $(M_1,\vph_1,G_1,(t_{1\al})_{\al\in\Mj})$ and $(M_2,\vph_2,G_2,(t_{2\al})_{\al\in\Mj})$ are two latticed $F$-isocrystals with a group and an emphasized family of tensors (indexed by the same set $\Mj$) over $k$, by an isomorphism between them we mean an isomorphism $f:(M_1,\vph_1,G_1)\tilde\to (M_2,\vph_2,G_2)$ such that the $W(k)$-linear isomorphism $\Mt(M_1)\tilde\to\Mt(M_2)$ induced by $f$, takes $t_{1\al}$ into $t_{2\al}$ for all $\al\in\Mj$.
\smallskip
{\bf (c)} By a {\it principal bilinear quasi-polarized latticed $F$-isocrystal with a group} over $k$ we mean a quadruple $(M,\vph,G,\lambda_M)$, where $(M,\vph,G)$ is a latticed $F$-isocrystal with a group over $k$ and where $\lambda_M:M\otimes_{W(k)} M\to W(k)$ is a perfect bilinear form with the properties that the $W(k)$-span of $\lambda_M$ is normalized by $G$ and that there exists $c\in{\bf Z}$ such that we have $\lambda_M(\vph(x),\vph(y))=p^c\sg(\lambda_M(x,y))$ for all $x$, $y\in M$. We refer to $\lambda_M$ as a principal bilinear quasi-polarization of $(M,\vph,G)$, $(M,\vph)$, and $(M[{1\over p}],\vph)$. Let $G^0$ be the Zariski closure in ${\bf GL}_M$ of the identity component of the subgroup of $G_{B(k)}$ that fixes $\lambda_M$. We refer to $(M,\vph,G^0)$ as the latticed $F$-isocrystal with a group over $k$ of $(M,\vph,G,\lambda_M)$. The quotient group $G_{B(k)}/G^0_{B(k)}$ is either trivial or isomorphic to ${\bf G}_m$. 
\smallskip
By an isomorphism between two principal bilinear quasi-polarized latticed $F$-iso-crystals with a group $(M_1,\vph_1,G_1,\lambda_{M_1})$ and $(M_2,\vph_2,G_2,\lambda_{M_2})$ over $k$ we mean an isomorphism $f:(M_1,\vph_1,G_1)\tilde\to (M_2,\vph_2,G_2)$ such that we have $\lambda_{M_1}(x,y)=\lambda_{M_2}(f(x),f(y))$ for all $x$, $y\in M_1$. We speak also about principal bilinear quasi-polarized latticed $F$-isocrystals with a group and an emphasized family of tensors over $k$ and about isomorphisms between them; notation $(M,\vph,G,(t_{\al})_{\al\in\Mj},\lambda_M)$. 
\smallskip
If the form $\lambda_M$ is alternating, we drop the word bilinear (i.e. we speak about {\it principal quasi-polarized latticed $F$-isocrystals with a group} over $k$, etc.). 
\smallskip
{\bf (d)} We say the {\it $W$-condition} holds for the latticed $F$-isocrystal with a group $(M,\vph,G)$ over $k$ if there exists a direct sum decomposition $M=\oplus_{i=a}^b \tilde F^i(M)$, where $a$, $b\in{\bf Z}$ with $b\ge a$, such that $M=\oplus_{i=a}^b \vph(p^{-i}\tilde F^i(M))$ and the cocharacter $\mu:{\bf G}_m\to {\bf GL}_M$ defined by the property that $\be\in{\bf G}_m(W(k))$ acts on $\tilde F^i(M)$ through $\mu$ as the multiplication by $\be^{-i}$, factors through $G$. In such a case we also refer to $(M,\vph,G)$ as a {\it $p$-divisible object with a group} over $k$. We refer to the factorization $\mu:{\bf G}_m\to G$ of $\mu$ as a {\it Hodge cocharacter} of $(M,\vph,G)$. For $i\in S(a,b)$ let $F^i(M):=\oplus_{j=b}^i \tilde F^j(M)$. We refer to the decreasing and exhaustive filtration $(F^i(M))_{i\in S(a,b)}$ of $M$ as a {\it lift} of $(M,\vph,G)$. If $G={\bf GL}_M$, we also refer to $(M,\vph)$ as a {\it $p$-divisible object} over $k$. 
\smallskip
Here ``$W$" stands to honor [42, p. 512] while the notion ``$p$-divisible object" is a natural extrapolation of the terminology ``object" introduced in [11, \S2].
\smallskip
{\bf (e)} By the {\it shifting number} (to be abbreviated as the {\it $s$-number}) of a latticed $F$-isocrystal $(M,\vph)$ over $k$ we mean the smallest number $s\in{\bf N}\cup\{0\}$ such that $\vph(p^sM)\subseteq M$ (equivalently such that $\vph(M)\subseteq p^{-s}M$). By the {\it greatest Hodge slope} (to be abbreviated as the {\it $h$-number}) of $(M,\vph)$ we mean the greatest Hodge slope $h$ of $(M,p^s\vph)$, i.e. the unique number $h\in{\bf N}\cup\{0\}$ such that we have $p^{h-s}M\subseteq\vph(M)$ and $p^{h-s-1}M\not\subseteq\vph(M)$. 
\smallskip
We have $s=0$ if and only if $(M,\vph)$ is an $F$-crystal over $k$; in this case $h$ is the number $h_{(M,\vph)}$ defined in 2.1. We have $s=0$ and $h\in\{0,1\}$ if and only if $(M,\vph)$ is a Dieudonn\'e module over $k$. 
\smallskip
Let $s^*$ and $h^*$ be the $s$-number and the $h$-number (respectively) of $(M^*,\vph)$. We have $\vph(M^*)=\vph(M)^*\subseteq p^{s-h}M^*$ but $\vph(M^*)\not\subseteq p^{s-h+1}M^*$. Thus $s^*={\rm max}\{0,h-s\}$. As $(M,\vph)$ is the dual of $(M^*,\vph)$, we also have $s={\rm max}\{0,h^*-s^*\}$. So if $s=0$, then $s^*=h$ and $h^*\in S(0,h)$. If $s>0$, then $s=h^*-s^*$ and thus $h^*=s+s^*={\rm max}(s,h)$. 
\smallskip
If $s=0$, then the $s$-number and the $h$-number of $({\rm End}(M),\vph)=(M,\vph)\otimes (M^*,\vph)$ are at most $s+s^*=h$ and $h+h^*\le 2h$ (respectively).   
\medskip
{\bf 2.2.2. Definitions.} {\bf (a)} Let $(M,\vph)$ be a $p$-divisible object $(M,\vph)$ over $k$. We say $(M,\vph)$ is a {\it cyclic Dieudonn\'e--Fontaine $p$-divisible object} over $k$ if there exists a $W(k)$-basis $\{e_1,\ldots,e_{r_M}\}$ of $M$ such that for $i\in S(1,r_M)$ we have an identity $\vph(e_i)=p^{n_i}e_{i+1}$, where $n_1,\ldots,n_{r_M}$ are integers that are either all non-negative or all non-positive. We refer to $\{e_1,\ldots,e_{r_M}\}$ as a {\it standard $W(k)$-basis} of $(M,\vph)$. 
\smallskip
We say $(M,\vph)$ is an {\it elementary Dieudonn\'e--Fontaine $p$-divisible object} over $k$ if it is a cyclic Dieudonn\'e--Fontaine $p$-divisible object over $k$ that is not the direct sum of two or more non-trivial cyclic Dieudonn\'e--Fontaine $p$-divisible objects over $k$. 
\smallskip
We say $(M,\vph)$ is an {\it elementary Dieudonn\'e $p$-divisible object} over $k$ if there exists a $W(k)$-basis $\{e_1,\ldots,e_{r_M}\}$ of $M$ such that for $i\in S(2,r_M)$ we have an identity $\vph(e_i)=e_{i+1}$ and moreover $\vph(e_1)=p^{n_1}e_2$ for some integer $n_1$ that is relatively prime to $r_M$.
\smallskip
We say $(M,\vph)$ is a {\it Dieudonn\'e--Fontaine (resp. a Dieudonn\'e) $p$-divisible object} over $k$ if it is a direct sum of elementary Dieudonn\'e--Fontaine (resp. of elementary Dieudonn\'e) $p$-divisible objects over $k$. 
\smallskip
{\bf (b)} By the {\it Dieudonn\'e--Fontaine torsion} (resp. {\it volume)} of a latticed $F$-isocrystal $(M,\vph)$ over $k$ we mean the smallest number 
$${\bf T}(M,\vph)\in{\bf N}\cup\{0\}$$ 
(resp. ${\bf V}(M,\vph)\in{\bf N}\cup\{0\}$) such that there exists a Dieudonn\'e--Fontaine  $p$-divisible object $(M_1,\vph_1)$ over $\bar k$ for which we have an isogeny $f:(M_1,\vph_1)\hookrightarrow (M\otimes_{W(k)} W(\bar k),\vph\otimes\sg_{\bar k})$ with the property that $p^{{\bf T}(M,\vph)}M\subseteq f(M_1)$ (resp. that $M/f(M_1)$ has length ${\bf V}(M,\vph)$). By replacing Dieudonn\'e--Fontaine with Dieudonn\'e, in a similar way we define the {\it Dieudonn\'e torsion} ${\bf T}_+(M,\vph)\in{\bf N}\cup\{0\}$ and the {\it Dieudonn\'e volume} ${\bf V}_+(M,\vph)\in{\bf N}\cup\{0\}$ of $(M,\vph)$.
\medskip 
{\bf 2.2.2.1. Remarks.} {\bf (a)} Any (elementary) Dieudonn\'e $p$-divisible object over $k$ is also an (elementary) Dieudonn\'e--Fontaine $p$-divisible object over $k$. Moreover, any Dieudonn\'e--Fontaine $p$-divisible object over $k$ is definable over ${\bf F}_p$. 
\smallskip
{\bf (b)} The existence of ${\bf V}_+(M,\vph)$ (and thus also of ${\bf V}(M,\vph)$, ${\bf T}_+(M,\vph)$, and ${\bf T}(M,\vph)$) is equivalent to Dieudonn\'e's classification of $F$-isocrystals over $\bar k$. This and the fact that suitable reductions (modulo powers of $p$) of $p$-divisible objects over $k$ are studied systematically for the first time in [14] and [15], explains our terminology. 
\smallskip
{\bf (c)} Classically one works only with Dieudonn\'e $p$-divisible objects (as they are uniquely determined by their Newton polygons) and with Dieudonn\'e volumes (as they keep track of degrees of isogenies); see [7], [28], [6], [10], etc. But working with Dieudonn\'e--Fontaine $p$-divisible objects and torsions one can get considerable improvements for many practical calculations or upper bounds (like the ones we will encounter in \S3). 
\medskip
{\bf 2.2.3. Lemma.} {\it Let $K$ be an algebraically closed field that contains $k$. Let $(M,\vph)$ be a Dieudonn\'e--Fontaine  $p$-divisible object over $k$ with the property that $\vph(M)\subseteq M$. Let $h$ be the $h$-number of $(M,\vph)$, let $e_M:=\max\{r_M,[{r_M^2\over 4}]\}$, let $k_1$ be the composite field of $k$ and ${\bf F}_{p^{r_M!}}$, and let $m\in{\bf N}$. We have the following two properties:
\medskip
{\bf (a)} For any endomorphism $f_{he_M+m}$ of $(M\otimes_{W(k)} W_{he_M+m}(K),\vph\otimes\sg_{K})$, the reduction $f_{m}$ mod $p^m$  of $f_{he_M+m}$ is the scalar extension of an endomorphism of $(M\otimes_{W(k)} W_m(k_1),\vph\otimes\sg_{k_1})$. If $(M,\vph)$ is a Dieudonn\'e $p$-divisible object over $k$, then the previous sentence holds with $e_M$ being substituted by $r_M$. 
\smallskip
{\bf (b)} Each endomorphism of $(M\otimes_{W(k)} W(K),\vph\otimes\sg_{K})$ is the scalar extension of an endomorphism of $(M\otimes_{W(k)} W(k_1),\vph\otimes\sg_{k_1})$.}
\medskip
{\it Proof:} We write $(M,\vph)=\oplus_{i=1}^s (M_i,\vph)$ as a direct sum of elementary Dieudonn\'e--Fontaine  $p$-divisible objects over $k$. Let $\{e_1^{(i)},\ldots,e_{r_{M_i}}^{(i)}\}$ be a standard $W(k)$-basis of $(M_i,\vph)$. We check that (a) holds. Let $i_0\in S(1,s)$ and let $j_0\in S(1,r_{M_{i_0}})$. We write $f_{he_M+m}(e_{j_0}^{(i_0)}\otimes 1)=\sum_{i=1}^s\sum_{j=1}^{r_{M_i}} e_j^{(i)}\otimes \be_{j_0j}^{(i_0i)}$, where all $\be_{j_0j}^{(i_0i)}$'s belong to $W_{he_M+m}(K)$. Let $r_{M_{i_0i}}:=l.c.m.\{r_{M_{i_0}},r_{M_i}\}$; it is a divisor of $r_M!$. If $i=i_0$, then $r_{M_{i_0i}}=r_{M_{i_0}}\le r_M\le e_M$. If $i\neq i_0$, then $r_{M_{i_0}}+r_{M_i}\le r_M$ and thus we have $r_{M_{i_0i}}\le r_{M_{i_0}}r_{M_i}\le [{{r_M^2}\over 4}]\le e_M$. 
\smallskip
As $f_{he_M+m}(\vph^{r_{M_{i_0i}}}(e_{j_0}^{(i_0)})\otimes 1)=(\vph\otimes\sg_K)^{r_{M_{i_0i}}}(f_{he_M+m}(e_{j_0}^{(i_0)}\otimes 1))$, we have an equality
$$p^{m_{j_0}^{(i_0i)}}\be_{j_0j}^{(i_0i)}=p^{q_{j}^{(i_0i)}}\sg_{K}^{r_{M_{i_0i}}}(\be_{j_0j}^{(i_0i)})\in W_{he_M+m}(K),\leqno (1)$$ 
where $m_{j_0}^{(i_0i)}\in{\bf N}\cup\{0\}$ is such that $\vph^{r_{M_{i_0i}}}(e_{j_0}^{(i_0)})=p^{m_{j_0}^{(i_0i)}}e_{j_0}^{(i_0)}$ and where $q_{j}^{(i_0i)}\in{\bf N}\cup\{0\}$ is such that $\vph^{r_{M_{i_0i}}}(e_j^{(i)})=p^{q_{j}^{(i_0i)}}e_j^{(i)}$. The numbers $m_{j_0}^{(i_0i)}$ and $q_j^{(i_0i)}$ are at most $hr_{M_{i_0i}}$ and so at most $he_M$. Let $s_{j_0j}^{(i_0i)}\in S(0,he_M+m)$ be the unique number such that we can write $\be_{j_0j}^{(i_0i)}=p^{s_{j_0j}^{(i_0i)}}\tilde\be_{j_0j}^{(i_0i)}$, with $\tilde\be_{j_0j}^{(i_0i)}\in{\bf G}_m(W_{he_M+m}(K))$. From (1) we easily get that 
$$t_{j_0j}^{(i_0i)}:=\min\{he_M+m,m_{j_0}^{(i_0i)}+s_{j_0j}^{(i_0i)}\}\;\;{\rm equals}\;\;{\rm to}\;\;\min\{he_M+m,q_{j}^{(i_0i)}+s_{j_0j}^{(i_0i)}\}$$ 
and that $\sg_K^{r_{M_{i_0i}}}(\tilde\be_{j_0j}^{(i_0i)})$ and $\tilde\be_{j_0j}^{(i_0i)}$ coincide mod $p^{he_M+m-t_{j_0j}^{(i_0i)}}$. Thus $\tilde\be_{j_0j}^{(i_0i)}$ mod $p^{he_M+m-t_{j_0j}^{(i_0i)}}$ belongs to $W_{he_M+m-t_{j_0j}^{(i_0i)}}({\bf F}_{p^{r_{M_{i_0i}}}})$ and therefore $\be_{j_0j}^{(i_0i)}$ mod $p^{he_M+m-t_{j_0j}^{(i_0i)}+s_{j_0j}^{(i_0i)}}$ belongs to $W_{he_M+m-t_{j_0j}^{(i_0i)}+s_{j_0j}^{(i_0i)}}({\bf F}_{p^{r_{M_{i_0i}}}})$. As $-m_{j_0}^{(i_0i)}\le -t_{j_0j}^{(i_0i)}+s_{j_0j}^{(i_0i)}$, we get that $\be_{j_0j}^{(i_0i)}$ mod $p^{he_M+m-m_{j_0}^{(i_0i)}}$ belongs to $W_{he_M+m-m_{j_0}^{(i_0i)}}({\bf F}_{p^{r_{M_{i_0i}}}})$ and thus also to $W_{he_M+m-m_{j_0}^{(i_0i)}}(k_1)$. So due to the inequality $m\le he_M+m-m_{j_0}^{(i_0i)}$, we have $f_m(e_{j_0}^{(i_0)}\otimes 1)\in M\otimes_{W(k)} W_m(k_1)$ for any pair $(i_0,j_0)\in S(1,s)\times S(1,r_{M_{i_0}})$. Thus $f_m$ is the scalar extension of an endomorphism of $(M\otimes_{W(k)} W_m(k_1),\vph\otimes\sg_{k_1})$. 
\smallskip
If $(M,\vph)$ is a Dieudonn\'e $p$-divisible object over $k$, then the Hodge slopes of $(M_{i_0},\vph)$ are $0,\ldots,0$, and some integer in $S(0,h)$; thus $m_{j_0}^{(i_0i)}\le h{r_{M_{i_0i}}\over r_{M_{i_0}}}\le hr_{M_i}\le hr_M$. A similar argument shows that $q_j^{(i_0i)}\le hr_M$. Thus in the previous paragraph we can substitute $e_M$ by $r_M$. So (a) holds. 
\smallskip
Part (b) follows from (a) by taking $m\to\infty$. \endproof 
\medskip
{\bf 2.2.4. Deviations of tuples.} Let $l\in{\bf N}$. Let $\tau=(n_1,\ldots,n_l)\in{\bf Z}^l$. 
\medskip
{\bf (a)} Suppose $\sum_{i=1}^l n_i$ is non-negative (resp. is non-positive). Let $P(\tau)$ be the set of pairs $(t,u)$, where $t\in S(1,l)$ and $u\in S(t,l+t-1)$ have the property that all sums $\sum_{i=v}^{u} n_i$ with $v\in S(t,u)$ are non-positive (resp. are non-negative). By the non-negative (resp. the non-positive) sign deviation of $\tau$ we mean the non-negative integer $\max\{0,-\sum_{i=t}^{u} n_i|(t,u)\in P(\tau)\}$ (resp. $\max\{0,\sum_{i=t}^{u} n_i|(t,u)\in P(\tau)\}$). 
\smallskip
{\bf (b)} If $\sum_{i=1}^l n_i$ is non-negative (resp. is non-positive), then by the non-negative (resp. the non-positive) value deviation of $\tau$ we mean the absolute value of the sum of all non-positive (resp. of all non-negative) entries of $\tau$. As a convention, this sum is $0$ if $\tau$ has no non-positive (resp. no non-negative) entries.
\smallskip
{\bf (c)} If $\sum_{i=1}^l n_i$ is positive (resp. is negative), then by the sign deviation ${\bf S}\tau$ of $\tau$ we mean its non-negative (resp. its non-positive) sign deviation. If $\sum_{i=1}^l n_i=0$, then by the sign deviation ${\bf S}\tau$ of $\tau$ we mean the smaller of its non-negative and non-positive sign deviations. We also use this definition with (sign, ${\bf S}$) replaced by (value, ${\bf W}$). 
\medskip
Samples: ${\bf S}(-1,1,-1,-1,1,1,0,-1)=1+1=2$, ${\bf W}(-1,1,-1,-1,1,1,0,-1)=3$, ${\bf S}(1,1,-2,1,3)={\bf W}(1,1,-2,1,3)=2$, and ${\bf S}(-1,1,-1)={\bf W}(-1,1,-1)=1$.
\medskip\smallskip
{\bf 2.3. Inequalities.} Let $(M,\vph)$ be a latticed $F$-isocrystal over $k$. Obviously ${\bf V}(M,\vph)\le {\bf V}_+(M,\vph)$
and ${\bf T}(M,\vph)\le {\bf T}_+(M,\vph)$. Moreover we have
$${\bf T}(M,\vph)\le {\bf V}(M,\vph)\le {\bf T}(M,\vph)r_M$$ 
and the same inequalities hold with $({\bf T},{\bf V})$ being replaced with $({\bf T}_+,{\bf V}_+)$. 
\medskip
{\bf 2.3.1. Lemma.} {\it Let $\tau=(n_1,\ldots,n_{r_M})\in{\bf Z}^{r_M}$. Suppose there exists  a $W(k)$-basis $\{e_1,\ldots,e_{r_M}\}$ of $M$ such that for $i\in S(1,r_M)$ we have $\vph(e_i)=p^{n_i}e_{i+1}$. Then we have the following sequence of three inequalities
$${\bf T}(M,\vph)\le {\bf S}\tau\le {\bf W}\tau\le\abs{n_1}+\abs{n_2}+\cdots+\abs{n_{r_M}}.\leqno (2)$$}
{\it Proof:}
The second and the third inequalities follow from very definitions.
\smallskip
We check the first inequality of (2) only in the case when $\sum_{i=1}^{r_M} n_i>0$ and at least one entry $n_i$ is negative, as in all other cases the first inequality of (2) is checked in the same way. We perform the following type of operation. 
\smallskip
Let $u\in{\bf N}\cup\{0\}$ be the greatest number such that there exists $t\in S(1,r_M)$ with the property that $n_{t-v,t}:=\sum_{i=t-v}^t n_i$ is non-positive for all $v\in S(0,u)$; we have $n_{t+1}>0$, $n_{t-u-1}>0$, and $\sum_{i=t-u-1}^t n_i>0$. For $v\in\ S(0,u)$ we replace $e_{t-v}$ by $\tilde e_{t-v}:=p^{-n_{t-v,t}}e_{t-v}$.
Up to a cyclic rearrangement of $\tau$, we can assume $t-u=1$; so $t=1+u$ and $n_{r_M}=n_0=n_{t-u-1}>0$. The $r_M$-tuple $(\tilde e_1,\ldots,\tilde e_{u+1},e_{u+2},\ldots,e_{r_M})$ is mapped by $\vph$ into the $r_M$-tuple $(\tilde e_2,\ldots,\tilde e_{u+1},e_{u+2},p^{n_{u+2}}e_{u+3},\ldots,p^{n_{r_M-1}}e_{r_M},p^{\sum_{i=0}^{u+1} n_i}\tilde e_1)$. We have $\sum_{i=0}^{u+1} n_i=\sum_{i=t-u-1}^t n_i> 0$. So if for all $i\in S(u+2,r_M-1)$ we have $n_i\ge 0$, then the pair $(<\tilde e_1,\ldots,\tilde e_{u+1},e_{u+2},\ldots,e_{r_M}>,\vph)$ is a cyclic Dieudonn\'e--Fontaine $p$-divisible object over $k$ and we are done as by very definitions we have $-n_{t-v,t}\in S(0,{\bf S}\tau)$ and thus $p^{{\bf S}\tau}$ annihilates the quotient $W(k)$-module $M/<\tilde e_1,\ldots,\tilde e_{u+1},e_{u+2},\ldots,e_{r_M}>$; if this is not the case, we next deal with the inoperated entries $n_{u+2},\ldots,n_{r_M}$.  
\smallskip
We repeat the operation as follows. Let $u_1\in{\bf N}\cup\{0\}$ be the greatest number such that there exists $t_1\in S(u+2,r_M)$ with the property that $n_{t_1-v_1,t_1}:=\sum_{i=t_1-v_1}^{t_1} n_i$ is non-positive for all $v_1\in S(0,u_1)$; we have $n_{t_1+1}>0$, $n_{t_1-u_1-1}>0$, and $\sum_{i=t_1-u_1-1}^{t_1} n_i>0$. Due to the ``greatest" property of $u$ we have $t_1-u_1> u+2$. For $v_1\in\ S(0,u_1)$ we replace $e_{t_1-v_1}$ by $\tilde e_{t_1-v_1}:=p^{-n_{t_1-v_1,t_1}}e_{t_1-v_1}$ and we repeat the operation for the inoperated entries $n_{u+2}$, $n_{u+3},\ldots,n_{t_1-u_1-1}$, $n_{t_1+1}$, $n_{t_1+2},\ldots,n_{r_M}$. By induction on the number of remaining inoperated entries (they do not have to be indexed by a set of consecutive numbers in $S(1,r_M)$), we get that the first inequality of (2) holds. \endproof
\medskip
{\bf 2.3.2. Example.} If for $i\in S(1,r_M)$ we have $n_i\in\{-1,0,1\}$, then from (2) we get 
$${\bf T}(M,\vph)\le {\bf W}\tau=\min\{n^-,n^+\},\leqno (3)$$
where $n^-$ (resp. $n^+$) is the number of $i$'s such that $n_i=-1$ (resp. such that $n_i=1$). 
\smallskip
We now consider the case when $r_M\ge 3$ and $(n_1,n_2,\ldots,n_{r_M})=(1,1,\ldots,1,-1)$. So $(M,\vph)$ has a unique slope ${{r_M-2}\over {r_M}}$ that is positive. As $\vph(M)\not\subseteq M$, we have ${\bf T}(M,\vph)\ge 1$. But ${\bf T}(M,\vph)\le 1$, cf. (3). Thus ${\bf T}(M,\vph)=1$. In fact $(<e_1,\ldots,pe_{r_M}>,\vph)$ is an elementary Dieudonn\'e--Fontaine $p$-divisible object  over $k$ whose Hodge slopes are $0$, $0$, $1,\ldots,1$. It is easy to see that ${\bf T}_+(M,\vph)=r_M-2$; so ${\bf T}_+(M,\vph)> {\bf T}(M,\vph)$ for $r_M>3$. 
\medskip\smallskip
{\bf 2.4. Estimates.} Let $(M,\vph)$ be a latticed $F$-isocrystal over $k$. Let $s$ and $h$ be the $s$-number and the $h$-number (respectively) of $(M,\vph)$. Let $\Mh$ be the set of slopes of $(M[{1\over p}],\vph)$. If $\al\in\Mh$, we write $\al={a_{\al}\over b_{\al}},$ 
where $(a_{\al},b_{\al})\in{\bf Z}\times{\bf N}$ with $g.c.d.(a_{\al},b_{\al})=1$. 
\medskip
{\bf 2.4.1. Lemma.} {\it Suppose $k=\bar k$. Let $(a,b,c)\in{\bf N}\times ({\bf N}\cup\{0\})\times ({\bf N}\cup\{0\})$. There exists a smallest number 
$$d(a,b,c)\in{\bf N}\cup\{0\}$$ 
(resp. $d_+(a,b,c)\in{\bf N}\cup\{0\}$) such that for any latticed  $F$-isocrystal ${\got C}$ over $k$ of rank $a$, $s$-number $b$, and $h$-number $c$, we have an inequality ${\bf T}({\got C})\le d(a,b,c)$ (resp. ${\bf T}_+({\got C})\le d_+(a,b,c)$). In particular, for any element $g\in$ ${\bf GL}_M(W(k))$ we have ${\bf T}(M,g\vph)\le d(r_M,s,h)$. Moreover upper bounds of $d(a,b,c)$ (resp. of $d_+(a,b,c)$) are effectively computable in terms of $a$, $b$, and $c$.}
\medskip
{\it Proof:} As ${\bf T}({\got C})\le {\bf T}_+({\got C})$, it suffices to prove the Lemma for $d_+(a,b,c)$. To ease the notations, we will assume that $(r_M,s,h)=(a,b,c)$ and that $(M,\vph)={\got C}$. We have
$${\bf T}_+(M,\vph)\le s(\max\{b_{\al}|\al\in \Mh\}-1)+{\bf T}_+(M,p^s\vph).\leqno (4)$$
To check this inequality we first remark that if $O$ is a $W(k)$-submodule of $M$ such that the pair $(O,p^s\vph)$ is an elementary Dieudonn\'e $p$-divisible object over $k$ and if $\{e_1,\ldots,e_{r_O}\}$ is a standard $W(k)$-basis of $(O,p^s\vph)$ such that we have $(p^s\vph)(e_i)=e_{i+1}$ for all $i\in S(1,r_O-1)$, then the following pair $(O^\prime,\vph):=(<p^{r_Os-s}e_1,p^{r_Os-2s}e_2,\ldots,p^se_{r_O-1},e_{r_O}>,\vph)$ is an elementary Dieudonn\'e $p$-divisible object over $k$. As $(O[{1\over p}],\vph)$ is a simple $F$-isocrystal over $k$ whose unique slope belongs to $\Mh$, we have  $r_O\le\max\{b_{\al}|\al\in \Mh\}$. From this and the fact that $O/O^\prime$ is annihilated by $p^{s(r_O-1)}$, we easily get that (4) holds. 
\smallskip
Thus it suffices to prove the existence of a number $d_+(a,b,c)$ that has all the required properties under the extra assumption $b=s=0$; as $s=0$, we have $\Mh\subseteq [0,h]$. We will use an induction on $a=r_M$. The case $a=1$ is trivial. To accomplish for $a\ge 2$ the inductive passage from $a-1=r_M-1$ to $a=r_M$, we consider two disjoint Cases. 
\medskip
{\bf Case 1.} Suppose the $F$-isocrystal $(M[{1\over p}],\vph)$ over $k$ is not simple. Let $\al\in \Mh$. We consider a short exact sequence $0\to (M_1,\vph)\to (M,\vph)\to (M_2,\vph)\to 0$ such that the $F$-isocrystal $(M_2[{1\over p}],\vph)$ over $k$ is simple of $\al$. For $i\in S(1,2)$, the $h$-number of $(M_i,\vph)$ is at most $h$ and we have $r_{M_i}<a=r_M$. By induction, there exists $d_i\in {\bf N}\cup\{0\}$ that has upper bounds effectively computable in terms of $r_{M_i}$ and $c=h$ and such that there exists a $W(k)$-submodule $O_i$ of $M_i$ with the properties that $r_{O_i}=r_{M_i}$, that $p^{d_i}M_i\subseteq O_i$, and that $(O_i,\vph)$ is a Dieudonn\'e $p$-divisible object over $k$. The map $\sg^{b_{\al}}-1_{W(k)}:W(k)\to W(k)$ is onto. This implies that $p^{a_{\al}}M_1\subseteq(\vph^{b_{\al}}-p^{a_{\al}}1_{M_1})(M_1)$. Let $x\in O_2$ be such that $\vph^{b_{\al}}(x)=p^{a_{\al}}x$ and $\vph^{b_{\al}-1}(x)\in O_2\setminus pO_2$. If $\tilde x\in M$ maps into $x$, then there exists $y\in M_1$ such that $\vph^{b_{\al}}(y)-p^{a_{\al}}(y)$ is $p^{a_{\al}}[\vph^{b_{\al}}(\tilde x)-p^{a_{\al}}(\tilde x)]\in p^{a_{\al}}M_1$. Thus $z:=-y+p^{a_{\al}}\tilde x\in M$ maps into $p^{a_{\al}}x$ and we have $\vph^{b_{\al}}(z)=p^{a_{\al}}(z)$. By choosing $x$ to belong to a standard $W(k)$-basis of $(O_2,\vph)$, we get that the monomorphism $i_2:(p^{a_{\al}}O_2,\vph)\hookrightarrow (M_2,\vph)$ lifts to a monomorphism $j_2:(p^{a_{\al}}O_2,\vph)\hookrightarrow (M,\vph)$. As $(O_1+j_2(p^{a_{\al}}O_2),\vph)$ is a Dieudonn\'e $p$-divisible object over $k$ and as $p^{d_1+d_2+a_{\al}}$ annihilates $M/O_1+j_2(p^{a_{\al}}O_2)$, we get
$${\bf T}_+(M,\vph)\le d_1+d_2+a_{\al}\le d_1+d_2+hb_{\al}\le d_1+d_2+hr_M.\leqno (5)$$ 
\indent 
{\bf Case 2.} Suppose the $F$-isocrystal $(M[{1\over p}],\vph)$ over $k$ is simple of slope $\al$. Thus $a=r_M=b_{\al}$. Let $\{e_1,\ldots,e_a\}\subseteq M$ be a $B(k)$-basis  of $M[{1\over p}]$ such that $e_1\in M\setminus pM$, for $i\in S(1,a-1)$ we have $\vph(e_i)=e_{i+1}$, and $\vph(e_a)=p^{a_{\al}}e_1$. For $t\in S(1,a)$, let $M_t:=<e_1,\ldots,e_t>$ and $\tilde M_t:=M_t[{1\over p}]\cap M$. We have $M_1=\tilde M_1$ and $\tilde M_a=M$.
\medskip
{\bf 2.4.1.1. Claim.} {\it There exists a strictly increasing sequence $(c_t)_{t\in S(1,a)}$ of non-negative integers that depends only on $a=r_M$ and $c=h$, that is effectively computable, and that has the property that for any $t\in S(1,a)$ we have inclusions $$p^{c_t}(\tilde M_t)\subseteq M_t\subseteq\tilde M_t.\leqno (6)$$ }
\indent
To check this Claim we use induction on $t\in S(1,a)$. Taking $c_1:=0$, (6) holds for $t=1$. Suppose there exists a number $r\in S(1,a-1)$ such that (6) holds for $t\in S(1,r)$. We now check that (6) holds for $t=r+1$. Thus we have to show that there exists an effectively computable natural number $c_{r+1}$ which is at least $c_r+1$ and for which we have
$$e_{r+1}=\vph(e_r)\notin p^{1+c_{r+1}}M+\tilde M_r.\leqno (7)$$
\indent
We write $e_{r+1}=\vph(e_r)=p^{n_r}x_r+y_r$, where $x_r\in M\setminus pM$, $n_r\in{\bf N}\cup\{0\}$, and $y_r\in\tilde M_r$. By our initial induction (on ranks), we can speak about an effectively computable number $d_r\in{\bf N} $ that is at least ${\rm max}\{d_+(r,0,l)|l\in S(0,h)\}$. Let 
$$
c_{r+1}:=c_r+d_r+r!ah. \leqno (8)
$$
\indent 
We show that the assumption $n_r>c_{r+1}$ leads to a contradiction. Let $\tilde M_0:=0$. Let $\tilde e_r\in \tilde M_r$ be such that we have a direct sum decomposition $\tilde M_r=\tilde M_{r-1}\oplus <\tilde e_r>$. Based on (6) (applied with $t=r$), we can write $e_r=y_{r-1}+l_r\tilde e_r$, where $y_{r-1}\in\tilde M_{r-1}$ and $l_r\in W(k)\setminus p^{c_r+1}W(k)$. Let $w_r:=\sg(l_r)^{-1}(y_r-\vph(y_{r-1}))\in\tilde M_r[{1\over p}]$. We have $\vph(\tilde e_r)=w_r+\sg(l_r)^{-1}p^{n_r}x_r\in M$. As $n_r>c_{r+1}>c_r$, we have $\sg(l_r)^{-1}p^{n_r}x_r\in p^{n_r-c_r}M\subseteq M$; thus $w_r\in \tilde M_r=\tilde M_r[{1\over p}]\cap M$. Let $\eta_r$ be the $\sg$-linear endomorphism of $\tilde M_r$ that acts on $\tilde M_{r-1}$ as $\vph$ does and that takes $\tilde e_r$ into $w_r$; thus $\eta_r(e_r)=y_r$. The difference $\vph(\tilde e_r)-\eta_r(\tilde e_r)$ is $\sg(l_r)^{-1}p^{n_r}x_r\in p^{n_r-c_r}M$. Thus $\vph$ restricted to $\tilde M_r$ and $\eta_r$, when viewed as maps from $\tilde M_r$ to $M$, coincide mod $p^{n_r-c_r}$. From this and the inequality $n_r-c_r\ge h+1$, we get that the pair $(\tilde M_r,\eta_r)$ is an $F$-crystal over $k$ whose $h$-number is at most $h$. 
\smallskip
Let $O_r$ be a $W(k)$-submodule of $\tilde M_r$ such that $p^{{\bf T}_+(\tilde M_r,\eta_r)}\tilde M_r\subseteq O_r$ and $(O_r,\eta_r)$ is a Dieudonn\'e $p$-divisible object over $k$. Let $t_r\in S(0,{\bf T}_+(\tilde M_r,\eta_r))$ be the smallest number such that $p^{t_r}e_1\in O_r\setminus pO_r$. We consider a direct sum decomposition $(O_r,\eta_r)=(O_{r,1},\eta_r)\oplus\cdots\oplus (O_{r,s_r},\eta_r)$ into elementary Dieudonn\'e $p$-divisible objects over $k$; we can assume that the indices are such that the component $e_{1,1}\in O_{r,1}$ of $p^{t_r}e_1\in O_r$ with respect to this direct sum decomposition of $O_r$, is not divisible inside $O_{r,1}$ by $p$ (i.e. we have $e_{1,1}\notin pO_{r,1}$). Let $u_r\in [0,h]$ be the unique slope of $(O_{r,1},\eta_r)$. 
\smallskip
The element $p^{r!a\al}e_1=\vph^{r!a}(e_1)$  is congruent mod $p^{n_r-c_r}$ to $\eta_r^{r!a}(e_1)$. As we have $p^{{\bf T}_+(\tilde M_r,\eta_r)}\tilde M_r\subseteq O_r$, we get that $\eta_r^{r!a}(e_{1,1})-p^{r!a\al}e_{1,1}\in p^{t_r+n_r-c_r-{\bf T}_+(\tilde M_r,\eta_r)}O_{r,1}$; thus $\eta_r^{r!a}(e_{1,1})-p^{r!a\al}e_{1,1}\in p^{n_r-c_r-{\bf T}_+(\tilde M_r,\eta_r)}O_{r,1}$. As $(O_{r,1},\eta_r)$ is an elementary Dieudonn\'e $p$-divisible object over $k$ whose rank divides $r!a$ and as $e_{1,1}\in O_{r,1}\setminus pO_{r,1}$, there exists $z_r\in O_{r,1}\setminus pO_{r,1}$ such that $\eta_r^{r!a}(e_{1,1})=p^{r!au_r}z_r$. Thus $p^{r!au_r}z_r-p^{r!a\al}e_{1,1}\in p^{n_r-c_r-{\bf T}_+(\tilde M_r,\eta_r)}O_{r,1}$. As $d_r\ge {\bf T}_+(\tilde M_r,\eta_r)$ and $n_r>c_{r+1}$, from (8) we get $n_r-c_r-{\bf T}_+(\tilde M_r,\eta_r)\ge r!ha+1$. As $h\ge {\rm max}\{u_r,\al\}$, we have
$r!ah+1>\max\{r!au_r,r!a\al\}$. Thus 
$$n_r-c_r-{\bf T}_+(\tilde M_r,\eta_r)>\max\{r!au_r,r!a\al\}.\leqno (9)$$ 
\indent
From (9) and the relations $p^{r!au_r}z_r-p^{r!a\al}e_{1,1}\in p^{n_r-c_r-{\bf T}_+(M_r,\eta_r)}O_{r,1}$ and $z_r,e_{1,1}\in O_{r,1}\setminus pO_{r,1}$, we get that $r!au_r=r!a\al$. Thus $\al=u_r$ and therefore $a=r_{O_{r,1}}\le r$. This contradicts the fact that $r\in S(1,a-1)$. Thus $n_r\le c_{r+1}$ and so (7) holds. Thus (6) holds for $t=r+1$. As $c_{r+1}$ depends only on $r$, $a$, $h$, $c_r$, and $d_r$ and as (by induction) $c_r$ and $d_r$ depend only on $r$, $a$, and $h$, we get that $c_{r+1}$ depends only on $a$ and $c=h$. This ends our second induction on $t\in S(1,a)$. Thus the Claim holds. 
\smallskip
We have $p^{c_a}M=p^{c_a}\tilde M_a\subseteq M_a\subseteq M$, cf. (6). As $(M_a,\vph)$ is an elementary Dieudonn\'e $p$-divisible object over $k$, we get ${\bf T}_+(M,\vph)\le c_a$. This ends Case 2.
\medskip
The above two Cases imply that ${\bf T}_+(M,\vph)$ has upper bounds that are effectively computable in terms of $a=r_M$ and $c=h$. Thus the number $d_+(a,0,c)$ exists and has effectively computable upper bounds in terms of $a$ and $c$. This ends the initial induction (on ranks $a$) and so it also ends the proof of Lemma 2.4.1.\endproof 
\medskip
{\bf 2.4.2. An interpretation.} The estimates of the proof of 2.4.1.1 are different from the ones of [22, \S I]. It seems to us that loc. cit. can be used in order to improve these estimates. Accordingly, we now make the connection between loc. cit., [10, \S 2], and 2.4.1.1. We situate ourselves in the context of Case 2 of 2.4.1. Let $O(\al)$ be the ${\bf Z}_p$-algebra of endomorphisms of the unique elementary Dieudonn\'e $p$-divisible object over $k$ of slope $\al$. It is known that $O(\al)$ is an order of the central division algebra over ${\bf Q}_p$ whose invariant is the image of the non-negative rational number $\al$ in ${\bf Q}/{\bf Z}$, cf. [5, Ch. IV, \S3]. There exists a $W(k)$-submodule $\tilde M$ of $M$ that contains $p^{-[-(r_M-1)\al]}M$ and such that the pair $(\tilde M,\vph)$ is an $\al$-divisible $F$-crystal over $k$ (i.e.  for all $n\in {\bf N}$ we can write $\vph^n$ as $p^{[n\al]}$ times a $\sg^n$-linear endomorphism $\vph_n$ of $\tilde M$), cf. [22, pp. 151--152]. As we have $r_M\al\in {\bf N}\cup\{0\}$, all slopes of $(M,\vph_{r_M})$ are $0$. Thus $\vph_{r_M}(\tilde M)=\tilde M$. Triples of the form $(\tilde M,\vph,\vph_{r_M})$ are easily classified. Their isomorphism classes are in one-to-one correspondence to isomorphism classes of torsion free $O(\al)$-modules which by inverting $p$ become free $O(\al)[{1\over p}]$-modules of rank 1, see [10, 2.4 and 2.5]. It is easy to see that under this correspondence, Claim 2.4.1.1 is equivalent to the following well known result.
\medskip
{\bf 2.4.2.1. Claim.} {\it There exists a smallest number $N(\al)\in{\bf N}$ which has effectively computable upper bounds and for which the following property holds: for any element $x\in O(\al)\setminus pO(\al)$, the length of the artinian ${\bf Z}_p$-module $O(\al)/O(\al)x$ is at most $N(\al)$.}
\medskip\smallskip
{\bf 2.5. Standard ${\bf Z}_p$ structures.} Let $(M,\vph,G,(t_{\al})_{\al\in\Mj})$ be a latticed $F$-isocrystal with a group and an emphasized family of tensors over $k$ such that the $W$-condition holds for $(M,\vph,G)$ (see 2.2.1 (b) and (d)). Let $M=\oplus_{i=a}^b \tilde F^i(M)$, $(F^i(M))_{i\in S(a,b)}$, and $\mu:{\bf G}_m\to G$ be as in 2.2.1 (d). Each tensor $t_{\al}\in\Mt(M)$ is fixed by both $\mu$ and $\vph$. Let $\mu_{\rm can}:{\bf G}_m\to {\bf GL}_M$ be the inverse of the canonical split cocharacter of $(M,(F^i(M))_{i\in S(a,b)},\vph)$ defined in [42, p. 512]. Let $M=\oplus_{i=a}^b \tilde F^i_{\rm can}(M)$ be the direct sum decomposition such that the cocharacter $\mu_{\rm can}$ acts on $\tilde F^i_{\rm can}(M)$ via the $-i$-th power of the identity character of ${\bf G}_m$. We have $F^i(M)=\oplus_{j=b}^i \tilde F^j_{\rm can}(M)$ for all $i\in S(a,b)$ and $M=\oplus_{i=a}^b \vph(p^{-i}\tilde F^i_{\rm can}(M))$, cf. loc. cit. The cocharacter $\mu_{\rm can}$ fixes each $t_{\al}$ (cf. the functorial aspects of [42, p. 513]) and so it factors through $G$. As $M=\oplus_{i=a}^b \vph(p^{-i}\tilde F^i_{\rm can}(M))$, the resulting cocharacter $\mu_{\rm can}:{\bf G}_m\to G$ is also a Hodge cocharacter of $(M,\vph,G)$ in the sense of 2.2.1 (d). 
\smallskip
Let $\sg_0:=\vph\mu(p)$. We have $\sg_0(M)=\vph(\oplus_{i=a}^b p^{-i}\tilde F^i(M))=M$. Thus $\sg_0$ is a $\sg$-linear automorphism of $M$ and so also of $\Mt(M)$. For $\al\in\Mj$ we have $\sg_0(t_{\al})=t_{\al}$. Let 
$M_{{\bf Z}_p}:=\{m\in M|\sg_0(m)=m\}$.
We now assume $k=\bar k$. So $M_{{\bf Z}_p}$ is a free ${\bf Z}_p$-module such that we have $M=M_{{\bf Z}_p}\otimes_{{\bf Z}_p} W(k)$ and $t_{\al}\in\Mt(M_{{\bf Z}_p})$ for all $\al\in\Mj$. Let $G_{{\bf Q}_p}$ be the subgroup of ${\bf GL}_{M_{{\bf Z}_p}[{1\over p}]}$ that fixes $t_{\al}$ for all $\al\in\Mj$; its pull back to ${\rm Spec}(B(k))$ is $G_{B(k)}$. Let $G_{{\bf Z}_p}$ be the Zariski closure of $G_{{\bf Q}_p}$ in ${\bf GL}_{M_{{\bf Z}_p}}$. As $G$ is the Zariski closure of $G_{B(k)}$ in ${\bf GL}_M$, we get that $G$ is the pull back to ${\rm Spec}(W(k))$ of $G_{{\bf Z}_p}$. If moreover we have a principal bilinear quasi-polarization $\lambda_M:M\otimes_{W(k)} M\to W(k)$ of $(M,\vph,G)$, then $\lambda_M$ is also the extension to $W(k)$ of a perfect bilinear form $\lambda_{M_{{\bf Z}_p}}$ on $M_{{\bf Z}_p}$.  
\medskip\smallskip
{\bf 2.6. Exponentials.} Let $H={\rm Spec}(A)$ be an integral, affine group scheme of finite type over ${\rm Spec}(W(k))$. Let $O$ be a free $W(k)$-module of finite rank such that we have a closed embedding homomorphism $H\hookrightarrow {\bf GL}_O$; one constructs $O$ as a $W(k)$-submodule of $A$ (cf. [6, Vol. I, Exp. VI${}_B$, 11.11.1]). If $p\ge 3$, let $E_{O}:=p{\rm End}(O)$. If $p=2$, let $E_{O}$ be the sum of $p^2{\rm End}(O)$ and of the set of nilpotent elements of $p{\rm End}(O)$.
Let 
$${\rm exp}:E_{O}\to {\bf GL}_O(W(k))$$ 
be the exponential map that takes $x\in E_{O}$ into $\sum_{i=0}^{\infty} {x^i\over {i!}}$; here $x^0:=1_O$. 
\smallskip
Let $l\in{\bf N}$. Here are the well known properties of the map ${\rm exp}$ we will often use.
\medskip
{\bf (a)} If $p\ge 3$ and $x\in p^l{\rm End}(O)$, then ${\rm exp}(x)$ is congruent mod $p^{2l}$ to $1_{O}+x$.
\smallskip
{\bf (b)} If $p=2$, $l\ge 2$, and $x\in p^l{\rm End}(O)$, then ${\rm exp}(x)$ is congruent mod $p^{2l}$ to $1_{O}+x+{{x^2}\over 2}$ and is congruent mod $p^{2l-1}$ to $1_{O}+x$.
\smallskip
{\bf (c)} If $x\in {\rm Lie}(H_{B(k)})\cap E_{O}$, then ${\rm exp}(x)\in H(W(k))$.  
\medskip
To check (c) it is enough to show that ${\rm exp}(x)\in H(B(k))$. It suffices to check this under the extra assumption that the transcendental degree of $k$ is countable. Fixing an embedding $W(k)\hookrightarrow{\bf C}$, we can view $H({\bf C})$ as a Lie subgroup of ${\bf GL}_O({\bf C})$; so the relation ${\rm exp}(x)\in H(B(k))$ follows easily from [20, Ch. II, \S1, 3].
\medskip
{\bf 2.6.1. Lemma.} {\it Suppose $H$ is smooth over ${\rm Spec}(W(k))$. Let $l\in{\bf N}$. Let $g_l\in H(W(k))$ be congruent mod $p^l$ to $1_{O}$. Then for any $i\in S(1,l)$ there exists $z_{i,l}\in {\rm Lie}(H)$ such that $g_l$ is congruent mod $p^{i+l}$ to $1_{O}+p^lz_{i,l}$.}
\medskip
{\it Proof:} We use induction on $i$. The case $i=1$ is trivial. Let $\bar z_{1,l}$ be the reduction mod $p$ of $z_{1,l}$. The passage from $i$ to $i+1$ goes as follows. We first consider the case when either $p\ge 3$ or $p=2$ and $i+1<l$. Let $g_{l+1}:=g_l{\rm exp}(-p^lz_{1,l})\in H(W(k))$; it is congruent mod $p^{l+1}$ to $(1_{O}+p^lz_{1,l})(1_{O}-p^lz_{1,l})$ (cf. 2.6 (a) and (b)) and so also to $1_O$. By replacing $g_l$ with $g_{l+1}$, the role of the pair $(i+1,l)$ is replaced by the one of the pair $(i,l+1)$. As $1_{O}+p^lz_{1,l}$ and ${\rm exp}(p^lz_{1,l})$ are congruent mod $p^{i+1+l}$ (cf. 2.6 (a) and (b)), by induction we get that $g_l=g_{l+1}{\rm exp}(p^lz_{1,l})$ is congruent mod $p^{i+1+l}$ to $(1_O+p^{l+1}z_{i,l+1})(1_{O}+p^lz_{1,l})$ and so also to $1_O+p^l(z_{1,l}+pz_{i,l+1})$. Thus as $z_{i+1,l}$ we can take the sum $z_{1,l}+pz_{i,l+1}$. 
\smallskip
Let now $p=2$ and $i+1=l\ge 2$. We have $\bar z_{1,l}^2\in {\rm Lie}(H_{k})$, cf. [1, Ch. II, 3.1, 3.5, Lemma 3 of 3.19]. Thus there exists $\tilde z_{1,l}\in {\rm Lie}(H)$ that is congruent mod $2$ to $z_{1,l}^2$. But $1_{O}-2^lz_{1,l}$ is congruent mod $2^{2l}$ to ${\rm exp}(-2^lz_{1,l}){\rm exp}(-2^{2l-1}\tilde z_{1,l})$, cf. 2.6 (b). The existence of $z_{i+1,l}$ is now argued as in the previous paragraph but working with $g_{l+1}:=g_l{\rm exp}(-2^lz_{1,l}){\rm exp}(-2^{2l-1}\tilde z_{1,l})\in H(W(k))$. This ends the induction.\endproof
\medskip
{\bf 2.6.2. Lemma.} {\it Suppose $H$ is smooth over ${\rm Spec}(W(k))$. Let $l\in{\bf N}$. If $z_l\in {\rm Lie}(H)$, then the reduction mod $p^{2l}$ of $1_{O}+p^lz_l$ belongs to $H(W_{2l}(k))$.}
\medskip
{\it Proof:} We can assume $p=2$ (cf. 2.6 (a) and (c) applied with $x=p^lz_l$) and $l\ge 2$ (as $H$ is smooth). By replacing $1_{O}+2^lz_l$ with $(1_{O}+2^lz_l){\rm exp}(-2^lz_l)$, we can assume (cf. 2.6 (b)) that $z_l\in 2^{l-1}{\rm Lie}(H)$. But this case is obvious (as $H$ is smooth).\endproof 
\medskip\smallskip
{\bf 2.7. Dilatations.} In this Subsection we study an arbitrary integral, closed subgroup scheme $G={\rm Spec}(R_G)$ of ${\bf GL}_M$. Let $W(k)^{\rm sh}$ be the strict henselization of $W(k)$. If $a:{\rm Spec}(W(k)^{\rm sh})\to G$ is a morphism, then the N\'eron measure of the defect of smoothness $\delta(a)\in{\bf N}\cup\{0\}$ of $G$ at $a$ is the length of the torsion part of the coherent $\Mo_{{\rm Spec}(W(k)^{\rm sh})}$-module $a^*(\Omega_{G/{\rm Spec}(W(k))})$. Here $\Omega_{G/{\rm Spec}(W(k))}$ is the coherent $\Mo_G$-module of relative differentials of $G$ with respect to ${\rm Spec}(W(k))$. As $G$ is a group scheme, the value of $\delta(a)$ does not depend on $a$ and so we denote it by $\delta(G)$. We have $\delta(G)\in{\bf N}$ if and only if $G$ is not smooth over ${\rm Spec}(W(k))$, cf. [2, 3.3, Lemma 1]. Let $S(G)$ be the Zariski closure in $G_k$ of all special fibres of $W(k)^{\rm sh}$-valued points of $G$. It is a reduced subgroup of $G_k$. We write $S(G)={\rm Spec}(R_G/J_G)$, where $J_G$ is the ideal of $R_G$ that defines $S(G)$. 
\smallskip
By the {\it canonical dilatation} of $G$ we mean the affine $G$-scheme $G_1={\rm Spec}(R_{G_1})$, where $R_{G_1}$ is the $R_G$-subalgebra of $R_G[{1\over p}]$ generated by ${x\over p}$, with $x\in J_G$. The ${\rm Spec}(W(k))$-scheme $G_1$ is integral and has a canonical group scheme structure with respect to which the morphism $G_1\to G$ is a homomorphism of group schemes over ${\rm Spec}(W(k))$, cf. [2, 3.2, p. 63 and (d) of p. 64]. The morphism $G_1\to G$ has the following universal property (cf. [2, 3.2, Prop. 1]): any morphism $Y\to G$ of flat ${\rm Spec}(W(k))$-schemes whose special fibre factors through the closed embedding $S(G)\hookrightarrow G_k$, factors uniquely through $G_1\to G$. Either $\delta(G_1)=0$ (i.e. $G_1$ is smooth over ${\rm Spec}(W(k))$) or (cf. [2, 3.3, Prop. 5]) we have $0<\delta(G_1)<\delta(G)$. 
\smallskip
By using at most $\delta(G)$ canonical dilatations, we get the existence of a unique smooth, affine group scheme $G^\prime$ over ${\rm Spec}(W(k))$ that is endowed with a homomorphism $G^\prime\to G$ whose fibre over ${\rm Spec}(B(k))$ is an isomorphism and that has the following universal property (cf. [2, 7.1, Thm. 5]): any morphism $Y\to G$ of ${\rm Spec}(W(k))$-schemes with $Y$ smooth, factors uniquely through $G^\prime\to G$. In particular, we can identify $G^\prime(W(k)^{\rm sh})$ with $G(W(k)^{\rm sh})$. The homomorphism $G^\prime\to G$ is called the {\it group smoothening} of $G$. Let 
$$n(G)\in S(0,\dl(G))$$ 
be the smallest number of canonical dilatations one has to perform in order to construct $G^\prime$. We have $n(G)=0$ if and only if $G$ is smooth over ${\rm Spec}(W(k))$. 
\smallskip
The closed embedding $i_G:G\hookrightarrow {\bf GL}_M$ gets replaced by a canonical homomorphism $i_{G^\prime}:G^\prime\to {\bf GL}_M$ that factors through $i_G$. We identify ${\rm Lie}(G^\prime)$ with a $W(k)$-lattice of ${\rm Lie}(G_{B(k)})$ contained in ${\rm End}(M)$. Let $d_{\rm sm}\in {\bf N}\cup\{0\}$ be the smallest number such that we have $p^{d_{\rm sm}}({\rm Lie}(G_{B(k)})\cap {\rm End}(M))\subseteq {\rm Lie}(G^\prime)\subseteq {\rm Lie}(G_{B(k)})\cap {\rm End}(M)$. 
\smallskip
We fix a closed embedding homomorphism $G^{\prime}\hookrightarrow {\bf GL}_{M^\prime}$, where $M^\prime$ is a free $W(k)$-module of finite rank (see beginning of 2.6). Let $g\in G^\prime(W(k))=G(W(k))$. 
\medskip
{\bf 2.7.1. Definition.} Let $n\in{\bf N}$. We say $g$ is congruent mod $p^n$ to $1_{M^\prime}$ (resp. to $1_M$) if and only if the image of $g$ in $G^\prime(W_n(k))$ (resp. in $G(W_n(k))$) is the identity element. 
\medskip
{\bf 2.7.2. Lemma.} {\it We have the following three properties:
\medskip
{\bf (a)} If $g$ is congruent mod $p^n$ to $1_{M^\prime}$, then $g$ is also congruent mod $p^n$ to $1_M$.
\smallskip
{\bf (b)} If $g$ is congruent mod $p^{n+n(G)}$ to $1_M$, then $g$ is also congruent mod $p^n$ to $1_{M^\prime}$.
\smallskip
{\bf (c)} We have an inequality $d_{\rm sm}\le n(G)$.}
\medskip
{\it Proof:}
Part (a) is trivial. We write $G^\prime={\rm Spec}(R_{G^\prime})$ and ${\bf GL}_M={\rm Spec}(R_M)$. Let $I_G$, $I_{G_1}$, $I_{G^\prime}$, and $I_M$ be the ideals of $R_G$, $R_{G_1}$, $R_{G^\prime}$, and $R_M$ (respectively) that define the identity sections. We have $I_{G_1}=I_G[{1\over p}]\cap R_{G_1}$ and $I_{G^\prime}=I_G[{1\over p}]\cap R_{G^\prime}$. 
\smallskip
We check (b). Let $m_g:R_G\to W(k)$ be the homomorphism that defines $g$; we have $m_g(I_G)\subseteq p^{n+n(G)}W(k)$. Let $m_{1g}:R_{G_1}\to W(k)$ be the homomorphism through which $m_g$ factors. We have $m_{1g}(I_{G_1})\subseteq p^{n+n(G)-1}W(k)$, cf. the very definition of $R_{G_1}$. Part (b) follows from a repeated application of this fact to the sequence of $n(G)$ dilatations performed to construct $G^\prime$. The cokernel of the cotangent map (computed at $W(k)$-valued identity elements) $I_G/I_G^2\to I_{G_1}/I_{G_1}^2$ is annihilated by $p$, cf. the very definition of $R_{G_1}$. By applying this repeatedly, we get that the cokernel of the cotangent map $I_M/I_M^2\to I_{G^\prime}/I_{G^\prime}^2$ is annihilated by $p^{n(G)}$. Taking duals we get that the cokernel of the $W(k)$-linear Lie monomorphism ${\rm Lie}(G^\prime)\hookrightarrow {\rm Lie}(G_{B(k)})\cap {\rm Lie}({\bf GL}_M)$ is also annihilated by $p^{n(G)}$. As ${\rm Lie}({\bf GL}_M)$ is the Lie algebra associated to ${\rm End}(M)$, we get that (c) holds.\endproof 
\medskip\smallskip
{\bf 2.8. Complements on $\Mm(W_q(S))$.} Let $q\in{\bf N}$ and let $l\in S(0,q)$. Let $f:S_1\to S$ be a morphism of ${\rm Spec}({\bf F}_p)$-schemes. Let $f_q:W_q(S_1)\to W_q(S)$ be the natural morphism of ${\rm Spec}({\bf Z}/p^q{\bf Z})$-schemes defined by $f$. Let ${\got C}$ be an $F$-crystal over $S$. In this Subsection we include four complements on the category $\Mm(W_q(S))$.
\medskip
{\bf 2.8.1. Pulls back.} Let $f_q^*:\Mm(W_q(S))\to \Mm(W_q(S_1))$ be the natural pull back functor. So if $S={\rm Spec}(R)$ and $S_1={\rm Spec}(R_1)$ are affine and if $h:(O,\vph_O)\to (O^\prime,\vph_{O^\prime})$ is a morphism of $\Mm(W_q(S))$, then $f_q^*(h)$ is the morphism 
$$h\otimes 1_{W_q(R_1)}:(O\otimes_{W_q(R)} W_q(R_1),\vph_O\otimes\Phi_{R_1})\to (O^\prime\otimes_{W_q(R)} W_q(R_1),\vph_{O^\prime}\otimes\Phi_{R_1}).$$ 
\indent
In general $W_q(S)\times_{W_{q+1}(S)} W_{q+1}(S_1)$ is not $W_q(S_1)$. Thus, in general the restriction of $f_{q+1}^*$ to $\Mm(W_q(S))$ and $f_q^*$ do not coincide as functors from $\Mm(W_q(S))$ to $\Mm(W_q(S_1))$ and therefore the sequence of functors $(f_q^*)_{q\in{\bf N}}$ does not define a pull back functor from $\Mm(W(S))$ to $\Mm(W(S_1))$. If the Frobenius endomorphism of $\Mo_{S_1}$ is surjective, then regardless of who $S$ is we have $W_q(S)\times_{W_{q+1}(S)} W_{q+1}(S_1)=W_q(S_1)$ and thus the sequence of functors $(f_q^*)_{q\in{\bf N}}$ does define a pull back functor $f^*:\Mm(W(S))\to\Mm(W(S_1))$. 
\smallskip
If $u$ is an object (or a morphism) of $\Mm(W_q(S))$, then by its pull back to an object (or a morphism) of $\Mm(W_q(S_1))$ we mean $f_q^*(u)$. If $t\in{\bf N}$ and if $f_{q+t}^*(u)$ is an object (or a morphism) of $\Mm(W_q(S_1))$, then we have $f_{q+t}^*(u)=f_q^*(u)$. If $S_1$ is the spectrum of a perfect field, we also speak simply of the pull back of $u$ via $f$, to be often denoted as $f^*(u)$ (instead of either $f_q^*(u)$ or $f_{q+t}^*(u)$). 
\smallskip
If $S$ is integral, if $k_S$ is the  field of fractions of $S$, and if $u$ is a morphism of $\Mm(W_q(S))$, then we say ${\rm Coker}(u)$ is {\it generically annihilated by $p^l$} if the pull back of $u$ to a morphism of $\Mm(W_q(k_S))=\Mm(W_q({\rm Spec}(k_S)))$ has a cokernel annihilated by $p^l$. 
\medskip
{\bf 2.8.2. The evaluation functor ${\bf E}$.} Let $\dl_q(S)$ be the canonical divided power structure of the ideal sheaf of $\Mo_{W_q(S)}$ that defines the closed embedding $S\hookrightarrow W_q(S)$. The evaluation of the $F$-crystal ${\got C}$ at the thickening $(S\hookrightarrow W_q(S),\dl_q(S))$ is a triple $(\Mf_q,\vph_{\Mf_q},\nabla_{\Mf_q})$, where $\Mf_q$ is a locally free $\Mo_{W_q(S)}$-module of 
finite rank, where $\vph_{\Mf_q}:\Mf_q\to\Mf_q$ is a $\Phi_S$-linear endomorphism, and where $\nabla_{\Mf_q}$ is an integrable and topologically nilpotent connection on $\Mf_q$, that satisfies certain axioms. In this paper, connections as $\nabla_{\Mf_q}$ will play no role; on the other hand, we will often use the following object 
$${\bf E}({\got C};W_q(S)):=(\Mf_q,\vph_{\Mf_q})$$ 
of $\Mm(W_q(S))$. A morphism $v:{\got C}\to {\got C}_1$ of $F$-crystals over $S$ defines naturally a morphism 
$${\bf E}(v;W_q(S)):{\bf E}({\got C};W_q(S))\to {\bf E}({\got C}_1;W_q(S)).$$ 
The association $v\to {\bf E}(v;W_q(S))$ defines a ${\bf Z}_p$-linear (evaluation) functor from the category of $F$-crystals over $S$ into the category $\Mm(W_q(S))$. 
\smallskip
To ease notations, let ${\bf E}({\got C};W_q(S_1)):={\bf E}({\got C}_{S_1};W_q(S_1))$ and, in the case when $S_1={\rm Spec}(R_1)$ is affine, let ${\bf E}({\got C};W_q(R_1)):={\bf E}({\got C};W_q(S_1))$. 
\smallskip
The functorial morphism $f_q:(S_1\hookrightarrow W_q(S_1),\dl_q(S_1))\to (S\hookrightarrow W_q(S),\dl_q(S))$ gives birth to a canonical isomorphism (to be viewed as an identity)
$$c_{f;q}:f_q^*({\bf E}({\got C};W_q(S)))\tilde\to{\bf E}({\got C};W_q(S_1)).\leqno (10)$$
If $e:S_2\to S_1$ is another morphism of ${\rm Spec}({\bf F}_p)$-schemes, then we have identities
$$(f\circ e)^*_q=e^*_q\circ f_q^*\;\; {\rm and}\;\; c_{e;q}\circ e_q^*(c_{f;q})=c_{f\circ e;q}.\leqno (11)$$
In what follows we will use without any extra comment the identities (10) and (11).   
\medskip
{\bf 2.8.3. Inductive limits.} Let $V\hookrightarrow V_1$ be a monomorphism of commutative ${\bf F}_p$-algebras. Suppose we have an inductive limit $V_1={\rm ind.}\,{\rm lim.}\,_{\al\in\Lambda} V_{\al}$ of commutative $V$-subalgebras of $V_1$ indexed by the set of objects $\Lambda$ of a filtered, small category. For $\al\in\Lambda$, let $f^{\al}:{\rm Spec}(V_{\al})\to {\rm Spec}(V)$ be the natural morphism. 
\smallskip
Let $(O,\vph_{O})$ and $(O^\prime,\vph_{O^\prime})$ be objects of $\Mm(W_q(V))$ such that $O$ and $O^\prime$ are free $W_q(V)$-modules of finite rank. Let $(O_1,\vph_{O_1})$ and $(O_1^\prime,\vph_{O_1^\prime})$ be the pulls back of $(O,\vph_O)$ and $(O^\prime,\vph_{O^\prime})$ (respectively) to objects of $\Mm(W_q(V_1))$. We consider a morphism 
$$u_1:(O_1,\vph_{O_1})\to (O_1^\prime,\vph_{O_1^\prime})$$ 
of $\Mm(W_q(V_1))$ whose cokernel is annihilated by $p^l$. We fix ordered $W_q(V)$-bases $\Mb_O$ and $\Mb_{O^\prime}$ of $O$ and $O^\prime$ (respectively). Let $B_1$ be the matrix representation of $u_1$ with respect to the ordered $W_q(V_1)$-basis of $O_1$ and $O_1^\prime$ defined naturally by $\Mb_O$ and $\Mb_{O^\prime}$ (respectively). As $p^l{\rm Coker}(u_1)=0$, for $x^\prime\in\Mb_{O^\prime}$ we can write $p^lx^\prime\otimes 1=u_1(\sum_{x\in\Mb_O} x\otimes\be_{xx^\prime})$, where each $\be_{xx^\prime}\in W_q(V_1)$. Let $V_{u_1}$ be the $V$-subalgebra of $V_1$ generated by the components of the Witt vectors of length $q$ with coefficients in $V_1$ that are either  entries of $B_1$ or $\be_{xx^\prime}$ for some pair $(x,x^\prime)\in\Mb_O\times\Mb_{O^\prime}$. As $V_{u_1}$ is a finitely generated $V$-algebra, there exists $\al_0\in \Lambda$ such that $V_{u_1}\hookrightarrow V_{\al_0}$. This implies that $u_1$ is the pull back of a morphism 
$$u_{\al_0}:f^{\al_0 *}_q(O,\vph_{O})\to f^{\al_0 *}_q(O^\prime,\vph_{O^\prime})$$ 
of $\Mm(W_q(V_{\al_0}))$ whose cokernel is annihilated by $p^l$. Here are four special cases of interest.
\medskip
{\bf (a)} If $V$ is a field and $V_1$ is an algebraic closure of $V$, then as $V_{\al}$'s we can take the finite field extensions of $V$ that are contained in $V_1$. 
\smallskip
{\bf (b)} If $V_1$ is a local ring of an integral domain $V$, then as $V_{\al}$'s we can take the $V$-algebras of global functions of open, affine subschemes of ${\rm Spec}(V)$ that contain ${\rm Spec}(V_1)$. 
\smallskip
{\bf (c)} We consider the case when $V$ is a discrete valuation ring that is an N-2 ring in the sense of [29, (31.A)], when $V_1$ is a faithfully flat $V$-algebra that is also a discrete valuation ring, and when each $V_{\al}$ is a $V$-algebra of finite type. The flat morphism $f^{\al_0}:{\rm Spec}(V_{\al_0})\to {\rm Spec}(V)$ has quasi-sections, cf. [18, Ch. IV, Cor. (17.16.2)]. In other words, there exists a finite field extension $\tilde k$ of $k$ and a $V$-subalgebra $\tilde V$ of $\tilde k$ such that: (i) $\tilde V$ is a local, faithfully flat $V$-algebra of finite type whose field of fractions is $\tilde k$, and (ii) we have a morphism $\tilde f^{\al_0}:{\rm Spec}(\tilde V)\to {\rm Spec}(V_{\al_0})$ such that $\tilde f:=f^{\al_0}\circ \tilde f^{\al_0}$ is the natural morphism ${\rm Spec}(\tilde V)\to {\rm Spec}(V)$. As $V$ is an N-2 ring, its normalization in $\tilde k$ is a finite $V$-algebra and so a Dedekind domain. This implies that we can assume $\tilde V$ is a discrete valuation ring. For future use, we recall that any excellent ring is a Nagata ring (cf. [29, (34.A)]) and so also an N-2 ring (cf. [29, (31.A)]). Let 
$$\tilde u:{\tilde f_q^{*}}(O,\vph_{O})=\tilde f_q^{\al_0 *}(f^{\al_0 *}_q(O,\vph_{O}))\to \tilde f_q^{*}(O^\prime,\vph_{O^\prime})=\tilde f_q^{\al_0 *}(f^{\al_0 *}_q(O^\prime,\vph_{O^\prime}))$$ 
be the pull back of $u_{\al_0}$ to a morphism of $\Mm(W_q(\tilde V))$; its cokernel is annihilated by $p^l$. 
\smallskip
If $V$ is the local ring of an integral ${\rm Spec}({\bf F}_p)$-scheme $U$, then $\tilde V$ is a local ring of the normalization of $U$ in $\tilde k$. So from (b) we get that there exists an open subscheme $\tilde U$ of this last normalization that has $\tilde V$ as a local ring and that has the property that $\tilde u$ extends to a morphism of $\Mm(W_q(\tilde U))$ whose cokernel is annihilated by $p^l$.
\smallskip
{\bf (d)} If $V$ is reduced and $V_1=V^{\rm perf}$, we can take $\Lambda={\bf N}$ and $V_n=V^{(p^n)}$ ($n\in{\bf N}$).  
\medskip
{\bf 2.8.4. Hom schemes.} Let $\Mo_1$ and $\Mo_2$ be two objects of $\Mm(W_q(S))$ such that their underlying $\Mo_{W_q}(S)$-modules are locally free of finite ranks. We consider the functor 
$${\bf Hom}(\Mo_1,\Mo_2):{\rm Sch}^S\to {\rm SET}$$ 
from the category ${\rm Sch}^S$ of $S$-schemes to the category ${\rm SET}$ of sets, with the property that ${\bf Hom}(\Mo_1,\Mo_2)(S_1)$ is the set underlying the ${\bf Z}/p^q{\bf Z}$-module of morphisms of $\Mm(W_q(S_1))$ that are between $f_q^*(\Mo_1)$ and $f_q^*(\Mo_2)$; here $f:S_1\to S$ is as in the beginning of 2.8. 
\medskip
{\bf 2.8.4.1. Lemma.} {\it The functor ${\bf Hom}(\Mo_1,\Mo_2)$ is representable by an affine $S$-scheme which locally is of finite presentation.} 
\medskip
{\it Proof:} Localizing, we can assume that $S={\rm Spec}(R)$ is affine and that $\Mo_1=(O_1,\vph_{O_1})$ and $\Mo_2=(O_2,\vph_{O_2})$ are such that $O_1$ and $O_2$ are free $W_q(R)$-modules. For $i\in\{1,2\}$ let $r_i$ be the rank of $O_i$. Let ${\bf Hom}(O_1,O_2)$ be the affine space (of relative dimension $qr_1r_2$) over ${\rm Spec}(R)$ with the property that for any commutative $R$-algebra $R_1$, ${\bf Hom}(O_1,O_2)(R_1)$ is the set of $W_q(R_1)$-linear maps $x:O_1\otimes_{W_q(R)}W_q(R_1)\to O_2\otimes_{W_q(R)}W_q(R_1)$. We have an identity $(\vph_{O_2}\otimes\Phi_{R_1})\circ x=x\circ (\vph_{O_1}\otimes\Phi_{R_1})$ if and only if $x$ belongs to the subset of ${\bf Hom}(O_1,O_2)(R_1)$ that is naturally identified with ${\bf Hom}(\Mo_1,\Mo_2)({\rm Spec}(R_1))$. As the relation $(\vph_{O_2}\otimes\Phi_{R_1})\circ x=x\circ (\vph_{O_1}\otimes\Phi_{R_1})$ defines a closed subscheme of ${\bf Hom}(O_1,O_2)$ that is of finite presentation, the Lemma follows.\endproof
\medskip\smallskip
{\bf 2.9. On two results of commutative algebra.} In 5.4 and \S6 we will use the following two geometric variations of well known results of commutative algebra.\medskip
{\bf 2.9.1. Lemma.} {\it Let $X$ and $Y$ be two integral, normal, locally noetherian schemes. Let $u:X\to Y$ be an affine morphism that is birational; let $K$ be the field of fractions of either $X$ or $Y$. Let $D(X)$ and $D(Y)$ be the set of local rings of $X$ and $Y$ (respectively) that are discrete valuation rings (of $K$). If $D(Y)\subseteq D(X)$, then $u$ is an isomorphism.}
\medskip
{\it Proof:} Working locally in the Zariski topology of $Y$, we can assume $X={\rm Spec}(R_X)$ and $Y={\rm Spec}(R_Y)$ are also affine and noetherian. Thus (inside $K$) we have 
$$R_Y\hookrightarrow R_X=\cap_{V\in D(X)} V\hookrightarrow \cap_{V\in D(Y)} V=R_Y$$ 
(cf. [29, (17.H)] for the two identities). So $R_Y=R_X$. Thus $u$ is an isomorphism.\endproof
\medskip
{\bf 2.9.2. Lemma.} {\it Let $X^\prime={\rm Spec}(R^\prime)\to X={\rm Spec}(R)$ be a morphism between affine schemes which at the level of rings is defined by an integral (i.e. an ind-finite) monomorphism $R\hookrightarrow R^\prime$. Then an open subscheme $U$ of $X$ is affine if and only if its pull back $U^\prime:=U\times_X X^\prime$ to $X^\prime$ is affine.}
\medskip
{\it Proof:} It is enough to show that $U$ is affine if $U^\prime$ is affine. The morphism $U^\prime\to U$ is surjective (see [29, (5.E)]). So as $U^\prime$ is quasi-compact (being affine), $U$ is also quasi-compact. Thus $X^{\rm top}\setminus U^{\rm top}$ is the zero locus in $X^{\rm top}$ of a finite number of elements of $R$. So there exists a finitely generated ${\bf Z}$-subalgebra $R_0$ of $R$ such that $U$ is the pull back of an open subscheme $U_0$ of ${\rm Spec}(R_0)$ through the natural morphism ${\rm Spec}(R)\to {\rm Spec}(R_0)$. 
\smallskip
Let $\Lambda$ (resp. $\Lambda^\prime$) be the set of finite subsets of $R$ (resp. of $R^\prime$). For $\al\in\Lambda$ (resp. $\al^\prime\in\Lambda^\prime$), let $R_{\al}$ (resp. $R_{\al^\prime}^\prime$) be the $R_0$-subalgebra of $R$ (resp. of $R^\prime$) generated by $\al$ (resp. by $\al^\prime$). Let $X_{\al}:={\rm Spec}(R_{\al})$ and $X^\prime_{\al^\prime}:={\rm Spec}(R^\prime_{\al^{\prime}})$. Let $U_{\al}$ and  $U^\prime_{\al^\prime}$ be the pulls back of $U_0$ to $X_{\al}$ and $X^\prime_{\al^\prime}$ (respectively). As $U_{\al^\prime}^\prime$ is a quasi-compact, open subscheme of $X^\prime_{\al^\prime}$, it is an $X_{\al^\prime}^\prime$-scheme of finite presentation. As the scheme $U^\prime$ is affine, by applying [18, Ch. IV, (8.10.5)] to the projective limit $U^\prime\hookrightarrow X^\prime={\rm proj.}\,{\rm lim.}\,_{\al^\prime\in \Lambda^\prime} U_{\al^\prime}^\prime\hookrightarrow X_{\al^\prime}^\prime$ of open embeddings of finite presentation, we get that there exists $\be^\prime\in\Lambda^\prime$ such that $U_{\be^\prime}^\prime$ is affine. Let $\be\in\Lambda$ be such that $R^{\prime}_{\be^\prime\cup\be}$ is a finite $R_{\be}$-algebra. As $U_{\be^\prime\cup\be}^\prime=U^\prime_{\be^\prime}\times_{X^\prime_{\be^\prime}} X_{\be^\prime\cup\be}^\prime $ is affine, the scheme $U_{\be}$ is also affine (cf. Chevalley theorem of [16, Ch. II, (6.7.1)] applied to the finite, surjective morphism $U_{\be^\prime\cup\be}^\prime\to U_{\be}$). Thus $U=U_{\be}\times_{X_{\be}} X$ is affine.\endproof\bigskip\smallskip
\centerline{\bigsll {\bf \S3 Proof of Main Theorem A and complements}}
\bigskip\smallskip
In 3.1 we prove Main Theorem A stated in 1.2. See formula (18) of 3.1.3 for a concrete expression of the number $n_{\rm fam}$ mentioned in 1.2. In 3.2 we include interpretations and variants of 1.2 in terms of reductions modulo powers of $p$; in particular, see 3.2.4 for the passage from 1.2 to 1.3. See 3.3 for improvements of the estimates of 3.1.1 to 3.1.5 in many particular cases of interest. If $p\ge 3$ let $\vep_p:=1$. Let $\vep_2:=2$. 
\medskip\smallskip
{\bf 3.1. Proof of 1.2.} We start the proof of 1.2. Until 3.1.4 we will assume $k=\bar k$. Let $v:=\dim(G_{B(k)})$. It suffices to prove 1.2 under the extra hypothesis $v\ge 1$. Let $\dl(G)$, $n(G)\in {\bf N}\cup\{0\}$, the group smoothening $G^\prime\to G$ of $G$, and the closed embedding homomorphism $G^\prime\hookrightarrow {\bf GL}_{M^\prime}$ be as in 2.7. We have $G^{\prime}(W(k))=G(W(k))$. 
\smallskip
Let $m:={\bf T}({\rm Lie}(G^\prime),\vph)\in {\bf N}\cup\{0\}$. Based on definitions 2.2.2 (a) and (b), there exists a $B(k)$-basis $\Mb=\{e_1,\ldots,e_v\}$
of ${\rm Lie}(G_{B(k)})$ formed by elements of ${\rm Lie}(G^\prime)$ and there exists a permutation $\pi$ of $S(1,v)$, such that the following three things hold:
\medskip
{\bf (a)} denoting $E:=<e_1,\ldots,e_v>$, we have $p^m{\rm Lie}(G^\prime)\subseteq E\subseteq {\rm Lie}(G^\prime)$;
\smallskip
{\bf (b)} if $l\in S(1,v)$, then we have $\vph(e_l)=p^{n_l}e_{\pi(l)}$ for some $n_l\in{\bf Z}$;
\smallskip
{\bf (c)} for any cycle $\pi_0=(l_1,\ldots,l_q)$ of $\pi$, the integers $n_{l_1},\ldots,n_{l_q}$ are either all non-negative or all non-positive.
\medskip
If we have $n_{l_j}\ge 0$ for all $j\in S(1,q)$, let $\tau(\pi_0):=1$. If there exists $j\in S(1,q)$ such that $n_{l_j}<0$, let $\tau(\pi_0):=-1$. Let $n\in{\bf N}$ be such that 
$$n\ge 2m+\vep_p+n(G).\leqno (12)$$ 
\indent  
Let $g_n\in G(W(k))$ be congruent mod $p^{n}$ to $1_M$. So $g_n\in G^{\prime}(W(k))$ is congruent mod $p^{n-n(G)}$ to $1_{M^\prime}$ (cf. 2.7.2 (b)) and below we will only use this congruence. 
\medskip
{\bf 3.1.1. Claim.} {\it For any $t\in{\bf N}$ there exists $\tilde g_t\in G^\prime(W(k))=G(W(k))$ congruent mod $p^{n-n(G)+t-1-m}$ to $1_{M^\prime}$ and such that $\tilde g_tg_n\vph \tilde g_t^{-1}\vph^{-1}\in G^\prime(W(k))=G(W(k))$ is congruent mod $p^{n-n(G)+t}$ to $1_{M^\prime}$. Thus there exists $\tilde g_0\in G^\prime(W(k))=G(W(k))$ congruent mod $p^{n-n(G)-m}$ to $1_{M^\prime}$ and such that $\tilde g_0g_n\vph=\vph\tilde g_0$.}
\medskip
If $t\ge 2$ and if the element $\tilde g_{t-1}$ exists, then we can replace the triple $(n,t,g_n)$ by the triple $(n+t-1,1,\tilde g_{t-1}g_n\vph \tilde g_{t-1}^{-1}\vph^{-1})$. Thus using induction on $t\in {\bf N}$, to prove the first part of the Claim we can assume that $t=1$. As $W(k)$ is $p$-adically complete, the second part of the Claim follows from its first part; this is so as we can take $\tilde g_0$ to be an infinite product of the form $\cdots\tilde h_3\tilde h_2\tilde h_1$ that has the property that for all $c\in {\bf N}$ the element $\tilde h_c\in G^\prime(W(k))=G(W(k))$ is congruent mod $p^{n-n(G)+c-1-m}$ to $1_{M^\prime}$ and moreover $\tilde h_c\tilde h_{c-1}\cdots \tilde h_1g_n\vph \tilde h_1^{-1}\cdots \tilde h_c^{-1}\vph^{-1}\in G^\prime(W(k))=G(W(k))$ is congruent mod $p^{n-n(G)+c}$ to $1_{M^\prime}$. Thus to prove the Claim, it suffices to prove its first part for $t=1$. 
\smallskip
For $t=1$ we will use what we call the {\it stairs method} for $E$. Let $z_n\in {\rm Lie}(G^\prime)$ be such that $g_n$ is congruent mod $p^{n-n(G)+1}$ to $1_{M^\prime}+p^{n-n(G)}z_n$. As $n-n(G)\ge m+1$, based on 3.1 (a) we can write 
$$p^{n-n(G)}z_n=\sum_{l\in S(1,v)} p^{u_l}c_le_l,$$
where $u_l\in{\bf N}$ depends only on the cycle of $\pi$ to which $l$ belongs and where $c_l\in W(k)$. We take the $u_l$'s to be the maximal possible values subject to the last sentence. Thus 
$$u_l\ge n-n(G)-m\ge m+\vep_p\ge\vep_p\ge 1.\leqno (13a)$$
From (13a) and (12) we get
$${\rm min}\{u_l+u_{l^\prime}|l,l^\prime\in S(1,v)\}\ge 2(n-n(G)-m)\ge n-n(G)+\vep_p\ge n-n(G)+1.\leqno (13b)$$
\indent
Due to (13b), the product $\tilde g^{(1)}_n:=\prod_{l\in S(1,v)} (1_{M^\prime}+p^{u_l}c_le_l)\in {\bf GL}_{M^\prime}(W(k))$
is congruent mod $p^{2(n-n(G)-m)}$ and so also mod $p^{n-n(G)+1}$ to $1_{M^\prime}+\sum_{l\in S(1,v)} p^{u_l}c_le_l=1_{M^\prime}+p^{n-n(G)}z_n$. 
The element $\tilde g_n^{(2)}:=(\tilde g_n^{(1)})^{-1}g_n\in {\bf GL}_{M^\prime}(W(k))$ is congruent mod $p^{n-n(G)+1}$ to $(1_{M^\prime}+p^{n-n(G)}z_n)^{-1}(1_{M^\prime}+p^{n-n(G)}z_n)=1_{M^\prime}$. We have $g_n=\tilde g_n^{(1)}\tilde g_n^{(2)}$. \smallskip 
For $l\in S(1,v)$ let $q_l:=-{\rm min}\{0,n_l\}\in{\bf N}\cup\{0\}$. We will choose $\tilde g_1\in G^{\prime}(W(k))$ to be a product $\prod_{l\in S(1,v)} {\rm exp}(p^{u_l+q_l}x_le_l)$, with all $x_l$'s in $W(k)$. This last product makes sense, cf. 2.6 (a) and (b) and the fact that for $p=2$ we have $u_l+q_l\ge u_l\ge m+\vep_2\ge\vep_2=2$. 
\smallskip
For $l\in S(1,v)$ we have $u_l=u_{\pi(l)}$. Thus
$$\vph \tilde g_1^{-1}\vph^{-1}=\prod_{l\in S(1,v)} {\rm exp}(-p^{u_l+q_l}\sg(x_l)\vph(e_l))=\prod_{l\in S(1,v)} {\rm exp}(-p^{u_l+q_{\pi^{-1}(l)}+n_{\pi^{-1}(l)}}\sg(x_{\pi^{-1}(l)})e_l).\leqno (14)$$
These exponential elements are well defined even if $p=2$, as for $p=2$ we have inequalities $u_l+q_{\pi^{-1}(l)}+n_{\pi^{-1}(l)}\ge u_l\ge\vep_2\ge 2$. Thus $\vph \tilde g_1^{-1}\vph^{-1}\in G^{\prime}(W(k))=G(W(k))$, cf. 2.6 (c). 
\smallskip
We have $\vph \tilde g_1^{-1}\vph^{-1}\in G^{\prime}(W(k))=G(W(k))$, cf. 2.6 (c). We have to show that we can choose the $x_l$'s such that $\tilde g_1g_n\vph \tilde g_1^{-1}\vph^{-1}\in G^{\prime}(W(k))$ is congruent mod $p^{n-n(G)+1}$ to $1_{M^\prime}$. It suffices to show that if $g_n$ is not congruent mod $p^{n-n(G)+1}$ to $1_{M^\prime}$, then we can choose the $x_l$'s such that by replacing $g_n$ with $g_{n+1}\tilde g_1g_n\vph \tilde g_1^{-1}\vph^{-1}\in G^\prime(W(k))$, where $g_{n+1}\in G^{\prime}(W(k))$ is congruent mod $p^{n-n(G)+1}$ to $1_{M^\prime}$, we can also replace each $u_l$ by $u_l+t_l$, where $t_l\in{\bf N}$ depends only on the cycle of $\pi$ to which $l$ belongs. 
\smallskip
The element $\tilde g_1g_n\vph \tilde g_1^{-1}\vph^{-1}\in G^{\prime}(W(k))$ is congruent mod $p^{n-n(G)+1}$ to the product $\tilde g_1\tilde g_n^{(1)}\vph\tilde g_1^{-1}\vph^{-1}$. From (13a), (13b), (14), and 2.6 (a) and (b), we get that $\vph \tilde g_1^{-1}\vph^{-1}$ is congruent mod $p^{n-n(G)+1}$ to $\prod_{l\in S(1,v)} [1_{M^\prime}-p^{u_l+q_{\pi^{-1}(l)}+n_{\pi^{-1}(l)}}\sg(x_{\pi^{-1}(l)})e_l]$. A similar argument shows that $\tilde g_1$ is congruent mod $p^{n-n(G)+1}$ to $\prod_{l\in S(1,v)} (1_{M^\prime}+p^{u_l+q_l}x_le_l)$. Thus the product $\tilde g_1\tilde g_n^{(1)}\vph\tilde g_1^{-1}\vph^{-1}$ of the three elements $\tilde g_1$, $\tilde g_n^{(1)}$, and $\vph\tilde g_1^{-1}\vph^{-1}$, is congruent mod $p^{n-n(G)+1}$ to the following product of three elements 
$$\prod_{l\in S(1,v)} (1_{M^\prime}+p^{u_l+q_l}x_le_l)\prod_{l\in S(1,v)} (1_{M^\prime}+p^{u_l}c_le_l)\prod_{l\in S(1,v)} [1_{M^\prime}-p^{u_l+q_{\pi^{-1}(l)}+n_{\pi^{-1}(l)}}\sg(x_{\pi^{-1}(l)})e_l]$$ 
and so (cf. (13b)) also to $1_{M^\prime}+\sum_{l\in S(1,v)} p^{u_l}[p^{q_l}x_l+c_l-p^{q_{\pi^{-1}(l)}+n_{\pi^{-1}(l)}}\sg(x_{\pi^{-1}(l)})]e_l$. 
\smallskip
To show that we can take each $t_l$ to be at least 1, it suffices to show that we can choose the $x_l$'s such that we have 
$$p^{q_l}x_l+c_l-p^{q_{\pi^{-1}(l)}+n_{\pi^{-1}(l)}}\sg(x_{\pi^{-1}(l)})\in pW(k)\;\; \forall l\in S(1,v).\leqno (15)$$ 
In other words, by denoting with $\bar x\in k$ the reduction mod $p$ of an arbitrary element $x\in W(k)$, it suffices to show that for each cycle $\pi_0=(l_1,\ldots,l_q)$ of $\pi$ there exist solutions in $k$ of the following circular system of equations over $k$
$$\bar b_{l_j}\bar x_{l_j}+\bar c_{l_j}-\bar d_{l_j}\bar x_{l_{j-1}}^p=0\;\;\;{\rm with}\;\;\;j\in S(1,q),\leqno (16)$$
where $b_{l_j}:=p^{q_{l_j}}$ and $d_{l_j}:=p^{q_{l_{j-1}}+n_{l_{j-1}}}$ (here we have $l_0=l_q$, cf. end of 2.1). 
\smallskip
If $\tau(\pi_0)=1$, then $q_{l_j}=q_{l_{j-1}}=0$ and $n_{l_{j-1}}\ge 0$; so $p^{q_{l_j}}x_{l_j}+c_{l_j}-p^{q_{l_{j-1}}+n_{l_{j-1}}}\sg(x_{l_{j-1}})$ is $x_{l_j}+c_{l_j}-p^{n_{l_{j-1}}}\sg(x_{l_{j-1}})$. If $\tau(\pi_0)=-1$, then $q_{l_{j-1}}=-n_{l_{j-1}}\ge 0$ and so we have $p^{q_{l_j}}x_{l_j}+c_{l_j}-p^{q_{l_{j-1}}+n_{l_{j-1}}}\sg(x_{l_{j-1}})=p^{q_{l_j}}x_{l_j}+c_{l_j}-\sg(x_{l_{j-1}})$; moreover, there exists $j_0\in S(1,j)$ such that $q_{l_{j_0}}=-n_{l_{j_0}}>0$. Thus depending on the value of $\tau(\pi_0)$ we have:
\medskip
{\bf (+)} $\bar b_{l_j}=1$ and $\bar d_{l_j}\in\{0,1\}$ for all $j\in S(1,q)$, if $\tau(\pi_0)=1$;
\smallskip
{\bf (--)} $\bar d_{l_j}=1$ and $\bar b_{l_j}\in\{0,1\}$  for all $j\in S(1,q)$ and moreover there exists $j_0\in S(1,q)$ such that $\bar b_{l_{j_0}}=0$, if $\tau(\pi_0)=-1$. 
\medskip
If $\tau(\pi_0)=1$, then based on (+) we can eliminate the variables $\bar x_{l_q}$, $\bar x_{l_{q-1}},\ldots, \bar x_{l_3}$, and $\bar x_{l_2}$ one by one from the system (16). The resulting equation in the variable $\bar x_{l_1}$ is of the form $\bar x_{l_1}=\bar u_{l_1}+\bar v_{l_1}\bar x_{l_1}^{p^q}$, where $\bar u_{l_1}$ and $\bar v_{l_1}\in k$. This equation defines an \'etale $k$-algebra. Thus (as $k$ is separably closed) the system (16) has solutions in $k$ if $\tau(\pi_0)=1$. 
\smallskip
If $\tau(\pi_0)=-1$, then based on (--) (and on the fact that $k$ is perfect) the values of $\bar x_{l_{j_0-1}}$, $\bar x_{l_{j_0-2}},\ldots,\bar x_{l_1}$, $\bar x_{l_q}$, $\bar x_{l_{q-1}},\ldots,\bar x_{l_{j_0}}$ are one by one uniquely determined and so the system (16) has a unique solution. 
\smallskip
This ends the proof of the first part of the Claim for $t=1$ and so it also ends the proof of the Claim. 
\medskip
{\bf 3.1.2. Inequalities involving $s$- and $h$-numbers.} Let $s^\prime_L$ and $h^\prime_L$ be the $s$-number and the $h$-number  (respectively) of $({\rm Lie}(G^\prime),\vph)$. Let $s_L$ and $h_L$ be the $s$-number and the $h$-number  (respectively) of $({\rm Lie}(G_{B(k)})\cap {\rm End}(M),\vph)$. We recall from 2.7 that $d_{\rm sm}\in {\bf N}\cup\{0\}$ is the smallest number such that we have inclusions
$$p^{d_{\rm sm}}({\rm Lie}(G_{B(k)})\cap {\rm End}(M))\subseteq {\rm Lie}(G^\prime)\subseteq {\rm Lie}(G_{B(k)})\cap {\rm End}(M).\leqno (17a)$$ 
We have $p^{s_L}\vph({\rm Lie}(G_{B(k)})\cap {\rm End}(M))\subseteq {\rm Lie}(G_{B(k)})\cap {\rm End}(M)$. Thus $p^{s_L+d_{\rm sm}}\vph({\rm Lie}(G^\prime))\subseteq p^{s_L+d_{\rm sm}}\vph({\rm Lie}(G_{B(k)})\cap {\rm End}(M))\subseteq p^{d_{\rm sm}}({\rm Lie}(G_{B(k)})\cap {\rm End}(M))\subseteq {\rm Lie}(G^\prime)$. From the very definition of $s_L^\prime$ we get 
$$s_L^\prime\le s_L+d_{\rm sm}.\leqno (17b)$$ 
\noindent
The $h$-numbers of $({\rm Lie}(G_{B(k)})\cap {\rm End}(M),p^{{\rm max}\{s_L,s^\prime_L\}}\vph)$ and $({\rm Lie}(G^\prime),p^{{\rm max}\{s_L,s^\prime_L\}}\vph)$ are $h_L+{\rm max}\{s_L,s^\prime_L\}-s_L$ and $h_L^\prime+{\rm max}\{s_L,s^\prime_L\}-s_L^\prime$ (respectively). From this and $(17a)$ we easily get that 
$$h_L^\prime+{\rm max}\{s_L,s^\prime_L\}-s_L^\prime\le d_{\rm sm}+h_L+{\rm max}\{s_L,s^\prime_L\}-s_L.\leqno (17c)$$
From $(17b)$, $(17c)$, and the inequalities $d_{\rm sm}\le n(G)\le\dl(G)$ (see 2.7.2 (c) and 2.7), we get
$$h^\prime_L\le h_L+d_{\rm sm}+s_L^\prime-s_L\le h_L+d_{\rm sm}+d_{\rm sm}\le h_L+2\dl(G).\leqno (17d)$$
\indent
{\bf 3.1.3. End of the proof of 1.2.} As $v=\dim(G_{B(k)})$, $G^\prime$, and $n(G)$ depend only on $G$ and as $s^\prime_L$ and $h^\prime_L$ depend only on the family $\{(M,g\vph,G)|g\in G(W(k))=G^\prime(W(k))\}$ of latticed $F$-isocrystals with a group over $k$, the number 
$$
n_{\rm fam}:=2d(v,s^\prime_L,h^\prime_L)+\vep_p+n(G)\leqno (18)
$$
is not changed if $\vph$ gets replaced by $g\vph$ for some $g\in G(W(k))$. As $m\le d(v,s^\prime_L,h^\prime_L)$ (cf. 2.4.1), we have $n_{\rm fam}\ge 2m+\vep_p+n(G)$. So from (12) applied with $n=n_{\rm fam}$ and from 3.1.1, we get that for any $g_{n_{\rm fam}}\in G(W(k))$ congruent mod $p^{n_{\rm fam}}$ to $1_M$, there exists an isomorphism between $(M,\vph,G)$ and $(M,g_{n_{\rm fam}}\vph,G)$ defined by an element of $G^\prime(W(k))=G(W(k))$. This property holds even if $\vph$ gets replaced by $g\vph$.  As $n(G)\le\dl(G)$, we have 
$$n_{\rm fam}\le 2d(v,s^\prime_L,h^\prime_L)+\vep_p+\dl(G).\leqno (19)$$
So based on the effectiveness part of 2.4.1 and on $(17b)$ and $(17d)$, to end the proof of 1.2 it is enough to show that $\dl(G)$, $s_L$, and $h_L$ are effectively bounded from above. But $\dl(G)$ is effectively computable in terms of the ideal sheaf of $\Mo_{{\bf GL}_M}$ that defines the closed embedding homomorphism $G\hookrightarrow {\bf GL}_M$, cf. [2, 3.3, Lemma 2]. As the numbers $s_L$ and $h_L$ are effectively computable in terms of $({\rm Lie}(G_{B(k)}),\vph)$ and ${\rm End}(M)$ and as the connected group $G_{B(k)}$ is uniquely determined by its Lie algebra (cf. [1, Ch. II, 7.1]), the numbers $s_L$ and $h_L$ are also effectively computable in terms of the closed embedding homomorphism $G\hookrightarrow {\bf GL}_M$. This ends the proof of 1.2.
\medskip
{\bf 3.1.4. Definition.} Let $(M,\vph,G)$ be a latticed $F$-isocrystal with a group over $k$. By the {\it isomorphism number} (to be abbreviated as the {\it $i$-number}) of $(M,\vph,G)$ we mean the smallest number $n\in{\bf N}\cup\{0\}$ such that for any $g_n\in G(W(\bar k))$ congruent mod $p^n$ to $1_{M\otimes_{W(k)} W(\bar k)}$, there exists an isomorphism between $(M\otimes_{W(k)} W(\bar k),\vph\otimes\sg_{\bar k},G_{W(\bar k)})$ and $(M\otimes_{W(k)} W(\bar k),g_n(\vph\otimes\sg_{\bar k}),G_{W(\bar k)})$ which is an element of $G(W(\bar k))$. If $G={\bf GL}_M$, we also refer to $n$ as the $i$-number of $(M,\vph)$. 
\smallskip
If $(M,\vph,G,\lambda_M)$ and $G^0$ are as in 2.2.1 (c), then by the $i$-number of $(M,\vph,G,\lambda_M)$ we mean the $i$-number of its latticed $F$-isocrystal with a group $(M,\vph,G^0)$ over $k$.
\medskip
{\bf 3.1.5. Example.} Suppose $G$ is smooth over ${\rm Spec}(W(k))$ and $k=\bar k$. Thus $n(G)=d_{\rm sm}=0$, $s_L^\prime=s_L$, $h_L^\prime=h_L$, and $n_{\rm fam}=2d(v,s_L,h_L)+\vep_p$. Let $m:={\bf T}({\rm Lie}(G^\prime),\vph)$. If $n$ is as in 3.1.4, then we have $0\le n\le 2m+\vep_p\le 2d(v,s_L,h_L)+\vep_p$ (cf. 3.1.1 and 3.1.3). 
\medskip\smallskip
{\bf 3.2. Truncations.} In 3.2.1 to 3.2.7 we work in a context that pertains to Dieudonn\'e modules and to $p$-divisible groups. In particular, in 3.2.1 and 3.2.2 we define and study $D$-truncations that are the crystalline analogues (with a group) over $k$ of truncated Barsotti--Tate groups over ${\rm Spec}(k)$. In 3.2.8 we consider reductions modulo powers of $p$ of those $F$-crystals with a group $(M,\vph,G)$ over $k$ for which $G$ is smooth over ${\rm Spec}(W(k))$. In 3.2.9 we introduce the $F$-truncations that generalize the $D$-truncations.
\medskip
{\bf 3.2.1. On $D$-truncations.} Let $(r,d)\in {\bf N}\times ({\bf N}\cup\{0\})$. Let $(M,\vph,G)$ be an $F$-crystal with a group over $k$. Until 3.2.8 we assume that $(M,\vph)$ is a Dieudonn\'e module over $k$, that $r_M=r$, that $d$ is the dimension of the kernel of $\vph$ mod $p$, that $G$ is smooth over ${\rm Spec}(W(k))$, and that the $W$-condition holds for $(M,\vph,G)$. Let $M=\tilde F^0(M)\oplus \tilde F^1(M)$ be a direct sum decomposition such that $M=\vph(\tilde F^0(M)\oplus {1\over p}\tilde F^1(M))$ and the cocharacter $\mu:{\bf G}_m\to {\bf GL}_M$ that fixes $\tilde F^0(M)$ and that acts via the inverse of the identical character of ${\bf G}_m$ on $\tilde F^1(M)$, factors through $G$ (cf. 2.2.1 (d)). The rank of $\tilde F^1(M)$ is $d$. Let $\sg_0:=\vph\mu(p)$; it is a $\sg$-linear automorphism of either $M$ and $\Mt(M)$ (cf. 2.5). As $\sg_0$ normalizes ${\rm Lie}(G_{B(k)})$ and ${\rm End}(M)$, it also normalizes ${\rm Lie}(G)={\rm Lie}(G_{B(k)})\cap {\rm End}(M)$. As $G(W(k))=G(B(k))\cap {\bf GL}_M(W(k))$, we have $\sg_0G(W(k))\sg_0^{-1}=G(W(k))$. 
\smallskip
Let $q\in{\bf N}$. Let $\vth:M\to M$ be the $\sg^{-1}$-linear endomorphism that is the Verschiebung map of $(M,\vph)$; we have $\vth\vph=\vph \vth=p1_M$. We denote also by $\vth$ its reduction mod $p^q$. By the {\it $D$-truncation} mod $p^q$ of $(M,\vph,G)$ we mean the quadruple 
$$(M/p^qM,\vph,\vth,G_{W_q(k)}).$$ 
This quadruple determines (resp. is determined by) the reduction mod $p^q$ (resp. mod $p^{q+1}$) of $(M,\vph,G)$. We also refer to $(M/p^qM,\vph,\vth)$ as the $D$-truncation mod $p^q$ of $(M,\vph)$. If $(M_1/p^qM_1,\vph_1,\vth_1,G_{1W_q(k)})$ is a similar $D$-truncation mod $p^q$, then by an isomorphism $f:(M/p^qM,\vph,\vth,G_{W_q(k)})\tilde\to (M_1/p^qM_1,\vph_1,\vth_1,G_{1W_q(k)})$ 
we mean a $W_q(k)$-linear isomorphism $f:M/p^qM\tilde\to M_1/p^qM_1$ such that we have identities $f\vph=\vph_1f$ and $f\vth=\vth_1f$ and the isomorphism ${\bf GL}_{M/p^qM}\tilde\to {\bf GL}_{M_1/p^qM_1}$ induced by $f$, takes $G_{W_q(k)}$ onto $G_{1W_q(k)}$. 
\smallskip
If $(M_1,\vph_1,G_1)$ is $(M,g\vph,G)$ with $g\in G(W(k))$ and if $f\in G(W_q(k))$, then we say $f$ is an {\it inner isomorphism} between the $D$-truncation mod $p^q$ of $(M,\vph,G)$ and $(M,g\vph,G)$.
\medskip
{\bf 3.2.2. Lemma.} {\it For $g\in G(W(k))$ the following two statements are equivalent:
\medskip
{\bf (a)} the $D$-truncations mod $p^q$ of $(M,\vph,G)$ and $(M,g\vph,G)$ are inner isomorphic;
\smallskip
{\bf (b)} there exists an element $\tilde g\in G(W(k))$ such that $\tilde gg\vph \tilde g^{-1}=g_q\vph$, where $g_q\in G(W(k))$ is congruent mod $p^q$ to $1_M$.}
\medskip
{\it Proof:} As $G$ is smooth over ${\rm Spec}(W(k))$, the reduction homomorphism $G(W(k))\to G(W_q(k))$ is onto. Thus it suffices to check that (a) implies (b) under the extra assumptions that the $\sg$-linear endomorphisms $\vph$ and $g\vph$ coincide mod $p^q$ and that the $\sg^{-1}$-linear endomorphisms $\vth$ and $\vth g^{-1}$ coincide mod $p^q$. Let $g_0:=\sg_0^{-1}g\sg_0 \in G(W(k))$. Let $g_{0,q}\in G(W_q(k))$ be $g_0$ mod $p^q$. We have $\vph=\sg_0\mu(p^{-1})$ and $g\vph=\sg_0 g_0\mu(p^{-1})$. We get that the two endomorphisms $\mu(p^{-1})$ and  $g_0\mu(p^{-1})$ of $M$ coincide mod $p^q$. Thus \medskip
{\bf (i)} the element $g_{0,q}$ fixes $\tilde F^0(M)/p^q\tilde F^0(M)$ and is congruent mod $p^{q-1}$ to $1_{M/p^qM}$. 
\medskip
We have $\vth=p\mu(p)\sg_0^{-1}$ and $\vth g^{-1}=p\mu(p)g_0^{-1}\sg_0^{-1}$. We get that the two endomorphisms $p\mu(p)$ and $p\mu(p)g_0^{-1}$ of $M$ coincide mod $p^q$. Thus $p\mu(p)$ and $p\mu(p)g_0$ also coincide mod $p^q$ and therefore we have an inclusion 
\medskip
{\bf (ii)} 
 $(1_{M/p^qM}-g_{0,q})(\tilde F^1(M)/p^q\tilde F^1(M))\subseteq p^{q-1}\tilde F^0(M)/p^q\tilde F^0(M)$.
\medskip
Let $\tilde F^{-1}({\rm End}(M))$ be the maximal direct summand of ${\rm End}(M)$ on which $\mu$ acts via the identity cocharacter of ${\bf G}_m$. Thus $\tilde F^{-1}({\rm End}(M))$ is the ${\rm Hom}(\tilde F^1(M),\tilde F^0(M))$ factor of the following direct sum decomposition ${\rm End}(M)={\rm End}(\tilde F^1(M))\oplus {\rm End}(\tilde F^0(M))\oplus {\rm Hom}(\tilde F^1(M),\tilde F^0(M))\oplus {\rm Hom}(\tilde F^0(M),\tilde F^1(M))$ into $W(k)$-submodules. 
\smallskip
Let $U_{-1}^{\rm big}$ be the closed subgroup scheme of ${\bf GL}_M$ defined by the rule: if $A$ is a commutative $W(k)$-algebra, then $U_{-1}^{\rm big}(A):=1_{M\otimes_{W(k)} A}+\tilde F^{-1}({\rm End}(M))\otimes_{W(k)} A\leqslant {\bf GL}_M(A)$. So $U_{-1}^{\rm big}$ is the maximal subgroup scheme of ${\bf GL}_M$ that fixes both $\tilde F^0(M)$ and $M/\tilde F^0(M)$; it is smooth over ${\rm Spec}(W(k))$. We have ${\rm Lie}(U_{-1}^{\rm big})=\tilde F^{-1}({\rm End}(M))$. Let $U_{-1}$ be the closed subgroup scheme of $U_{-1}^{\rm big}$ defined by the rule: 
$$U_{-1}(A):=1_{M\otimes_{W(k)} A}+({\rm Lie}(G)\cap \tilde F^{-1}({\rm End}(M)))\otimes_{W(k)} A.$$ 
The group scheme $U_{-1}$ is smooth, unipotent, has connected fibres and its Lie algebra is the direct summand ${\rm Lie}(G)\cap \tilde F^{-1}({\rm End}(M))$ of $\tilde F^{-1}({\rm End}(M))$. As $U_{-1B(k)}$ is connected and ${\rm Lie}(U_{-1})\subseteq {\rm Lie}(G)$, the group $U_{-1B(k)}$ is a closed subgroup of $G_{B(k)}$ (cf. [1, Ch. II, 7.1]). Thus $U_{-1}$ is a closed subgroup scheme of $G$. As $\mu$ factors through $G$, we have two identities ${\rm Lie}(G)\cap {\rm Lie}(U_{-1}^{\rm big})={\rm Lie}(U_{-1})$ and ${\rm Lie}(G_k)\cap {\rm Lie}(U_{-1k}^{\rm big})={\rm Lie}(U_{-1k})$. Thus the intersection $U_{-1}^\prime:=G\cap U_{-1}^{\rm big}$ has smooth fibres and has $U_{-1}$ as its identity component. As the group $U_{-1B(k)}^{\rm big}/U_{-1B(k)}$ has no non-trivial finite subgroups, we have $U_{-1B(k)}=U_{-1B(k)}^\prime$. As $U_{-1}^{\rm big}$ is a complete intersection in ${\bf GL}_M$, $U_{-1}^\prime$ has dimension at least equal to $1+\dim(U_{-1k})$ at each $k$-valued point. So as $\dim(U^\prime_{-1k})=\dim(U_k)$, the group $U_{-1B(k)}^\prime=U_{-1B(k)}$ is Zariski dense in $U_{-1}^\prime$. Thus $U_{-1}=U_{-1}^\prime=G\cap U_{-1}^{\rm big}$.
\smallskip
We have $g_{0,q}\in {\rm Ker}(U_{-1}^{\rm big}(W_q(k))\to U_{-1}^{\rm big}(W_{q-1}(k)))$, cf. (i) and (ii). As $U_{-1}=G\cap U_{-1}^{\rm big}$, in fact we have $g_{0,q}\in {\rm Ker}(U_{-1}(W_q(k))\to U_{-1}(W_{q-1}(k)))$. Thus, up to a replacement of $g$ by a ${\rm Ker}(G(W(k))\to G(W_q(k)))$-multiple of it, we can assume that $g_0\in U_{-1}(W(k))$. We write $g_0=1_M+p^{q-1}u_{-1}$, where $u_{-1}\in {\rm Lie}(U_{-1})$. 
\smallskip
Let $\tilde g:=1_M+p^qu_{-1}\in U_{-1}(W(k))\leqslant G(W(k))$. As we have $\tilde g^{-1}=1_M-p^qu_{-1}$ and $g_0^{-1}=1_M-p^{q-1}u_{-1}$, we get $\mu(p^{-1})\tilde g^{-1}\mu(p)=g_0^{-1}$. Thus we compute that $\tilde gg\vph \tilde g^{-1}$ is 
$\tilde gg\sg_0\mu(p^{-1})\tilde g^{-1}\mu(p)\sg_0^{-1}\vph=\tilde gg\sg_0 g_0^{-1}\sg_0^{-1}\vph=\tilde gg(\sg_0g\sg_0^{-1})^{-1}\vph=\tilde ggg^{-1}\vph=\tilde g\vph$. So as $g_q$ we can take $\tilde g$. Thus (a) implies (b). Obviously (b) implies (a).\endproof
\medskip
{\bf 3.2.3. Corollary.} {\it Suppose $k=\bar k$ and $(M,\vph)$ is the Dieudonn\'e module of a $p$-divisible group $D$ over ${\rm Spec}(k)$ of height $r=r_M$ and dimension $d$. Let $n\in {\bf N}\cup\{0\}$ be the $i$-number of $(M,\vph)$. Then $n$ is the smallest number $t\in {\bf N}\cup\{0,\infty\}$ for which the following statement holds:
\medskip
{\bf (*)}  if $D_1$ is a $p$-divisible group over ${\rm Spec}(k)$ of height $r$ and dimension $d$ and if $D_1[p^t]$ is isomorphic to $D[p^t]$, then $D_1$ is isomorphic to $D$.}
\medskip
{\it Proof:}
The Dieudonn\'e module of $D_1$ is isomorphic to $(M,g\vph)$ for some $g\in {\bf GL}_M(W(k))$; moreover any such pair $(M,g\vph)$ is isomorphic to the Dieudonn\'e module of some $p$-divisible group over ${\rm Spec}(k)$ of height $r$ and dimension $d$. For $q\in{\bf N}$, the classical Dieudonn\'e theory achieves also a natural one-to-one and onto correspondence between the isomorphism classes of truncated Barsotti--Tate groups of level $q$ over $k$ and the isomorphism classes of $D$-truncations mod $p^q$ of Dieudonn\'e modules over $k$ (see [14, pp. 153 and 160]). So (*) holds if and only if for any $g\in {\bf GL}_M(W(k))$, the fact that the $D$-truncations mod $p^t$ of $(M,\vph)$ and $(M,g\vph)$ are isomorphic implies that $(M,\vph)$ and $(M,g\vph)$ are isomorphic. From this and 3.2.2 (applied with $G={\bf GL}_M$ and $q=t$), we get that (*) holds if and only if $(M,\vph)$ and $(M,g\vph)$ are isomorphic for all elements $g\in {\bf GL}_M(W(k))$ congruent mod $p^t$ to $1_M$. Thus (*) holds if and only if $t\ge n$.\endproof
\medskip
{\bf 3.2.4. Proof of 1.3.} To prove 1.3, we can assume that $k=\bar k$ and that $(M,\vph)$ is the Dieudonn\'e module of $D$. Let $n$ be the $i$-number of $(M,\vph)$. Let $n_{\rm fam}$ be as in 3.1.5 for $G={\bf GL}_M$. We have $n\le n_{\rm fam}$, cf. 3.1.5. From 3.2.3 we get that $D$ is uniquely determined up to isomorphism by $D[p^n]$ and so also by $D[p^{n_{\rm fam}}]$. Thus the number $T(r,d)$ exists and we have $T(r,d)\le n_{\rm fam}$. From 3.2.3 we also get that $T(r,d)\ge n$. If $d\in\{0,r\}$, then $T(r,d)=n_{\rm fam}=n=0$. If $d\notin\{0,r\}$, then the $s$-number and the $h$-number of $({\rm End}(M),\vph)$ are $1$ and $2$ (respectively) and ${\rm End}(M)$ has rank $r^2$.  Thus for $d\in S(1,r-1)$ we have (cf. 3.1.5)
$$n\le T(r,d)\le n_{\rm fam}\le 2d(r^2,1,2)+\vep_p.$$ 
\indent
So $T(r,d)$ is effectively bounded from above in terms of $r$, cf. 2.4.1. This proves 1.3.
\medskip
{\bf 3.2.5. Principal quasi-polarizations.} 
Suppose $r=2d=r_M$, $k=\bar k$, and we have a principal quasi-polarization $\lambda_M$ of $(M,\vph)$. We refer to the triple $(M,\vph,\lambda_M)$ as a principally quasi-polarized Dieudonn\'e module over $k$. Let $G:={\bf Sp}(M,\lambda_M)$. Let $n\in {\bf N}\cup\{0\}$ be the $i$-number of $(M,\vph,G)$. Let $(D,\lambda_D)$ be a principally quasi-polarized $p$-divisible group over ${\rm Spec}(k)$ whose principally quasi-polarized Dieudonn\'e module is isomorphic to $(M,\vph,\lambda_M)$. The principally quasi-polarized Dieudonn\'e module of any other principally quasi-polarized $p$-divisible group over ${\rm Spec}(k)$ of height $r=r_M$ is of the form $(M,g\vph,\lambda_M)$ for some $g\in G(W(k))$. So as in the proof of 3.2.3 we argue that $n$ is the smallest number $t\in {\bf N}\cup\{0,\infty\}$ for which the following statement holds:
\medskip
{\bf (*)} {\it if $(D_1,\lambda_{D_1})$ is a principally quasi-polarized $p$-divisible group over ${\rm Spec}(k)$ of height $r=2d$ and if the principally quasi-polarized truncated Barsotti--Tate groups of level $t$ of $(D_1,\lambda_{D_1})$ and $(D,\lambda_D)$ are isomorphic, then $(D_1,\lambda_{D_1})$ is isomorphic to $(D,\lambda_D)$}.
\medskip
As in 3.2.4 we argue that there exists a smallest number $T(d)\in {\bf N}$ such that any principally quasi-polarized $p$-divisible group over ${\rm Spec}(k)$ of height $r=2d$ is uniquely determined up to isomorphism by its principally quasi-polarized truncated Barsotti--Tate group of level $T(d)$. The number $T(d)$ is effectively bounded from above in terms of the relative dimension $2d^2+d$ of $G={\bf Sp}(M,\lambda_M)$ over ${\rm Spec}(W(k))$ and so also of $d$ itself. 
\smallskip
From the very definition of $T(d)$ we get:
\medskip
{\bf 3.2.6. Corollary.} {\it Suppose $k=\bar k$. Let $D$ be a $p$-divisible group over ${\rm Spec}(k)$ of height $r=2d$ and dimension $d$. Let $D^{\rm t}$ be the Cartier dual of $D$. Then the number of isomorphism classes of principally quasi-polarized $p$-divisible groups of the form $(D,\lambda_D)$ is bounded from above by the finite number of distinct truncations of level $T(d)$ of isomorphisms $D\tilde\to D^{\rm t}$, i.e. by the number of elements of the following finite set of isomorphisms ${\rm Im}({\rm Isom}(D,D^{\rm t})\to {\rm Isom}(D[p^{T(d)}],D^{\rm t}[p^{T(d)}]))$.} 
\medskip
{\bf 3.2.7. Proposition.} {\it Suppose $k=\bar k$. Let $R$ be the normalization of $k[[w]]$ in an algebraic closure $K$ of $k((w))$. For $*\in\{k,K\}$ let $D_*$ be a $p$-divisible group over ${\rm Spec}(*)$ of height $r$ and dimension $d$. Then $D_k$ is the specialization (via ${\rm Spec}(R)$) of a $p$-divisible group over ${\rm Spec}(K)$ which is isomorphic to $D_K$ if and only if $D_k[p^{T(r,d)}]$ is the specialization (via ${\rm Spec}(R)$) of a truncated Barsotti--Tate group of level $T(r,d)$ over ${\rm Spec}(K)$  which is isomorphic to $D_K[p^{T(r,d)}]$.}
\medskip
{\it Proof:} It suffices to check the if part. Let $\Mg_R$ be a truncated Barsotti--Tate group of level $T(r,d)$ over ${\rm Spec}(R)$ that lifts $D_k[p^{T(r,d)}]$ and that has the property that its fibre over ${\rm Spec}(K)$ is isomorphic to $D_K[p^{T(r,d)}]$. Let $R_0$ be a finite $k[[w]]$-subalgebra of $R$ such that $\Mg_R$ is the pull back of a truncated Barsotti--Tate group $\Mg_{R_0}$ of level $T(r,d)$ over ${\rm Spec}(R_0)$. As $R_0$ is a complete discrete valuation ring, there exists a $p$-divisible group $D_{R_0}^\prime$ over ${\rm Spec}(R_0)$ that lifts both $D_k$ and $\Mg_{R_0}$ (cf. [21, 4.4 f)]). The pull back of $D_{R_0}^\prime$ to ${\rm Spec}(K)$ is isomorphic to $D_K$ (cf. 1.3) and it specializes (via ${\rm Spec}(R)$) to $D_k$.\endproof 
\medskip
We have a logical variant of 3.2.7 in the principally quasi-polarized context. We next consider $F$-crystals with a group over $k$. 
\medskip
{\bf 3.2.8. Theorem.} {\it Suppose $k=\bar k$. Let $(M,\vph,G)$ be an $F$-crystal with a group over $k$; so $\vph(M)\subseteq M$. Let $h$ be the $h$-number of $(M,\vph)$. Let $m:={\bf T}({\rm Lie}(G),\vph)$ and let $n:=2m+\vep_p$. Let $g\in G(W(k))$. Let $t\in{\bf N}\cup\{0\}$. 
\medskip
{\bf (a)} Suppose $G$ is smooth over ${\rm Spec}(W(k))$ and there exists $\tilde g_{h+n+t}\in G(W_{h+n+t}(k))$ which is an isomorphism between the reductions mod $p^{h+n+t}$ of $(M,g\vph)$ and $(M,\vph)$. Then there exists $\tilde g_0\in G(W(k))$ which is an isomorphism between $(M,g\vph,G)$ and $(M,\vph,G)$ and whose image in $G(W_{h+n+t}(k))$ is congruent mod $p^{n-m+t}$ to $\tilde g_{h+n+t}$.
\smallskip
{\bf (b)} If $G={\bf GL}_M$, then the images of the two reduction homomorphisms ${\rm Aut}(M,\vph)\to {\rm Aut}(M/p^{n-m+t}M,\vph)$ and ${\rm Aut}(M/p^{n+h+t}M,\vph)\to {\rm Aut}(M/p^{n-m+t}M,\vph)$ are the same.}
\medskip
{\it Proof:} Part (b) is a practical application of (a) for the case when $G={\bf GL}_M$. As $G(W(k))$ surjects onto $G(W_{h+n+t}(k))$, it suffices to prove (a) under the extra assumption that $\tilde g_{h+n+t}=1_{M/p^{h+n+t}M}$. So $g$ mod $p^{h+n+t}$ fixes ${\rm Im}(\vph(M)\to M/p^{h+n+t}M)$. But $p^hM\subseteq\vph(M)$ and so $g$ mod $p^{h+n+t}$ fixes also the $W_{h+n+t}(k)$-submodule $p^hM/p^{h+n+t}M$ of $M/p^{h+n+t}M$. Thus $g$ is congruent mod $p^{n+t}$ to $1_M$. So the element $\tilde g_0$ exists, cf. 3.1.1 (applied with $G=G^\prime$, $M=M^\prime$, and $n(G)=0$).\endproof
\medskip 
{\bf 3.2.9. On $F$-truncations.} One can generalize the $D$-truncations of 3.2.1 as follows. Let $(M,\vph,G)$ be a $p$-divisible object with a group over $k$ and let $(F^i(M))_{i\in S(a,b)}$ be a lift of it (see 2.2.1 (b) and (d)). Let $\vph_i:F^i(M)\to M$ be the $\sg$-linear map defined by the rule $\vph_i(x)=p^{-i}\vph(x)$, where $x\in F^i(M)$. We denote also by $\vph_i$ its reduction mod $p^q$. By an {\it inner isomorphism} between $(M/p^qM,(F^i(M)/p^qF^i(M))_{i\in S(a,b)},(\vph_i)_{i\in S(a,b)},G_{W_q(k)})$ and a similar quadruple $(M/p^qM,(F^i_1(M)/p^qF^i_1(M))_{i\in S(a,b)},((g\vph)_i)_{i\in S(a,b)},G_{W_q(k)})$ defined by some lift $(F^i_1(M))_{i\in S(a,b)}$ of $(M,g\vph,G)$ (with $g\in G(W(k))$), we mean an arbitrary element $f\in {\rm Im}(G(W(k))\to G(W_q(k)))$ that has the following two properties:
\medskip
{\bf (i)} it takes $F^i(M)/p^qF^i(M)$ onto $F^i_1(M)/p^qF^i_1(M)$ for all $i\in S(a,b)$; 
\smallskip
{\bf (ii)} we have $f\vph_i=(g\vph)_if$ for all $i\in S(a,b)$.
\medskip
By the {\it $F$-truncation} mod $p^q$ of $(M,\vph,G)$ we mean the set $F_q(M,\vph,G)$ of inner isomorphism classes of quadruples $(M/p^qM,(F^i(M)/p^qF^i(M))_{i\in S(a,b)},(\vph_i)_{i\in S(a,b)},G_{W_q(k)})$ 
we obtain by allowing $(F^i(M))_{i\in S(a,b)}$ to run through all lifts of $(M,\vph,G)$. 
\smallskip
If $(M,\vph)$ is a Dieudonn\'e module over $k$ and if $G$ is smooth, then it is easy to see that the $D$-truncations mod $p^q$ of $(M,\vph,G)$ and $(M,g\vph,G)$ are inner isomorphic if and only if the $F$-truncations mod $p^q$ of $(M,\vph,G)$ and $(M,g\vph,G)$ are inner isomorphic.
\medskip\smallskip
{\bf 3.3. Refinements of 3.1.1.} In many particular cases we can choose the $W(k)$-span $E$ of 3.1.1 (a) to be stable under products and this can lead to significant improvements of the inequalities we obtained in 3.1.1 to 3.1.5. For the sake of generality, we now formalize such improvements in a relative context. 
\smallskip
Let $(M,\vph,G)$ be a latticed $F$-isocrystal with a group over $k$. Until \S4 we assume $k=\bar k$. Until \S4 we also assume that there exists an integral, closed subgroup scheme $G_1$ of ${\bf GL}_M$ which contains $G$ and for which the following two conditions hold:
\medskip
{\bf (i)} the triple $(M,\vph,G_1)$ is a latticed $F$-isocrystal with a group over $k$;
\smallskip
{\bf (ii)} there exists a $B(k)$-basis $\Mb_1:=\{e_1,\ldots,e_{v_1}\}$ of ${\rm Lie}(G_{1B(k)})$ and a permutation $\pi_1$ of $S(1,v_1)$, that have the following four properties:
\medskip\noindent
{\bf (ii.a)} $E_1:=<e_1,\ldots,e_{v_1}>$ is a $W(k)$-submodule of ${\rm End}(M)$ such that $E_1E_1\subseteq E_1$;
\smallskip\noindent
{\bf (ii.b)} each cycle $(l_1,\ldots,l_{q})$ of $\pi_1$ has the property that for $j\in S(1,q)$ we have $\vph(e_{l_j})=p^{n_{l_j}}e_{l_{j+1}}$, where the integers $n_{l_j}$'s are either all non-negative or all non-positive;
\smallskip\noindent
{\bf (ii.c)} if $v:=\dim(G_{B(k)})$ (so $v\in S(1,v_1)$), then the intersection $E:=E_1\cap {\rm Lie}(G_{B(k)})$ has $\{e_1,\ldots, e_v\}$ as a $W(k)$-basis and moreover the permutation $\pi_1$ normalizes $S(1,v)$;
\smallskip\noindent
{\bf (ii.d)} for any $t\in{\bf N}$ and every element of $G(W(k))$ that has the form $1_M+\sum_{l\in S(1,v_1)} p^tx_le_l$, where all $x_l$'s belong to $W(k)$, we have $x_l\in pW(k)$ for all $l\in S(v+1,v_1)$.
\medskip
Let $$E_2:=E_1[{1\over p}]\cap {\rm End}(M)={\rm Lie}(G_{1B(k)})\cap {\rm End}(M).$$ 
From (ii.a) we get that $1_M+E_2$ is a semigroup with identity contained in ${\rm End}(M)$.
\medskip
{\bf 3.3.1. Lemma.} {\it There exists a closed subgroup scheme $G_2$ of ${\bf GL}_M$ that is defined by the rule: if $A$ is a commutative $W(k)$-algebra, then $G_2(A)$ is the group of invertible elements of the semigroup with identity $1_{M\otimes_{W(k)} A}+E_2\otimes_{W(k)} A$. We have $G_1=G_2$.}
\medskip
{\it Proof:} We show that $G_2$ is an integral, closed subgroup scheme of ${\bf GL}_M$ and that ${\rm Lie}(G_{2B(k)})\subseteq E_2[{1\over p}]=E_1[{1\over p}]$. If $1_M\in E_2$, then $1_M+E_2=E_2$ is a $W(k)$-subalgebra of ${\rm End}(M)$ which as a $W(k)$-submodule is a direct summand; so obviously $G_2$ is an integral, closed subgroup scheme of ${\bf GL}_M$ and we have ${\rm Lie}(G_{2B(k)})\subseteq E_2[{1\over p}]=E_1[{1\over p}]$. 
\smallskip
We now consider the case when $1_M\notin E_2$. Let $E_3:={\rm End}(M)\cap (E_2[{1\over p}]+B(k)1_M)$; it is a $W(k)$-subalgebra of ${\rm End}(M)$ that has $E_2$ as a two-sided ideal. The finite $W(k)$-algebra $E_3/E_2$ is isomorphic to a $W(k)$-subalgebra of $B(k)1_M$ and so to $W(k)1_M$. So as $W(k)$-modules, we have a direct sum decomposition $E_3=E_2\oplus W(k)1_M$. Let $G_3$ be the integral, closed subgroup scheme of ${\bf GL}_M$ of invertible elements of $E_3$. Let $x\in 1_{M\otimes_{W(k)} A}+E_2\otimes_{W(k)} A$ be an element that has an invertible determinant. The inverse $x^{-1}\in {\rm End}(M)$ of $x$ belongs to $G_3(A)$ and moreover modulo the ideal $E_2\otimes_{W(k)} A$ of $E_3\otimes_{W(k)} A$ it is $1_{M\otimes_{W(k)} A}$. Thus $x^{-1}\in 1_{M\otimes_{W(k)} A}+E_2\otimes_{W(k)} A$. This implies that the group $G_2(A)$ is the group of all elements of $1_{M\otimes_{W(k)} A}+E_2\otimes_{W(k)} A$ that have an invertible determinant. From this description of points of $G_2$, we get that $G_2$ is an integral, closed subgroup scheme of either ${\bf GL}_M$ or $G_3$ and that we have ${\rm Lie}(G_{2B(k)})\subseteq E_2[{1\over p}]=E_1[{1\over p}]$.
\smallskip
If $x\in E_1$, then we have $1_M+p^tx\in G_2(W(k))$ for all $t>>0$. Thus $E_1[{1\over p}]\subseteq {\rm Lie}(G_{2B(k)})$ and therefore we have identities ${\rm Lie}(G_{1B(k)})=E_1[{1\over p}]=E_2[{1\over p}]={\rm Lie}(G_{2B(k)})$. So $G_{1B(k)}=G_{2B(k)}$, cf. [1, Ch. II, 7.1]. Thus $G_1=G_2$. \endproof
\medskip
{\bf 3.3.2. Theorem.} {\it We recall that conditions 3.3 (i) and (ii) hold. Let $m_1\in{\bf N}\cup\{0\}$ be the smallest number such that $p^{m_1}(E_2)\subseteq E_1\subseteq E_2$. Let $g\in G(W(k))\cap (1_M+p^jE_1)$ with $j\in {\bf N}$. If $p=2$ we assume that either $G=G_1$ or $j\ge 2$. We have:
\medskip
{\bf (a)} There exists $\tilde g\in G(W(k))\cap (1_M+p^jE_1)$ which is an isomorphism between $(M,g\vph,G)$ and $(M,\vph,G)$.
\smallskip
{\bf (b)} The $i$-number of $(M,\vph,G)$ is at most $m_1+1$.}
\medskip
{\it Proof:} To prove (a), we write $g=1_M+\sum_{l\in S(1,v_1)} p^jx_l(j)e_l$, where all $x_l(j)$'s belong to $W(k)$. By induction on $t\in \{j,j+1,\ldots\}$ we construct an element $\tilde g_t\in G(W(k))$ such that there exist elements $x_l(t)$'s in $W(k)$ that satisfy the identity 
$$\tilde g_t\tilde g_{t-1}\cdots \tilde g_jg\vph \tilde g_j^{-1}\cdots \tilde g_{t-1}^{-1}\tilde g_t^{-1}\vph^{-1}=1_M+\sum_{l\in S(1,v_1)} p^tx_l(t)e_l.\leqno (20a)$$ 
If $t=j$ let $\tilde g_t:=1_M$. The passage from $t$ to $t+1$ goes as follows. For $l\in S(1,v_1)$ let $q_l:=-\min\{0,n_l\}$; we recall that $n_l\in{\bf Z}$ is such that $\vph(e_l)=p^{n_l}e_{\pi_1(l)}$ (cf. 3.3 (ii.b)). 
\smallskip
We first consider the case $p\ge 3$. Let 
$$\tilde g_{t+1}:=\prod_{l\in S(1,v)}{\rm exp}(p^{t+q_l}\tilde x_l(t+1)e_l)\in G(W(k)),\leqno (20b)$$ 
where all $\tilde x_l(t+1)$'s belong to $W(k)$. As $\pi_1$ normalizes $S(1,v)$, we have 
$$\vph\tilde g_{t+1}^{-1}\vph^{-1}=\prod_{l\in S(1,v)}{\rm exp}(-p^{t+q_l+n_l}\sg(\tilde x_l(t+1))e_{\pi_1(l)})\in G(W(k))\leqno (20c)$$ 
(to be compared with (14)). As $E_1E_1\subseteq E_1$ and $p\ge 3$, for any $e\in E_1$ the element ${\rm exp}(p^te)=1_M+\sum_{i=1}^{\infty} {p^{ti}\over {i!}}e^i$ belongs to $1_M+p^tE_1$. From this and the inequalities $q_l\ge 0$ and $q_l+n_l\ge 0$, we get that any exponential element of either (20b) or (20c) belongs to $1_M+p^tE_1$. Thus, as $1_M+p^tE_1$ is a semigroup contained in ${\rm End}(M)$, we get that $\tilde g_{t+1}\in 1_M+p^tE_1$ and that (cf. also (20a)) we can write
$$\tilde g_{t+1}(\tilde g_t\cdots\tilde g_jg\vph\tilde g_j^{-1}\cdots\tilde g_t^{-1}\vph^{-1})\vph \tilde g_{t+1}^{-1}\vph^{-1}=1_M+\sum_{l\in S(1,v_1)} p^tx_l^\prime(t+1)e_l,\leqno (20d)$$ 
where all $x_l^\prime(t+1)$'s belong to $W(k)$. If $i>v$, then $x_l^\prime(t+1)\in pW(k)$ (cf. 3.3 (ii.d)). 
\smallskip
Based on 2.6 (a), from $(20b)$, $(20c)$, and $(20d)$ we get that for any $l\in S(1,v)$ the Witt vector $x_l^\prime(t+1)\in W(k)$ is congruent mod $p$ to the sum (to be compared with (15))
$$p^{q_l}\tilde x_l(t+1)+x_l(t)-p^{n_{\pi_1^{-1}(l)}+q_{\pi_1^{-1}(l)}}\sg(\tilde x_{\pi_1^{-1}(l)}(t+1)).$$ 
As in the part of the proof of 3.1.1 that pertains to (+) and (--), we argue that we can choose the $\tilde x_l(t+1)$'s with $l\in S(1,v)$ such that we have $x_l^\prime(t+1)\in pW(k)$ for all $l\in S(1,v)$. Thus for all $l\in S(1,v_1)$ we can write $p^tx_l^\prime(t+1)=p^{t+1}x_l(t+1)$, where $x_l(t+1)\in W(k)$. This ends the induction for $p\ge 3$. 
\smallskip
Let now $p=2$. For $t\ge 2$, we have $2t-1\ge t+1$. So the above passage from $t$ to $t+1$ has to be modified only if $t=j=1$, cf. 2.6 (b). As $E_1E_1\subseteq E_1$, for any $e\in E_1$ the element $1_M+pe$ has an inverse in $1_M+pE_1$ and thus (cf. 3.3.1) it is an element of $G_2(W(k))=G_1(W(k))$. Thus, if $t=j=1$ and $G=G_1$, we can take 
$$\tilde g_{2}:=1_M+\sum_{l\in S(1,v)}p^{1+q_l}\tilde x_l(2)e_l\in G_2(W(k))=G_1(W(k))=G(W(k))\leqno (20e)$$ 
and we can proceed as above. This ends the induction. 
\smallskip
The infinite product $\tilde g:=\cdots \tilde g_{j+2}\tilde g_{j+1}\tilde g_j\in G(W(k))\cap (1_M+p^jE_1)$ is well defined (as $W(k)$ is $p$-adically complete). Passing to limit $t\to\infty$ in $(20a)$, we get $\tilde g g\vph\tilde g^{-1}\vph^{-1}=1_M$. Thus $\tilde g$ is an isomorphism between $(M,g\vph,G)$ and $(M,\vph,G)$. So (a) holds.
\smallskip
We prove (b). For $g\in G(W(k))\leqslant G_1(W(k))$, we have $g-1_M\in E_2$ (cf. 3.3.1). If $g$ is congruent mod $p^{m_1+1}$ to $1_M$, then $g-1_M\in p^{m_1+1}E_2\subseteq pE_1$. Thus $g\in G(W(k))\cap (1_M+pE_1)$. So there exists an isomorphism between $(M,g\vph,G)$ and $(M,\vph,G)$ which is an element of $G(W(k))$, cf. (a) applied with $j=1$. So (b) follows from definition 3.1.4.\endproof  
\medskip
{\bf 3.3.3. Variant.} Suppose $G=G_1$ and $E_1E_1=0$. Then 3.3.2 (a) holds even if $j=0$, i.e. even if we only have $g\in G(W(k))\cap (1_M+E)$. This is so as we have ${\rm exp}(x)=1+x$ for any $x\in E=E_1$. Thus the proof of 3.3.2 (a) holds even if $j=0$ and so the proof of 3.3.2 (b) can be adapted to get that the $i$-number of $(M,\vph,G)$ is at most $m_1$. 
\medskip
{\bf 3.3.4. Corollary.} {\it Let $m_1$ be as in 3.3.2. We assume that $G=G_1$ and that all slopes of $({\rm Lie}(G_{B(k)}),\vph)$ are $0$. Then the $i$-number of $(M,\vph,G)$ is at most $m_1$.}
\medskip
{\it Proof:} Let $E_{{\bf Z}_p}:=\{x\in {\rm Lie}(G_{B(k)})\cap {\rm End}(M)|\vph(x)=x\}$. We can assume $E=E_1$ is the $W(k)$-span of $E_{{\bf Z}_p}$. We have $E_{{\bf Z}_p}E_{{\bf Z}_p}\subseteq E_{{\bf Z}_p}$, cf. 3.3 (ii.a). Let $g_{m_1}\in G(W(k))$ be congruent mod $p^{m_1}$ to $1_M$. We write $g_{m_1}=1_M+e$, with $e\in E$ (to be compared with the proof of 3.3.2 (b)). Based on 3.3.2 (a), to prove the Corollary it suffices to check that there exists $\tilde g\in G(W(k))$ such that $\tilde gg_{m_1}\vph\tilde g^{-1}\vph^{-1}\in G(W(k))\cap (1_M+pE)$. 
\smallskip
The element $1+e=g_{m_1}\in G(W(k))$ normalizes $E_2$. Thus $(1_M+e)E\subseteq E\subseteq E_2=(1_M+e)E_2$. So by reasons of length of artinian modules we get that $(1_M+e)E=E$. Let $e^\prime\in E$ be such that $(1_M+e)e^\prime=-e$. Thus $(1_M+e)(1_M+e^\prime)=1_M+e-e=1_M$. So $(1_M+e)^{-1}=1_M+e^\prime\in 1_M+E=1_M+E_{{\bf Z}_p}\otimes_{{\bf Z}_p} W(k)$. 
\smallskip
Let  $H$ be the group scheme over ${\rm Spec}({\bf Z}_p)$ defined by the rule: if $A$ is a commutative ${\bf Z}_p$-algebra, then $H(A)$ is the group of invertible elements of the semigroup with identity $1_{M\otimes_{W(k)} A}+E_{{\bf Z}_p}\otimes_{{\bf Z}_p} A$ contained in ${\rm End}(M)\otimes_{W(k)} A$. The automorphism $\sg$ of $W(k)$ acts naturally on $H(W(k))$. If $u\in H(W(k))$, then $u\in G(W(k))$ (cf. 3.3.1) and $\sg(u)=\sg u\sg^{-1}$ is $\vph u\vph^{-1}\in H(W(k))$. Moreover we have $g_{m_1}\in H(W(k))$, cf. previous paragraph. 
\smallskip
Let $\bar g_{m_1}$ be the image of $g_{m_1}$ in $H(k)$. The scheme $H$ is an open subscheme of the affine space over ${\rm Spec}({\bf Z}_p)$ that is of relative dimension $v$ and that is defined naturally by $E_{{\bf Z}_p}$. Thus the group scheme $H$ over ${\rm Spec}({\bf Z}_p)$ has connected fibres and is smooth. Moreover, the special fibre $H_{{\bf F}_p}$ is a quasi-affine group and thus also an affine group over ${\rm Spec}({\bf F}_p)$ (cf. [6, Vol. I, Exp. VI${}_B$, 11.11]). Let $\bar{\tilde g}\in H(k)$ be such that $\bar{\tilde g}\bar g_{m_1}\sg(\bar{\tilde g})^{-1}$ is the identity element of $H(k)$, cf. Lang theorem (see [1, Ch. V, 16.4]). Let $\tilde g\in H(W(k))$ be an element that lifts $\bar{\tilde g}$. The element $\tilde gg_{m_1}\vph\tilde g^{-1}\vph^{-1}=\tilde gg_{m_1} \sg(\tilde g^{-1})\in H(W(k))\leqslant G(W(k))$
has a trivial image in $H(k)$ and thus it belongs to $1_M+pE$.\endproof 
\medskip
{\bf 3.3.5. Example.} Let $c\in{\bf N}$ be such that $g.c.d.(c,r_M-c)=1$. We assume $(M,\vph)$ is the Dieudonn\'e module of a $p$-divisible group $D$ over ${\rm Spec}(k)$ of height $r:=r_M$, dimension $d:=r-c=r_M-c$, and (unique) slope ${d\over {r}}$. Let $m:={\bf T}({\rm End}(M),\vph)$. Let $E$ be the $W(k)$-subalgebra of ${\rm End}(M)$ generated by elements of ${\rm End}(M)$ fixed by $\vph$. As all slopes of $({\rm End}(M[{1\over p}]),\vph)$ are $0$, any $W(k)$-submodule $O$ of ${\rm End}(M)$ with the property that $(O,\vph)$ is a Dieudonn\'e--Fontaine $p$-divisible object over $k$, is $W(k)$-generated by elements fixed by $\vph$ and so is contained in $E$. Thus $m\in{\bf N}\cup\{0\}$ is the smallest number such that $p^m{\rm End}(M)\subseteq E$. As $c\in{\bf N}$, we have $E\neq {\rm End}(M)$ and so $m\ge 1$. The $i$-number $n$ of $(M,\vph)$ is at most $m$, cf. 3.3.4. But $D$ is uniquely determined up to isomorphism by $D[p^n]$ (cf. 3.2.3) and thus also by $D[p^m]$. Using direct sums of $t\in{\bf N}$ copies of $(M,\vph)$, a similar argument shows that $D^t$ is uniquely determined up to isomorphism by $D^t[p^m]$.  
\medskip
{\bf 3.3.6. Example.} 
The case $m=1$ of 3.3.5 also solves positively the isoclinic case of [35, Conjecture 5.7] as one can easily check this starting from [10, 5.3 and 5.4]. For reader's convenience here is the version of the last sentence in the spirit of this paper. In this Subsubsection, we use the notations of 3.3.5 and we moreover assume that $D$ is such that there exists a $W(k)$-basis $\{e_0,e_1,\ldots,e_{r-1}\}$ of $M$ with the property that $\vph(e_i)$ is $e_{i+d}$ if $i\in S(0,c-1)$ and is $pe_{i+d}$ if $i\in S(c,r-1)$. Here and below, all the lower right indices of the form ${}_i$ and ${}_{i,j}$ are mod $r$. For $i,j\in S(0,r-1)$ let $e_{i,j}\in {\rm End}(M)$ be such that it annihilates $e_{j^\prime}$ if $j^\prime\neq j$ and takes $e_j$ into $e_i$. We have $\vph(e_{i,j})=p^{n_{i,j}}e_{i+d,j+d}$, where the integer $n_{i,j}$ is defined by the rule:
\medskip
{\bf (*)} {\it it is $0$ if $(i,j)\in S(0,c-1)^2\cup S(c,r-1)^2$, it is $1$ if $(i,j)\in S(c,r-1)\times S(0,c-1)$, and it is $-1$ if $(i,j)\in S(0,c-1)\times S(c,r-1)$.}
\medskip
We check that $m=1$. Replacing $D$ by $D^{\rm t}$ if necessary, we can assume $c>d=r-c$; so $c>{r\over 2}$. Based on (2) and the inequality $m\ge 1$, to show that $m=1$ it is enough to show that for all pairs $(i,j)\in S(0,r-1)^2$ we have an equality ${\bf S}\tau_{i,j}=1$, where $\tau_{i,j}:=(n_{i,j},n_{i+d,j+d},n_{i+2d,j+2d},\ldots,n_{i+(r-1)d,j+(r-1)d})$. 
\smallskip
So it suffices to show that none of the $\tau_{i,j}$'s is of the form $(-1,0,0,\ldots,0,-1,\ldots)$, cf. definitions 2.2.4 (b) and (c). Thus it suffices to show that for any pair $(i_0,j_0)\in S(0,c-1)\times S(c,r-1)$, the first non-zero number of the sequence $n_{i_0+d,j_0+d},\ldots,n_{i_0+rd,j_0+rd}$ is $1$. We write $j_0=c+d_0$, with $d_0\in S(0,d-1)$. As $d_0\in S(0,c-1)$, we can assume $n_{i_0+d,j_0+d}=n_{i_0+d,d_0}$ is $0$. Thus $(i_1,j_1):=(i_0+d,d_0)\in S(0,c-1)^2$ and we have $i_1\ge j_1$. 
\smallskip
We have to show that the first non-zero number of the sequence $n_{i_1,j_1},n_{i_1+d,j_1+d},\ldots$, $n_{i_1+(r-1)d,j_1+(r-1)d}$ is 1. We can  assume $n_{i_1+d,j_1+d}\neq 1$. As $i_1\ge j_1$ we have $n_{i_1+d,j_1+d}\neq -1$. Thus $n_{i_1+d,j_1+d}=0$. If $i_1+d\le c-1$, let $(i_2,j_2):=(i_1+d,j_1+d)\in S(0,c-1)^2$; if $i_1+d\ge c$, then from (*) and the equality $n_{i_1+d,j_1+d}=0$ we get $j_1+d\ge c$ and thus we have $(i_2,j_2):=(i_1+2d-r,j_1+2d-r)\in S(0,c-1)^2$. We conclude that $(i_2,j_2)\in S(0,c-1)^2$ and $i_2\ge j_2$. We have $n_{i_2,j_2}=0$. We have to show that the first non-zero number of the sequence $n_{i_2,j_2},\ldots,n_{i_2+(r-1)d,j_2+(r-1)d}$ is 1. As in this way we can not construct indefinitely pairs $(i_u,j_u)\in S(0,c-1)^2$ with $i_u\ge j_u$ (here $u\in{\bf N}$), we get that the first non-zero number of the sequence $n_{i_2,j_2},\ldots,n_{i_2+(r-1)d,j_2+(r-1)d}$ is 1. 
\smallskip
So ${\bf S}\tau_{i,j}=1$ and $m=1$. Thus for $t\in {\bf N}$, $D^t$ is uniquely determined up to isomorphism by $D^t[p]$ (cf. end of 3.3.5). This was predicted by [35, Conjecture 5.7].
\medskip
{\bf 3.3.7. Example.} We assume that all slopes of $({\rm End}(M),\vph)$ are $0$ and that $G_1={\bf GL}_M$. Thus 3.3 (i) holds. Let $E_{1{\bf Z}_p}$ be the ${\bf Z}_p$-subalgebra of ${\rm End}(M)$ formed by elements fixed by $\vph$. Let $E_1:=E_{1{\bf Z}_p}\otimes_{{\bf Z}_p} W(k)\subseteq{\rm End}(M)$; we have $v_1=r_M^2$ and $E_2={\rm End}(M)$. We also assume that $p\ge 3$ and that $G=$ {\bf Sp}$(M,\lambda_M)$ (resp. and that $G=$ {\bf SO}$(M,\lambda_M)$), where $\lambda_M$ is a perfect alternating (resp. perfect symmetric) bilinear form on $M$ which is a principal (resp. a principal bilinear) quasi-polarization of $(M,\vph)$. As $p\ge 3$, we have a direct sum decomposition ${\rm End}(M)={\rm Lie}(G)\oplus {\rm Lie}(G)^{\perp}$, where ${\rm Lie}(G)^{\perp}$ is the perpendicular of ${\rm Lie}(G)$ with respect to the trace form on ${\rm End}(M)$. As $\vph$ normalizes ${\rm Lie}(G)[{1\over p}]$, it also normalizes ${\rm Lie}(G)^{\perp}[{1\over p}]$. Thus $E_1$ has a $W(k)$-basis $\Mb_1=\{e_1,\ldots,e_{v_1}\}$ that is the disjoint union of a ${\bf Z}_p$-basis $\{e_1,\ldots,e_v\}$ of $E_{1{\bf Z}_p}\cap {\rm Lie}(G)$ and of a ${\bf Z}_p$-basis $\{e_{v+1},\ldots,e_{v_1}\}$  of $E_{1{\bf Z}_p}\cap {\rm Lie}(G)^{\perp}$. Let $\pi_1:=1_{S(1,v_1)}$. 
\smallskip
Properties 3.3 (ii.a) to (ii.c) hold, cf. constructions. We check that 3.3 (ii.d) holds. Let $g\in G(W(k))$ be of the form $1_M+\sum_{l\in S(1,v_1)} p^tx_le_l$, where all $x_l$'s belong to $W(k)$. The involution of ${\rm End}(M)$ defined by $\lambda_M$ fixes ${\rm Lie}(G)^{\perp}$ and acts as $-1$ on ${\rm Lie}(G)$. Thus the product 
$(1_M-\sum_{l=1}^v p^tx_le_l+\sum_{l=v+1}^{v_1} p^tx_le_l)(1_M+\sum_{l\in S(1,v_1)} p^tx_le_l)$ 
is $1_M$ (as $g\in G(W(k))$) and belongs to $1_M+2(\sum_{l=v+1}^{v_1} p^tx_le_l)+p^{t+1}E_1$. As $p\ge 3$, for $l\in S(v+1,v_1)$ we have $x_l\in pW(k)$. So 3.3 (ii.d) holds. Thus 3.3 (ii) holds. So 3.3.2 applies. In particular, the $i$-number of $(M,\vph,G)$ is at most $m_1+1$, where $m_1\in {\bf N}\cup\{0\}$ is the smallest number such that we have $p^{m_1}{\rm End}(M)\subseteq E_1=E_{1{\bf Z}_p}\otimes_{{\bf Z}_p} W(k)$ (cf. 3.3.2 (b)).  
\bigskip\smallskip
\centerline{\bigsll {\bf \S4 Four examples}}
\bigskip\smallskip
In this section we include four examples that pertain to Subsections 3.1 to 3.3. Let $\vep_p\in\{1,2\}$ be as before 3.1. Let $(M,\vph,G)$ be a $p$-divisible object with a group over $k$. Let $n\in {\bf N}\cup\{0\}$ be the $i$-number of $(M,\vph,G)$, cf. 3.1.4. In this Section we will assume that $k=\bar k$ and that $G$ is a reductive group scheme over ${\rm Spec}(W(k))$. Thus the group scheme $G$ is smooth over ${\rm Spec}(W(k))$ and its fibres are connected and have trivial unipotent radicals. As $G$ is smooth over ${\rm Spec}(W(k))$, with the notations of 2.7 we have $G=G^\prime$ and $n(G)=0$. Let $M=\oplus_{i=a}^b \tilde F^i(M)$, $(F^i(M))_{i\in S(a,b)}$, and $\mu:{\bf G}_m\to G$ be as in 2.2.1 (d). Let $b_L\in S(0,b-a)$ be the smallest number such we have a direct sum decomposition
$${\rm Lie}(G):=\oplus_{i=-b_L}^{b_L} \tilde F^i({\rm Lie}(G))$$
with the property that $\be\in {\bf G}_m(W(k))$ acts on $\tilde F^i({\rm Lie}(G))$ through $\mu$ as the multiplication with $\be^{-i}$. As the group scheme $G$ is reductive, both $W(k)$-modules $\tilde F^{b_L}({\rm Lie}(G))$ and $\tilde F^{-b_L}({\rm Lie}(G))$ are non-zero. As in 2.5, we have a $\sg$-linear automorphism $\sg_0:=\vph\mu(p)$ of $M$. As $\vph=\sg_0\mu(p^{-1})$, the $s$-number of $({\rm Lie}(G),\vph)$ is $b_L$. If $b_L\le 1$, we say $(M,\vph,G)$ is a {\it Shimura $p$-divisible object} over $k$. Let $f_{-1}\in {\bf N}\cup\{0\}$ be the rank of $\tilde F^{-1}({\rm Lie}(G))$. 
\smallskip
In Subsections 4.1 to 4.5 we will consider four unrelated situations.
\medskip
{\bf 4.1. Example 1.} In this Subsection we assume that $b_L=1$ and that all slopes of $({\rm Lie}(G),\vph)$ are $0$. Let $f\in{\bf N}$ be the smallest number such that there exists a filtration 
$$0=\Me_0\subseteq \Me_1\subseteq \cdots \subseteq \Me_f:={\rm Lie}(G)$$ 
by $W(k)$-submodules that are direct summands, with the property that for any number $i\in S(1,f)$, the quotient $W(k)$-module $\Me_i/\Me_{i-1}$ is a maximal direct summand of ${\rm Lie}(G)/\Me_{i-1}$ normalized by $\vph$. For $i\in S(2,f)$ we choose $x_i\in \Me_i\setminus (\Me_{i-1}+p\Me_i)$ such that we have $p\vph(x_i)-px_i\in \Me_{i-1}\setminus p\Me_{i-1}$ and the images of $x_1,\ldots,x_{i-1}$, and $x_i$ in 
$${\rm Lie}(G)/(\tilde F^0({\rm Lie}(G))+\tilde F^1({\rm Lie}(G))+p{\rm Lie}(G))\tilde\to \tilde F^{-1}({\rm Lie}(G))/p\tilde F^{-1}({\rm Lie}(G))$$ 
are $k$-linearly independent. The possibility of making such choices is implied by the maximal property of $\Me_i/\Me_{i-1}$. By reasons of ranks we get $f-1\le f_{-1}$. By induction on $j\in S(1,f)$ we get $\vph$ normalizes $\Me_1+p\Me_2+\cdots +p^{j-1}\Me_j$. Taking $j=f$, we get that $E:=\Me_1+p\Me_2+\cdots +p^{f-1}\Me_f$ has a $W(k)$-basis formed by elements fixed by $\vph$. Thus $(E,\vph)$ is a Dieudonn\'e $p$-divisible object over $k$. As $p^{f-1}$ annihilates ${\rm Lie}(G)/E$ we get ${\bf T}({\rm Lie}(G),\vph)\le f-1$. Thus (cf. 3.1.5 for the first inequality)
$$n\le 2{\bf T}({\rm Lie}(G),\vph)+\vep_p\le 2(f-1)+\vep_p\le 2f_{-1}+\vep_p.$$ 
\indent
Often 3.3.4 provides (resp. 3.3.2 (b) and 3.3.7 applied with $G_1={\bf GL}_M$ provide) better upper bounds of $n$. For instance, if $G={\bf GL}_M$ (resp. if $p\ge 3$ and $G=$ {\bf Sp}$(M,\lambda_M)$ with $\lambda_M$ as a principal quasi-polarization of $(M,\vph)$) we get $n\le {\bf T}({\rm Lie}(G),\vph)\le f_{-1}$ (resp. $n\le {\bf T}({\rm Lie}(G),\vph)+1\le f_{-1}+1$). If $(M,\vph)$ is the Dieudonn\'e module of a supersingular $p$-divisible group over ${\rm Spec}(k)$ of height $2d$ and if $G={\bf GL}_M$ (resp. and if $p\ge 3$ and $G=$ {\bf Sp}$(M,\lambda_M)$), then $f_{-1}$ is $d^2$ (resp. is ${{d^2+d}\over 2}$). Thus we have the following concrete application of 3.2.3 (resp. of 3.2.5):
\medskip
{\bf 4.1.1. Proposition.} {\it Let $d\in {\bf N}$. If $p\ge 2$ (resp. if $p\ge 3$), then any (resp. any principally quasi-polarized) supersingular $p$-divisible group of height $r=2d$ over ${\rm Spec}(k)$ is uniquely determined up to isomorphism by its truncated (resp. its principally quasi-polarized truncated) Barsotti--Tate group of level $d^2$ (resp. of level ${{d^2+d}\over 2}+1$).}
\medskip\smallskip
{\bf 4.2. Root decompositions.} The image of $\mu$ is either trivial or a closed ${\bf G}_m$ subgroup of $G$ and thus its centralizer in $G$ is a reductive group scheme which has a maximal torus (cf. [6, Vol. III, Exp. XIX, 2.8 and 6.1]). Thus there exists a maximal torus $T$ of $G$ through which $\mu$ factors. We have ${\rm Lie}(T)\subseteq \tilde F^0({\rm Lie}(G))$. It is easy to check that there exists $g\in G(W(k))$ such that $g\vph$ normalizes ${\rm Lie}(T)$. Accordingly, for the next three Examples (i.e. until \S5) we will assume that we have $\vph({\rm Lie}(T))={\rm Lie}(T)$. Let 
$${\rm Lie}(G)={\rm Lie}(T)\bigoplus_{\ga\in\Phi} {\got g}_{\ga}$$ 
be the root decomposition relative to $T$. So $\Phi$ is a root system of characters of $T$ and each ${\got g}_{\ga}$ is a free $W(k)$-module of rank $1$ on which $T$ acts via the character $\ga$. 
\smallskip 
Let $\Delta$ be a basis of $\Phi$ such that $\oplus_{\ga\in\Delta}{\got g}_{\ga}\subseteq \oplus_{i=0}^{b_L} \tilde F^i({\rm Lie}(G))$. Let $\Phi^+$ and $\Phi^-$ be the sets of positive and negative (respectively) roots of $\Phi$ with respect to $\Delta$. Let $C$ be the unique Borel subgroup scheme of $G$ which contains $T$ and for which we have ${\rm Lie}(C)={\rm Lie}(T)\bigoplus_{\ga\in\Phi^+} {\got g}_{\ga}$, cf. [6, Vol. III, Exp. XXII,  5.5.1]. As ${\rm Lie}(C)[{1\over p}]$ is generated by the $B(k)$-vector subspace $\oplus_{\ga\in\Delta}{\got g}_{\ga}[{1\over p}]$ of the Lie algebra $\oplus_{i=0}^{b_L} \tilde F^i({\rm Lie}(G))[{1\over p}]$, we have an inclusion ${\rm Lie}(C)\subseteq \oplus_{i=0}^{b_L} \tilde F^i({\rm Lie}(G))$. 
\smallskip
As $\mu$ factors through $T$, for any root $\ga\in \Phi$ there exists an integer $n(\ga)\in S(-b_L,b_L)$ such that we have ${\got g}_{\ga}\subseteq \tilde F^{n(\ga)}({\rm Lie}(G))$. As $\vph=\sg_0\mu(p^{-1})$ and $\sg_0({\rm Lie}(T))=\vph({\rm Lie}(T))={\rm Lie}(T)$, there exists a permutation $\Pi$ of $\Phi$ such that we have
$$\sg_0({\got g}_{\ga})={\got g}_{\Pi(\ga)}\;\;{\rm and}\;\;\vph({\got g}_{\ga})=p^{n(\ga)}{\got g}_{\Pi(\ga)},\;\;\forall\ga\in\Phi.\leqno (21)$$
If $\ga\in\Phi^+$ (resp. if $\ga\in\Phi^-$), then $n(\ga)\in S(0,b_L)$ (resp. then $n(\ga)\in S(-b_L,0)$); this is a consequence of the inclusion ${\rm Lie}(C)\subseteq \oplus_{i=0}^{b_L} \tilde F^i({\rm Lie}(G))$. 
\smallskip
As ${\rm Lie}(T)$ is normalized by $\vph$, it has a $W(k)$-basis formed by elements fixed by $\vph$. Let $\Pi_0:=(\ga_1,\ldots,\ga_l)$ be a cycle of $\Pi$. For $j\in S(1,l)$ let $y_{\ga_j}\in {\got g}_{\ga_j}\setminus\{0\}$ be such that  (cf. (21)) we have $\vph(y_{\ga_j})=p^{m(\ga_j)}y_{\ga_{j+1}}$ (with $\ga_{l+1}:=\ga_1$), where $m(\ga_1),\ldots,m(\ga_l)$ are integers that are either all positive or all negative. Let $\Mb_0:=\{y_{\ga_1},\ldots,y_{\ga_l}\}$. 
\smallskip
Let $E$ be a $W(k)$-submodule of ${\rm Lie}(G)$ that contains ${\rm Lie}(T)$, that satisfies the identity $E[{1\over p}]={\rm Lie}(G)[{1\over p}]$, and that is maximal subject to the property that it has a $W(k)$-basis $\Mb$ which is the union of a ${\bf Z}_p$-basis of $\{x\in {\rm Lie}(T)|\vph(x)=x\}$ and of subsets $\Mb_0$ that are associated as above to some cycle $\Pi_0$ of $\Pi$. Let $\pi$ be the permutation of $\Mb$ which fixes $\Mb\cap {\rm Lie}(T)$ and which for $\ga\in\Phi$ takes the element of ${\got g}_{\ga}\cap\Mb$ into the element of ${\got g}_{\Pi(\ga)}\cap\Mb$. 
\medskip
{\bf 4.3. Example 2.} It is not difficult to check that there exists an element $g\in G(W(k))$ which normalizes $T$ and which has the property that $g\vph$ takes ${\rm Lie}(C)$ into ${\rm Lie}(C)$. Accordingly, in this Subsection we assume that there exists a Borel subgroup $\tilde C$ of $G$ which contains $T$ and which has the property that $\vph({\rm Lie}(\tilde C))\subseteq {\rm Lie}(\tilde C)$. So if we have ${\got g}_{\ga}\subseteq{\rm Lie}(\tilde C)$, then $n(\ga)\ge 0$. Thus we have ${\rm Lie}(\tilde C)\subseteq \oplus_{i=0}^{b_L} \tilde F^i({\rm Lie}(G))$. Thus, not to introduce extra notations, we can assume that $C=\tilde C$; so $\vph({\rm Lie}(C))\subseteq {\rm Lie}(C)$. 
\smallskip
The last inclusion implies that $\Pi$ normalizes both $\Phi^+$ and $\Phi^-$. As we have $n(\ga)\ge 0$ if $\ga\in\Phi^+$ and $n(\ga)\le 0$ if $\ga\in\Phi^-$, for any cycle $\Pi_0=(\ga_1,\ldots,\ga_l)$ of $\Pi$ we can choose the above elements $y_{\ga_j}\in {\got g}_{\ga_j}$ to be generators of ${\got g}_{\ga_j}$.  Thus, due to the maximal property of $E$, we have $E={\rm Lie}(G)$ (i.e. for any $\ga\in\Phi$ the intersection $\Mb\cap {\got g}_{\ga}$ is a $W(k)$-basis of ${\got g}_{\ga}$). Thus ${\bf T}({\rm Lie}(G),\vph)=0$, i.e. $({\rm Lie}(G),\vph)=(E,\vph)$ is a Dieudonn\'e--Fontaine $p$-divisible object over $k$. We have $n\le\vep_p$, cf. 3.1.5. Thus $n\le 1$ if $p\ge 3$ and $n\le 2$ if $p=2$. 
\medskip
{\bf 4.3.1. Proposition.} {\it We recall that $G$ is a reductive group scheme, that $T$ is a maximal torus of $G$ through which $\mu:{\bf G}_m\to G$ factors and whose Lie algebra is normalized by $\vph$, and that we have $\vph({\rm Lie}(C))\subseteq {\rm Lie}(C)$ for some Borel subgroup scheme $C$ of $G$ that contains $T$. Then the $i$-number $n$ of $(M,\vph,G)$ is at most $1$.} 
\medskip
{\it Proof:} We know that $n\le 1$ if $p\ge 3$. Thus we can assume $p=2$ (but as the below arguments work for all primes, we will keep the notation $p$ instead of $2$). Let $\Phi(0):=\Phi\cup\{0\}$, $n(0):=0$, and $G(0):=T$. For $\ga\in\Phi$ let $G(\ga)$ be the unique ${\bf G}_a$ subgroup scheme of $G$ that is normalized by $T$ and such that ${\rm Lie}(G(\ga))={\got g}_{\ga}$, cf. [6, Vol. III, Exp. XXII, 1.1]. If $x\in {\got g}_{\ga}$ is such that $x^p\neq 0$, then the torus $T$ acts on $<x^p>$ via the $p$-th power of the character $\ga$. The reduction $\bar x^p$ mod $p$ of $x^p$ belongs to ${\got g}_{\ga}$ mod $p$, cf. [1, Ch. II, 3.1, 3.5, Lemma 3 of 3.19]. From the last two sentences we get that $\bar x^p=0$. This implies that for each $\ga\in\Phi$ we have a bijection ${\rm exp}_{\ga}:{\got g}_{\ga}\tilde\to G(\ga)(W(k))$ which maps $x\in {\got g}_{\ga}$ into ${\rm exp}_{\ga}(x)=\sum_{i=0}^{\infty} {{x^i}\over {i!}}$ (we emphasize that ${\got g}_{\ga}$ is not included in the domain of the exponential map of 2.6 defined for $O=M$). 
\smallskip
Let $g_1\in G(W(k))$ be congruent mod $p$ to $1_M$. Let $l_1\in{\rm Lie}(G)$ be such that $g_1$ is congruent mod $p^{2}$ to $1_M+pl_1$. We show that there exists $g\in G(W(k))$ congruent mod $p$ to $1_M$ and such that $gg_1\vph g^{-1}\vph^{-1}\in G(W(k))$ is congruent mod $p^{2}$ to $1_M$. 
\smallskip
We take $g$ to be a product $\prod_{\ga\in\Phi(0)} g_1^{\ga}$ (taken in any order), where $g_1^{\ga}\in G(\ga)(W(k))$ is congruent mod $p^{1+\max\{0,-n(\ga)\}}$ to $1_M$. Let $l_1^{0}\in {\rm Lie}(T)$ be such that $g_1^0\in G(0)(W(k))$ is congruent mod $p^2$ to $1_M+pl_1^0$. The element $\vph g_1^0\vph^{-1}=\sg_0 g_1^0\sg_0^{-1}\in G(0)(B(k))\cap {\bf GL}_M(W(k))=G(0)(W(k))$ is congruent mod $p^2$ to $1_M+p\sg_0(l_1^0)$. For $\ga\in \Phi$ let $x_1^{\ga}\in p^{1+\max\{0,-n(\ga)\}}{\got g}_{\ga}$ and $l_1^{\ga}\in {\got g}_{\ga}$ be such that $g_1^{\ga}={\rm exp}_{\ga}(x_1^{\ga})$ is congruent mod $p^{2+\max\{0,-n(\ga)\}}$ to $1_M+p^{1+\max\{0,-n(\ga)\}}l_1^{\ga}$, cf. 2.6.1 and 2.6.2. Based on (21) we have 
$$\vph g_1^{\ga}\vph^{-1}={\rm exp}_{\Pi(\ga)}(p^{1+\max\{0,n(\ga)\}}\sg_0(x_1^{\ga}))\in G(\Pi(\ga))(W(k)).$$ 
So if $n(\ga)>0$ (resp. if $n(\ga)\le 0$), then from 2.6 (b) (resp. from the very definition of $l_1^{\ga}$) we get that $\vph g_1^{\ga}\vph^{-1}={\rm exp}_{\Pi(\ga)}(p^{1+\max\{0,n(\ga)\}}\sg_0(x_1^{\ga}))$ is congruent mod $p^2$ to $1_M$ (resp. to $1_M+p\sg_0(l_1^{\ga})$). Thus by replacing $g_1$ with the following product
$$gg_1\vph g^{-1}\vph^{-1}=(\prod_{\ga\in\Phi(0)} g_1^{\ga})g_1 (\prod_{\ga\in\Phi(0)} \vph g_1^{\ga}\vph^{-1})^{-1}\in G(W(k))$$ 
 of elements congruent mod $p$ to $1_M$, the role of $l_1$ gets replaced by the one of 
$$\tilde l_1:=l_1+\sum_{\ga\in\Phi,n(\ga)>0} l_1^{\ga}+\sum_{\ga\in\Phi(0),n(\ga)=0} [l_1^{\ga}-\sg_0(l_1^{\ga})]+\sum_{\ga\in\Phi,n(\ga)<0} -\sg_0(l_1^{\ga}).$$ 
By writing all elements defining $\tilde l_1$ as linear combinations of elements of the $W(k)$-basis $\Mb$ of $E={\rm Lie}(G)$, as in the part of 3.1.1 that involves (+) and (--) we argue that we can choose the $l_1^{\ga}$'s and so also the $g_1^{\ga}$'s, such that we have $\tilde l_1\in p{\rm Lie}(G)$; here $\ga\in\Phi(0)$. Thus $gg_1\vph g^{-1}\vph^{-1}$ is congruent mod $p^{2}$ to $1_M$. 
So if $n=2$, then the $i$-number of $(M,\vph,G)$ is at most $1$ and this contradicts the definition of $n$. As $n\le 2$, we get $n\le 1$. \endproof
\medskip
{\bf 4.3.2. Proposition.} {\it We continue to work under the hypotheses of 4.3.1. Let $E_0:=\{e\in {\rm Lie}(G)|\vph(e)=e\}$. We have:
\medskip
{\bf (a)} The cocharacter $\mu:{\bf G}_m\to G$ is the unique Hodge cocharacter of $(M,\vph,G)$ that centralizes $E_0$. 
\smallskip
{\bf (b)} The lift $(F^i(M))_{i\in S(a,b)}$ is the unique lift of $(M,\vph,G)$ such that for any $e\in E_0$ and every $i\in S(a,b)$ we have $e(F^i(M))\subseteq F^i(M)$.}
\medskip
{\it Proof:} As $\vph=\sg_0\mu(p^{-1})$, we have $n(\ga)=0$ if and only if ${\got g}_{\ga}\subseteq\tilde F^0({\rm Lie}(G))$. So as $\Pi$ normalizes $\Phi^+$ and $\Phi^-$ and as $n(\ga)\ge 0$ (resp. $n(\ga)\le 0$) if $\ga\in\Phi^+$ (resp. if $\ga\in\Phi^-$), we easily get that $E_0\subseteq \tilde F^0({\rm Lie}(G))$. Thus $\mu$ centralizes $E_0$ and for $e\in E_0$ we have $e(\tilde F^i(M))\subseteq \tilde F^i(M)$ for all $i\in S(a,b)$; thus $e(F^i(M))\subseteq F^i(M)$ for all $i\in S(a,b)$.
\smallskip
Let $\mu_1$ be another Hodge cocharacter of $(M,\vph,G)$ that centralizes $E_0$. As $\vph({\rm Lie}(T))={\rm Lie}(T)$, ${\rm Lie}(T)$ is $W(k)$-generated by elements of $E_0\cap {\rm Lie}(T)$. Thus $\mu_1$ centralizes ${\rm Lie}(T)$ and therefore it factors through $T$. So $\mu$ and $\mu_1$ commute. So to show that $\mu=\mu_1$ it is enough to show that $\mu_k=\mu_{1k}$. As $\vph^{-1}(M)=(\sg_0\mu(p^{-1}))^{-1}(M)=\oplus_{i=a}^b p^{-i}\tilde F^i(M)$, for $i\in S(a,b)$ we have $((p^{-i}\vph)^{-1}(M))\cap M=\sum_{t=0}^{i-a} p^tF^{i-t}(M)$.
This identity implies that the filtration $(F^i(M)/pF^i(M))_{i\in S(a,b)}$ of $M/pM$ is uniquely determined by $(M,\vph)$. Thus both cocharacters $\mu_k$ and $\mu_{1k}$ act on $F^i(M)/F^{i+1}(M)+pF^i(M)\tilde\to \tilde F^i(M)/p\tilde F^i(M)$ via the $-i$-th power of the identity character of ${\bf G}_m$. So as $\mu_k$ and $\mu_{1k}$ commute, by decreasing induction on $i\in S(a,b)$ we get that $\tilde F^i(M)/p\tilde F^i(M)$ is the maximal $k$-vector subspace of $M/pM$ on which both $\mu_k$ and $\mu_{1k}$ act via the $-i$-th power of the identity character of ${\bf G}_m$. This implies $\mu_k=\mu_{1k}$. Thus $\mu=\mu_1$. So (a) holds.
\smallskip
Let $(F_1^i(M))_{i\in S(a,b)}$ be another lift of $(M,\vph,G)$ such that for any $e\in E_0$ and $i\in S(a,b)$, we have $e(F^i_1(M))\subseteq F^i_1(M)$. The inverse of the canonical split cocharacter of $(M,(F_1^i(M))_{i\in S(a,b)},\vph)$ fixes all elements of $E_0$ (cf. the functorial aspects of [42, p. 513]), factors through $G$ (cf. 2.5), and thus it is $\mu$ (cf. (a)). Thus for $i\in S(a,b)$ we have $F^i(M)=\oplus_{j=i}^b\tilde F^j(M)=F^i_1(M)$. So (b) holds.\endproof 
\medskip
If $g\in G(W(k))$, then the Newton polygon of $(M,g\vph)$ is above the Newton polygon of $(M,\vph)$ (cf. [37, Thm. 4.2]). Thus Proposition 4.3.1 generalizes the well known fact that an ordinary $p$-divisible group $D$ over ${\rm Spec}(k)$ is uniquely determined up to isomorphism by $D[p]$. Proposition 4.3.2 generalizes the well known fact that the canonical lift of $D$ is the unique lift of $D$ to ${\rm Spec}(W(k))$ with the property that any endomorphism of $D$ lifts to it. The last two sentences motivate the next definition.
\medskip
{\bf 4.3.3. Definition.} We refer to $(M,\vph,G)$ of 4.3.1 as an {\it ordinary $p$-divisible object with a reductive group} over $k$ and to either $(F^i(M))_{i\in S(a,b)}$ or $(M,(F^i(M))_{i\in S(a,b)},\vph,G)$ as the {\it canonical lift} of $(M,\vph,G)$.
\medskip
{\bf 4.4. Example 3.} Let $c\in {\bf N}$. We assume that $b_L=1$ and that there exists a direct sum decomposition $M=\oplus_{i=1}^c M_i$ in $W(k)$-modules of rank $2$ such that $G=\prod_{i=1}^c {\bf GL}_{M_i}$ and we have $\vph(M_i[{1\over p}])=M_{i+1}[{1\over p}]$ for $i\in S(1,c)$, where $M_{c+1}:=M_1$. We have $r_M=2c$. As $b_L=1$ and $G\tilde\to {\bf GL}_2^c$, we have $f_{-1}\in S(1,c)$. For $i\in S(1,c)$ we have $\vph({\rm End}(M_i)[{1\over p}])={\rm End}(M_{i+1})[{1\over p}]$. Thus the permutation $\pi$ of $\Mb$ has at most two cycles formed by elements of $\Mb\setminus {\rm Lie}(T)$ (equivalently the permutation $\Pi$ of $\Phi$ has at most two cycles). If we have one such cycle, then its length is $2c$. If we have two such cycles, then their length are $c$. 
\smallskip
Let $\vep\in\{1,2\}$. Let $l=\vep c$ be such that we have a cycle $(y_1,\ldots,y_l)$ of $\pi$ formed by elements of $\Mb\setminus {\rm Lie}(T)$. For $i\in S(1,l)$ let $s_i\in{\bf N}\cup\{0\}$ be such that $e_i:=p^{-s_i}y_i$ generates ${\got g}_{\ga_i}$ for some $\ga_i\in\Phi$. Let $(n_1,\ldots,n_l):=(n({\ga_1}),\ldots,n({\ga_l}))\in {\bf Z}^l$; we have $\vph(e_i)=p^{n_i}e_{i+1}$ (cf. (21)). As $b_L=1$, for $i\in S(1,l)$ we have $n_i\in\{-1,0,1\}$. The number $n_+$ (resp. $n_-$) of those $i\in S(1,l)$ such that $n_i=1$ (resp. $n_i=-1$), is at most $f_{-1}$. Moreover $n_++n_-=\vep f_{-1}\le\vep c=l$. Thus $n_0:={\rm min}\{n_+,n_-\}\le {\vep\over 2}f_{-1}$. 
\smallskip
We have ${\bf T}(<e_1,\ldots,e_l>,\vph)\le n_0$ and there exist numbers $a_i\in S(0,n_0)$ such that $(<p^{a_1}e_1,\ldots,p^{a_l}e_l>,\vph)$ is a Dieudonn\'e--Fontaine $p$-divisible object over $k$, cf. (3) and the proof of (2). Based on the maximal property of $E$, we can assume $<p^{a_1}e_1,\ldots,p^{a_l}e_l>\subseteq <y_1,\ldots,y_l>$ (i.e. $s_i\le a_i$ for all $i\in S(1,l)$). If $l=c$, then the set $\{n_+,n_-\}$ (and so also $n_0$) does not depend on the choice of the cycle $(y_1,\ldots,y_l)$. Thus we can choose $E$ such that $p^{n_0}{\rm Lie}(G)\subseteq E$. So any element of $G(W(k))$ congruent mod $p^{n_0+1}$ to $1_M$ belongs to $1_M+pE$. As ${\rm Lie}(T)\subseteq E$, $E$ is a $W(k)$-subalgebra of $\prod_{i=1}^c {\rm End}(M_i)$ and so also of ${\rm End}(M)$. Thus (cf. 3.3.2 (b) applied with $G=G_1$ and with $m_1\in S(0,n_0)$) we have
$$n\le n_0+1={\rm min}\{n_+,n_-\}+1\le {\vep\over 2} f_{-1}+1\le c+1.$$ 
\indent
{\bf 4.5. Example 4.} Let $d\in{\bf N}\setminus\{1,2\}$. Let $D$ be a $p$-divisible group over ${\rm Spec}(k)$ of height $r=2d$, dimension $d$, and slopes ${1\over d}$ and ${{d-1}\over d}$. Let $(M,\vph_0)$ be the Dieudonn\'e module of $D$; we have $r_M=r$. It is easy to see that we have a short exact sequence
$$0\to D_2\to D\to D_1\to 0$$ 
 of $p$-divisible groups over ${\rm Spec}(k)$, where the slopes of $D_1$ and $D_2$ are ${{d-1}\over d}$ and ${1\over d}$ (respectively). This short exact sequence is different from the classical slope filtration of $D$ (see [41, \S3]) which is a short exact sequence $0\to D_1\to D\to D_2\to 0$. As $D_1$ and $D_2$ are uniquely determined up to isomorphisms (see [5, Ch. IV, \S8]), there exist a $W(k)$-basis $\{e_1,\ldots,e_{r}\}$ of $M$ and elements $x_1,\ldots,x_d\in <e_1,\ldots,e_d>$ such that $\vph_0$ takes the $r$-tuple $(e_1,\ldots,e_{r})$ into $(e_2,pe_3, \ldots,pe_d,pe_1,e_{d+2}+x_2,\ldots,e_{r}+x_d,pe_{d+1}+px_1)$. Let 
$$M_1:=<e_1,\ldots,e_d>\;\; {\rm and}\;\;M_2:=<e_{d+1},\ldots,e_{r}>.$$ 
The pairs $(M_1,\vph_0)$ and $(M/M_1,\vph_0)$ are the Dieudonn\'e modules of $D_1$ and $D_2$ (respectively). Let $\vph$ be the $\sg$-linear endomorphism of $M$ that takes $(e_1,\ldots,e_{r})$ into $(e_2,pe_3, \ldots,pe_d,pe_1,e_{d+2},\ldots,e_r,pe_{d+1})$. Let $G:={\bf GL}_M$.
\smallskip
Let $P_1$ and $P_2$ be the maximal parabolic subgroup schemes of $G$ that normalize $M_1$ and $M_2$ (respectively). Let $U_1$ and $U_2$ be the unipotent radicals of $P_1$ and $P_2$ (respectively). Let $u_0\in U_1(W(k))$ be the unique element such that $\vph_0=u_0\vph$. As $T$ we take the maximal torus of $G$ that normalizes $<e_i>$ for all $i\in S(1,r)$. Let $L_1={\bf GL}_{M_1}\times_{{\rm Spec}(W(k))} {\bf GL}_{M_2}$ be the unique Levi subgroup scheme of either $P_1$ or $P_2$ such that $T\leqslant L_1$, cf. [6, Vol. III, Exp. XXVI, 1.12 (ii)]. We have natural identifications ${\rm Lie}(L_1)={\rm End}(M_1)\oplus {\rm End}(M_2)$, ${\rm Lie}(U_1)={\rm Hom}(M_2,M_1)$, and ${\rm Lie}(U_2)={\rm Hom}(M_1,M_2)$. The triples $(M,\vph_0,P_1)$ and $(M,\vph_0,U_1)$ are latticed $F$-isocrystals with a group over $k$. 
\smallskip
As $U_1$ is commutative we have $({\rm Lie}(U_1),\vph)=({\rm Lie}(U_1),\vph_0)$. We easily get that ${\rm Lie}(U_1)$ has a $W(k)$-basis $\{e_i^{(j)}|1\le i,j\le d\}$ such that $\vph_0(e_i^{(j)})=p^{n_i^{(j)}}e_{i+1}^{(j)}$, where each $d$-tuple $(n_1^{(j)},n_2^{(j)},\ldots,n_d^{(j)})$ is either $(1,1,\ldots,1,1,-1)$ or some $d$-tuple of the form $(1,1,\ldots,1,0,1,\ldots,1,0)$. As ${\bf S}(1,1,\ldots,1,-1)=1$ and ${\bf S}(1,1,\ldots,1,0,1,\ldots,1,0)=0$, from (2) applied to all pairs $(<e_1^{(j)},\ldots,e_d^{(j)}>,\vph)$ we get ${\bf T}({\rm Lie}(U_1),\vph_0)\le 1$. As ${\rm Lie}(U_1)$ is a $W(k)$-submodule of ${\rm End}(M)$ whose product with itself is the zero $W(k)$-submodule, the $i$-numbers of $(M,\vph_0,U_1)$ and $(M,\vph,U_1)$ are at most $1$ (cf. 3.3.3). A similar argument shows that both ${\bf T}({\rm Lie}(U_2),\vph)$ and the $i$-number of $(M,\vph,U_2)$ are at most $1$. 
\medskip
{\bf 4.5.1. Proposition.} {\it The $p$-divisible group $D$ is uniquely determined up to isomorphism by $D[p^3]$.}
\medskip
{\it Proof:} Let $t\in{\bf N}\setminus\{1\}$. Let $g_t\in G$ be congruent mod $p^t$ to $1_M$. We will show that if $t\ge 3$, then there exists $g\in G(W(k))$ such that $gg_t\vph_0 g^{-1}\vph_0^{-1}\in G(W(k))$ is congruent mod $p^{t+1}$ to $1_M$. The product morphism $U_2\times_{{\rm Spec}(W(k))} L_1\times_{{\rm Spec}(W(k))} U_1\to G$ is an open embedding around the identity section, cf. [6, Vol. III, Exp. XXII, 4.1.2]. Thus we can write $g_t=u_2l_1u_1$, where the elements $u_1\in U_1(W(k))$, $l_1\in L_1(W(k))$, and $u_2\in U_2(W(k))$ are all congruent mod $p^t$ to $1_M$. As the $i$-number of $(M,\vph_0,U_1)$ is at most 1, to show the existence of $g$ we can replace $u_0$ by any other element of $U_1(W(k))$ that is congruent mod $p$ to $u_0$. Thus we can assume $u_1=1_M$. 
\medskip
For $i\in\{1,2\}$ let $E_i$ be the $W(k)$-span of the ${\bf Z}_p$-algebra of endomorphisms of $(M_i,\vph)$. Let $E:=E_1\oplus E_2$. We have $p{\rm Lie}(L_1)\subseteq E$ (cf. 3.3.6 applied to both $D_1$ and $D_2$) and thus $l_1\in L_1(W(k))\cap (1_M+p^{t-1}E)$. There exists $\tilde l_1\in L_1(W(k))$ congruent mod $p^{t-1}$ to $1_M$ and such that $\tilde l_1l_1\vph\tilde l_1^{-1}=\vph$, cf. 3.3.2 (a) applied with $j=t-1$ to $(M,\vph,L_1)$ and $E$. So, as $L_1$ normalizes both $U_1$ and $U_2$, by replacing $g_t\vph_0$ with $\tilde l_1g_t\vph_0\tilde l_1^{-1}$, $u_2$ with $\tilde l_1u_2 \tilde l_1^{-1}$, and $u_0$ with the element $\tilde l_1l_1u_0l_1^{-1}\tilde l_1^{-1}\in U_1(W(k))$ congruent mod $p^{t-1}$ to $u_0$, we can assume $l_1$ is congruent mod $p^{t+1}$ to $1_M$. This implies that $g_t$ and $u_2$ are congruent mod $p^{t+1}$. 
\smallskip
As ${\bf T}({\rm Lie}(U_2),\vph)\le 1$, from 3.3.2 (a) applied in a way similar to the one of the previous paragraph we deduce the existence of $\tilde u_2\in U_2(W(k))$ congruent mod $p^{t-1}$ to $1_M$ and such that $\tilde u_2u_2\vph\tilde u_2^{-1}=\vph$. As $g_t$ and $u_2$ are congruent mod $p^{t+1}$, the element 
$$g_t^\prime:=\tilde u_2g_t\vph_0\tilde u_2^{-1}\vph_0^{-1}=\tilde u_2g_tu_0\vph\tilde u_2^{-1}\vph^{-1} u_0^{-1}=\tilde u_2g_tu_0u_2^{-1}\tilde u_2^{-1}u_0^{-1}\in G(W(k))$$ 
is congruent mod $p^{t+1}$ to the commutator of $\tilde u_2u_2$ and $u_0$. A simple matrix computation of this commutator shows that for $t\ge 3$ we can write $g_t^\prime=u_2^{\prime}l_2^{\prime}u_1^\prime$, where $u_1^{\prime}\in U_1(W(k))$ and $l_1^{\prime}\in L_1(W(k))$ are congruent mod $p^{t-1}$ to $1_M$ and where $u_2^{\prime}\in U_2(W(k))$ is congruent mod $p^{2t-2}$ and so also mod $p^{t+1}$ to $1_M$ (here is the only place where we need $t\neq 2$). 
\smallskip
Repeating twice the above part that allowed us to assume that $u_1$ and $l_1$ are congruent mod $p^{t+1}$ to $1_M$, we get that for $t\ge 3$ we can assume $g_t^\prime$ is congruent mod $p^{t+1}$ to $u_2^\prime$ and thus also to $1_M$. This ends the argument for the existence of $g$. 
\smallskip
Thus for $t\ge 3$ we have $t+1\neq n$; so $n\le 3$. So the Proposition follows from 3.2.3.\endproof
\medskip
{\bf 4.5.2. Notations for $d=3$.}
Let $d$ be $3$. For $\al\in W(k)$ let $\vph_{\al}$ be the $\sg$-linear endomorphism of $M$ that takes $(e_1, \ldots,e_6)$ into $(e_2,pe_3,pe_1,e_5,e_6+\al e_1,pe_4)$. The slopes of $(M,\vph_{\al})$ are $1\over 3$ and $2\over 3$. For $(i,j)\in\{1,2,3\}\times\{4,5,6\}$ let $U_{ij}$ be the unique ${\bf G}_a$ subgroup scheme of $U_1$ that is normalized by $T$, that fixes $e_{j^\prime}$ for $j^\prime\in\{4,5,6\}\setminus\{j\}$, and that takes $e_j$ into $e_j+<e_i>$. Let $n_{\al}\in U_{16}(W(k))$ be the unique element such that $\vph_{\al}=n_{\al}\vph$.
\medskip
{\bf 4.5.3. Proposition.} {\it Suppose $d=3$. Let $\al_1$, $\al_2\in W(k)\setminus\{0\}$ be such that the ${\bf F}_{p^3}$-vector subspace of $k$ generated by $\al_1$ mod $p$ is different from the ${\bf F}_{p^3}$-vector subspace of $k$ generated by $\al_2$ mod $p$. Then $(M,\vph_{\al_1})$ and $(M,\vph_{\al_2})$ are not isomorphic.}
\medskip
{\it Proof:} 
We show that the assumption that there exists $g_1\in G(W(k))$ which is an isomorphism between $(M,\vph_{\al_1})$ and $(M,\vph_{\al_2})$ leads to a contradiction. As ${\rm Lie}(P_1)[{1\over p}]$ is the maximal direct summand of ${\rm End}(M[{1\over p}])$ normalized by $\vph_{\al}$ and such that all slopes of $({\rm Lie}(P_1)[{1\over p}],\vph_{\al})$ are non-negative, the element $g_1\in G(W(k))$ normalizes ${\rm Lie}(P_1)$ and so also $P_1$. Thus $g_1\in P_1(W(k))$. We write $g_1=u_1l_1$, where $u_1\in U_1(W(k))$ and $l_1\in L_1(W(k))$. We have $u_1l_1n_{\al_1}l_1^{-1}l_1\vph=n_{\al_2}\vph u_1\vph^{-1}\vph l_1$. As $u_1l_1n_{\al_1}l_1^{-1}$ and $n_{\al_2}\vph u_1\vph^{-1}$ belong to $U_1(W(k))$, the actions of $l_1\vph$ and $\vph l_1$ on both $M_1$ and $M/M_1$ are equal. So as $l_1\vph$ and $\vph l_1$ both normalize the direct supplement $M_2[{1\over p}]$ of $M_1[{1\over p}]$ in $M[{1\over p}]$, we get $l_1\vph=\vph l_1$; so we also have $u_1l_1n_{\al_1}l_1^{-1}=n_{\al_2}\vph u_1\vph^{-1}$. As $l_1\vph=\vph l_1$, a simple computation shows that $l_1$ takes $e_1$ (resp. $e_6$) into $a_1e_1+b_1e_2+c_1e_3$ (resp. $pa_6e_4+pb_6e_5+c_6e_6$), where $a_1$, $b_1$, $c_1\in W({\bf F}_{p^3})$ (resp. $a_6$, $b_6\in W({\bf F}_{p^3})$ and $c_6\in {\bf G}_m(W({\bf F}_{p^3}))$). As $U_1=\prod_{i=1}^3\prod_{j=4}^6 U_{ij}$, there exist unique elements $u_{ij}\in U_{ij}(W(k))$ such that we have $u_1=\prod_{i=1}^3\prod_{j=4}^6 u_{ij}$. We call $u_{ij}$ as the component of $u_1$ in $U_{ij}(W(k))$. Both $g_1$ mod $p$ and $l_1$ mod $p$ normalize the kernel of $\vph_{\al}$'s mod $p$, i.e. they normalize $<e_2,e_3,e_6>$ mod $p$. So $u_1$ mod $p$ also normalizes $<e_2,e_3,e_6>$ mod $p$. If $(i,j)\neq (1,6)$, then $u_{ij}$ mod $p$ normalizes $<e_2,e_3,e_6>$ mod $p$. Thus $u_{16}$ mod $p$ normalizes $<e_2,e_3,e_6>$ mod $p$ and therefore it is $1_{M/pM}$. As $U_{35}$ fixes both $<e_2,e_3,e_6>$ and $M/<e_2,e_3,e_6>$, $\vph(u_{35})$ mod $p$ is $1_{M/pM}$. Thus the component of $n_{\al_2}\vph u_1\vph^{-1}$ in $U_{16}(W(k))$ is congruent mod $p$ to $n_{\al_2}$. The component of $l_1n_{\al_1}l_1^{-1}$ in $U_{16}(W(k))$ is $n_{\al_3}$, where $\al_3:=a_1\al_1c_6^{-1}\in W(k)$. So $\al_3$ mod $p$ belongs to the ${\bf F}_{p^3}$-vector subspace of $k$ generated by $\al_1$ mod $p$. The component of $u_1l_1n_{\al_1}l_1^{-1}$ in $U_{16}(W(k))$ is congruent mod $p$ to $n_{\al_3}$. As $u_1l_1n_{\al_1}l_1^{-1}=n_{\al_2}\vph u_1\vph^{-1}$, we get that $n_{\al_2}$ and $n_{\al_3}$ are congruent mod $p$. So $\al_2$ mod $p$ belongs to the ${\bf F}_{p^3}$-vector subspace of $k$ generated by $\al_1$ mod $p$. This contradicts our hypothesis. Thus $g_1$ does not exist. \endproof 
\medskip
{\bf 4.5.4. Remark.} The set of isomorphism classes of $p$-divisible groups over ${\rm Spec}(k)$ of height $6$ and dimension $3$ has the same cardinality as $k$, cf. 4.5.3, 1.3, and the classical Dieudonn\'e theory. But the set of isomorphism classes of $p$-torsion subgroup schemes of such $p$-divisible groups over ${\rm Spec}(k)$ is finite (see [26]; to be compared with [34, \S1]). Based on this, 4.5.1, and 3.2.3, we get that the set of those elements $\al\in{\bf G}_m(W(k))$ for which the $i$-number of $(M,\vph_{\al})$ is either 2 or 3, has the same cardinality as $k$. Let $\lambda_M$ be the perfect alternating form on $M$ defined by the rule: if $1\le i<j\le 6$, then $\lambda_M(e_i,e_j)\in\{0,1\}$ and we have $\lambda_M(e_i,e_j)=1$ if and only if $(i,j)\in\{(1,6),(3,5),(2,4)\}$. The form $\lambda_M$ is a principal quasi-polarization of any $(M,\vph_{\al})$. 
\smallskip
From 4.1.1 and 4.5.1 we easily get that $T(6,3)\le 9$. 
\vfill\eject
\centerline{\bigsll {\bf \S 5 Four direct applications}}
\bigskip\smallskip
In this Section we continue to assume that $k=\bar k$.
\medskip\smallskip
{\bf 5.1. The homomorphism form of 1.2.} Let $(M_1,\vph_1)$ and $(M_2,\vph_2)$ be two $F$-crystals over $k$.  Let $(M,\vph):=(M_1,\vph_1)\oplus (M_2,\vph_2)$. For $i\in\{1,2\}$, let $h_i$ be the $h$-number of $(M_i,\vph_i)$. Let $h_{12}:=\max\{h_1,h_2\}$, let $m_{12}:={\bf T}({\rm End}(M),\vph)$, let $v_{12}:=m_{12}+h_{12}$, let $\vep_p\in\{1,2\}$ be as before 3.1, and let $n_{12}:=m_{12}+\vep_p$. The $h$-number of $(M,\vph)$ is $h_{12}$.
\medskip
{\bf 5.1.1. Theorem.} {\it We endow $({\bf N}\cup\{0\})^2$ with the lexicographic order. We have:
\medskip
{\bf (a)} For all $t\in{\bf N}\cup\{0\}$, the images defined via restrictions of the two groups ${\rm Hom}((M_1,\vph_1),(M_2,\vph_2))$ and 
${\rm Hom}((M_1/p^{n_{12}+v_{12}+t}M_1,\vph_1),(M_2/p^{n_{12}+v_{12}+t}M_2,\vph_2))$ in the group ${\rm Hom}((M_1/p^{n_{12}+t}M_1,\vph_1),(M_2/p^{n_{12}+t}M_2,\vph_2))$, coincide.
\smallskip
{\bf (b)} We fix a quadruple $(r_1,r_2,c_1,c_2)\in ({\bf N}\cup\{0\})^4$. There exists a smallest pair 
$$(v,n):=(v(r_1,r_2,c_1,c_2),n(r_1,r_2,c_1,c_2))\in ({\bf N}\cup\{0\})^2$$ 
with the properties that $n\in S(0,v+\vep_p-{\rm max}\{c_1,c_2\})$ and that for any $t\in{\bf N}\cup\{0\}$ and for every two $F$-crystals $(M_1,\vph_1)$ and $(M_2,\vph_2)$ over $k$ which satisfy $(r_{M_1},r_{M_2},h_1,h_2)=(r_1,r_2,c_1,c_2)$, the images defined via restrictions of the two groups ${\rm Hom}((M_1,\vph_1),(M_2,\vph_2))$ and ${\rm Hom}((M_1/p^{n+v+t}M_1,\vph_1),(M_2/p^{n+v+t}M_2,\vph_2))$ in the group ${\rm Hom}((M_1/p^{n+t}M_1,\vph_1),\break (M_2/p^{n+t}M_2,\vph_2))$, coincide. The number $v$ (and so also $n$) has upper bounds that depend only on $r_1$, $r_2$, and ${\rm max}\{c_1,c_2\}$.
\smallskip
{\bf (c)} Let $r_1$, $r_2\in{\bf N}$. There exists a smallest pair 
$$(v,n):=(v(r_1,r_2),n(r_1,r_2))\in ({\bf N}\cup\{0\})^2$$ 
with the properties that $n\in S(0,v+\vep_p-1)$ and that for any $t\in{\bf N}\cup\{0\}$ and for every two $p$-divisible groups $D_1$ and $D_2$ over ${\rm Spec}(k)$ of heights $r_1$ and $r_2$, a homomorphism $D_1[p^{n+t}]\to D_2[p^{n+t}]$ lifts to a homomorphism $D_1\to D_2$ if and only if it lifts to a homomorphism $D_1[p^{n+v+t}]\to D_2[p^{n+v+t}]$. The number $v$  (and so also $n$) has upper bounds that depend only on $r_1$ and $r_2$.}
\medskip\smallskip
{\it Proof:}  Let $e_{12}\in {\rm Hom}((M_1/p^{n_{12}+v_{12}+t}M_1,\vph_1),(M_2/p^{n_{12}+v_{12}+t}M_2,\vph_2))$. Let $g\in {\rm Aut}(M/p^{n_{12}+v_{12}+t}M,\vph)$ be such that it takes $x_2\in M_2/p^{n_{12}+v_{12}+t}M_2$ into $x_2$ and it takes $x_1\in M_1/p^{n_{12}+v_{12}+t}M_1$ into $x_1+e_{12}(x_1)$. Let $\tilde g\in {\rm Aut}(M,\vph)$ be such that it lifts the reduction mod $p^{n_{12}+t}$ of $g$, cf. 3.2.8 (b). Let $\tilde e_{12}:M_1\to M_2$ be the $W(k)$-linear map such that we have $\tilde g(x_1)-\tilde e_{12}(x_1)\in M_1$ for all $x_1\in M_1$. We have $\tilde e_{12}\in {\rm Hom}((M_1,\vph_1),(M_2,\vph_2))$ and moreover $e_{12}$ and $\tilde e_{12}$ have the same image in ${\rm Hom}((M_1/p^{n_{12}+t}M_1,\vph_1),(M_2/p^{n_{12}+t}M_2,\vph_2))$. This proves (a). 
\smallskip
We know that $m_{12}$ has an upper bound $b_{12}\in {\bf N}$ which is effectively computable in terms of $2h_{12}=2{\rm max}\{c_1,c_2\}$ and $r_M^2=(r_{M_1}+r_{M_2})^2=(r_1+r_2)^2$ (cf. 2.4.1 and end of 2.2.1 (e)) and that we have $v\le m_{12}+h_{12}$ (cf. (a)). From this (b) follows. 
\smallskip
To prove (c),  for $i\in\{1,2\}$ we take $(M_i,\vph_i)$ to be the Dieudonn\'e module of $D_i$. We have $h_{12}\in\{0,1\}$. If $h_{12}=0$, then $h_1=h_2=0$ and so both $D_1$ and $D_2$ are \'etale $p$-divisible groups; in such a  case any homomorphism $D_1[p^{n+t}]\to D_2[p^{n+t}]$ lifts to a homomorphism $D_1\to D_2$. Thus we can assume that $h_{12}=1$. Based on 3.2.3, the proofs of (a) and (b) can be adapted to the context of $p$-divisible groups; thus (c) follows from the particular case of (b) when $\max\{c_1,c_2\}$ is $h_{12}=1$.\endproof
\medskip
{\bf 5.1.2. Remark.} If $v_1,n_1\in {\bf N}$ are such that $v_1\ge v$ and $n_1\ge n$, then the homomorphisms parts of 5.1.1 (b) and (c) continue to hold if we replace $(v,n)$ with $(v_1,n_1)$. 
\medskip\smallskip
{\bf 5.2. Transcendental degrees of definition.}  
For simplicity, in this Subsection we work in a context without principal bilinear quasi-polarizations (but we emphasize that all of 5.2 can be adapted to the context of 2.2.1 (c)). Let $(M,\vph,G,(t_{\al})_{\al\in\Mj})$ be a latticed $F$-isocrystal with a group and an emphasized family of tensors over $k$ such that the $W$-condition holds for $(M,\vph,G)$. Let $(F^i(M))_{i\in S(a,b)}$ be a lift of $(M,\vph,G)$. 
\smallskip
Let $\mu:=\mu_{\rm can}:{\bf G}_m\to G$ be the inverse of the canonical split cocharacter of $(M,(F^i(M))_{i\in S(a,b)},\vph)$ (see 2.5). Let $(M_{{\bf Z}_p},G_{{\bf Z}_p})$ be the ${\bf Z}_p$ structure of $(M,G)$ that is defined as in 2.5 by the $\sg$-linear automorphism $\sg_0:=\vph\mu(p):M\tilde\to M$; thus $\sg_0$ fixes $M_{{\bf Z}_p}$ and normalizes $G(W(k))$ and moreover we have $t_{\al}\in\Mt(M_{{\bf Z}_p})$ for all $\al\in\Mj$. Let $n\in {\bf N}\cup\{0\}$ be the $i$-number of $(M,\vph,G)$. 
\medskip
{\bf 5.2.1. Fact.} {\it Let $k_1$ be the smallest subfield of $k$ with the property that $\mu_{W_n(k)}$ is the pull back of a cocharacter $\mu_{1W_n(k_1)}:{\bf G}_m\to G_{{\bf Z}_p}\times_{{\rm Spec}({\bf Z}_p)} {\rm Spec}(W_n(k_1))$. Then the field $k_1$ is finitely generated and its transcendental degree $t(k_1)$ is at most $nr_M^2$. Thus if $t(k_1)=0$, then $k_1$ is a finite field.}
\medskip
{\it Proof:} Let $\Mb$ be a $W(k)$-basis of $M$ such that $\mu$ normalizes the $W(k)$-spans of elements of $\Mb$. Let $\Mb_1$ be a ${\bf Z}_p$-basis of $M_{{\bf Z}_p}$; we also view it as a $W(k)$-basis of $M$. Let $B\in M_{r_M\times r_M}(W(k))$ be the change of coordinates matrix that changes $\Mb_1$-coordinates into $\Mb$-coordinates. Let $R_1$ be the ${\bf F}_p$-subalgebra of $k$ generated by the coordinates of the Witt vectors of length $n$ with coefficients in $k$ that are entries of $B$ mod $p^n$. Obviously $k_1$ is a subfield of the field of fractions of $R_1$. As $R_1$ is generated by $nr_M^2$ elements,  $k_1$ is finitely generated and we have $\dim({\rm Spec}(R_1))\le nr_M^2$. Thus $t(k_1)\le nr_M^2$.\endproof 
\medskip
Until 5.3 we take $G$ to be smooth over ${\rm Spec}(W(k))$. Thus for any $l\in{\bf N}$ there exists a cocharacter $\mu_{1,n+l}$ of $G_{{\bf Z}_p}\times_{{\rm Spec}({\bf Z}_p)} {\rm Spec}(W_{n+l}(k_1^{\rm perf}))$ that lifts the pull back $\mu_{1W_n(k_1^{\rm perf})}$ to ${\rm Spec}(W_n(k_1^{\rm perf}))$ of $\mu_{1W_n(k_1)}$, cf. [6, Vol. II, Exp. IX, 3.6]. We can assume that $\mu_{1,n+l+1}$ lifts $\mu_{1,n+l}$, cf. loc. cit. From [6, Vol. II, Exp. IX, 7.1] we get that there exists a unique cocharacter $\mu_1:{\bf G}_m\to G_{{\bf Z}_p}\times_{{\rm Spec}({\bf Z}_p)} {\rm Spec}(W(k_1^{\rm perf}))$ that lifts all $\mu_{1,n+l}$'s. 
\smallskip
The cocharacter $\mu_{1W(k)}:{\bf G}_m\to G$ is of the form $g_n{\mu}g_n^{-1}$ for some $g_n\in G(W(k))$ congruent mod $p^{n}$ to $1_M$, cf.  [6, Vol. II, Exp. IX, 3.6] and the fact that we have $G(W(k))={\rm proj.}\,{\rm lim.}\,_{l\in {\bf N}}\, G(W_{n+l}(k))$. Let $\tilde g_n:=g_n^{-1}\sg_0g_n\sg_0^{-1}\in G(W(k))$; it is congruent mod $p^n$ to $1_M$. The element $g_n^{-1}\in G(W(k))$ defines an isomorphism between $(M,\sg_0\mu_{1W(k)}({1\over p}),G,(t_{\al})_{\al\in\Mj})$ and $(M,\tilde g_n\vph,G,(t_{\al})_{\al\in\Mj})$. Moreover $(M,\vph,G,(t_{\al})_{\al\in\Mj})$ and $(M,\tilde g_n\vph,G,(t_{\al})_{\al\in\Mj})$ are isomorphic under an isomorphism defined by an element of $G(W(k))$, cf. the very definition of $n$. We conclude that:
\medskip
{\bf (*)}  {\it the quadruple $(M,\vph,G,(t_{\al})_{\al\in\Mj})$ is isomorphic to $(M,\sg_0\mu_{1W(k)}({1\over p}),G,(t_{\al})_{\al\in\Mj})$; thus $(M,\vph,G,(t_{\al})_{\al\in\Mj})$ is definable over $k_1^{\rm perf}$ and moreover its isomorphism class is uniquely determined by the triple $(M_{{\bf Z}_p},(t_{\al})_{\al\in\Mj},\mu_{1W_n(k_1)})$.} 
\medskip
This motivates the following definitions.
\medskip
{\bf 5.2.2. Definitions.} {\bf (a)} We say $k_1$ (resp. $t(k_1)$) is the field (resp. the transcendental degree) of definition of $(M,\vph,G)$ or of $(M,\vph,G,(t_{\al})_{\al\in\Mj})$ with respect to the lift $(F^i(M))_{i\in S(a,b)}$ of $(M,\vph,G)$. By the {\it transcendental degree of definition} $td\in{\bf N}\cup\{0\}$ 
of $(M,\vph,G)$ or of $(M,\vph,G,(t_{\al})_{\al\in\Mj})$ we mean the smallest number we get by considering transcendental degrees of definition of $(M,\vph,G)$ with respect to (arbitrary) lifts of it. 
\smallskip
{\bf (b)} If $td=0$, then by the {\it field of definition} of $(M,\vph,G)$ or of $(M,\vph,G,(t_{\al})_{\al\in\Mj})$ we mean the finite field that has the smallest number of elements and that is the field of definition of $(M,\vph,G)$ with respect to some lift of it.
\medskip
If $td=0$ we do not stop to study when the field of definition of $(M,\vph,G)$ is contained in all fields of definition of $(M,\vph,G)$ with respect to (arbitrary) lifts of $(M,\vph,G)$. 
\medskip
{\bf 5.2.3. Theorem (Atlas Principle).} {\it We recall that $k=\bar k$, that $G$ is smooth over ${\rm Spec}(W(k))$, and that the $W$-condition holds for $(M,\vph,G)$. Let $q\in{\bf N}$. Let $\Mi(q)\in{\bf N}\cup\{0,\infty\}$ be the number of isomorphism classes of latticed $F$-isocrystals with a group and an emphasized family of tensors over $k$ that have the form $(M,g\vph,G,(t_{\al})_{\al\in\Mj})$ for some $g\in G(W(k))$, that have transcendental degrees of definition $0$, and that have ${\bf F}_{p^q}$ as their fields of definition. Then we have $\Mi(q)\in{\bf N}\cup\{0\}$.}
\medskip
{\it Proof:} We check that the number $N_{\rm tor}\in{\bf N}\cup\{\infty\}$ of isomorphism classes of pairs of the form $(M_{{\bf Z}_p},(t_{\al})_{\al\in\Mj})$ that are obtained by replacing $\vph$ with some $g\vph$ and by considering some lift of $(M,g\vph,G)$, is finite. The ``difference" between two such pairs is measured by a torsor $\Th$ of $G_{{\bf Z}_p}$ in the flat topology of ${\rm Spec}({\bf Z}_p)$. So $\Th$ is smooth over ${\rm Spec}({\bf Z}_p)$. Thus $\Th$ is a trivial torsor if and only if $\Th_{{\bf F}_p}$ is a trivial torsor. As the set $H^1({\bf F}_p,G_{{\bf F}_p})$ is finite (cf. [38, Ch. III, \S4, 4.2 and 4.3]), the number of isomorphism classes of torsors of $G_{{\bf F}_p}$ is finite. From the last two sentences we get that $N_{\rm tor}\in{\bf N}$.
\smallskip
Let $m\in{\bf N}$. We check that the number $N(\mu,m,q,G_{{\bf Z}_p})\in{\bf N}\cup\{0,\infty\}$ of cocharacters of $G_{{\bf Z}_p}\times_{{\rm Spec}({\bf Z}_p)} {\rm Spec}(W_m({\bf F}_{p^q}))$ that over ${\rm Spec}(W_m(k))$ are $G(W_m(k))$-conjugate to $\mu_{W_m(k)}$, is also finite. Based on the infinitesimal liftings of [6, Vol. II, Exp. IX, 3.6] and on the fact that the group $G_{{\bf Z}_p}(W_m({\bf F}_{p^q}))$ is finite, it is enough to prove that the number $N(\mu,1,q,G_{{\bf Z}_p})$ is finite. It suffices to prove that $N(\mu,1,q,{\bf GL}_{M_{{\bf Z}_p}})\in {\bf N}$. The number of maximal split tori of ${\bf GL}_{M_{{\bf Z}_p}\otimes_{{\bf Z}_p} {\bf F}_{p^q}}$ is finite and each such torus has precisely $(b-a+1)^{r_M}$ cocharacters that act on $M_{{\bf Z}_p}\otimes_{{\bf Z}_p} {\bf F}_{p^q}$ via those $-i$-th powers of the identity character of ${\bf G}_m$ that satisfy $i\in S(a,b)$. Thus $N(\mu,1,q,{\bf GL}_{M_{{\bf Z}_p}})\in {\bf N}$. 
\smallskip
Let $n_{\rm fam}$ be as in 3.1.5; we have $n_{\rm fam}\ge n$. Based on 5.2 (*), we get that $\Mi(q)$ is bounded from above by a sum of $N_{\rm tor}$ numbers of the form $N(\mu,n_{\rm fam},q,G_{{\bf Z}_p})$ and in particular that $\Mi(q)\le N_{\rm tor}N(\mu,n_{\rm fam},q,{\bf GL}_{M_{{\bf Z}_p}})$. Thus $\Mi(q)\in{\bf N}\cup\{0\}$.\endproof
\medskip\smallskip
{\bf 5.3. Groupoids and stratifications.} Main Theorem A has many reformulations in terms of (stacks of) groupoids. Not to increase the length of the paper, we postpone to future work the introduction of Shimura (stacks of) groupoids that parametrize isomorphism classes of Shimura $p$-divisible objects we defined in the beginning of \S4. Presently, this Shimura context is the most general context to which we can extend the classical deformation theories of $p$-divisible groups (see [30, Ch. 4 and 5], [21, Ch. 3 and 4], [11, 7.1], [8, Main Thm. of Introd.], [9, Main Thm. 1], and [12, \S7]). Here, as an anticipation of the numerous possibilities offered by 1.2, we work only with (principally quasi-polarized) $p$-divisible groups. However, we point out that based on 3.2.2, by using [40, 5.4] as a substitute for the deformation theory of [21, 4.8], the below proof of 5.3.1 can be adapted to contexts that involve arbitrary Shimura varieties of Hodge type and thus involve (principally quasi-polarized)  Dieudonn\'e modules equipped with smooth groups as in 3.2 (like the context of [40, \S5]).
\smallskip 
Let $S$ be a reduced ${\rm Spec}({\bf F}_p)$-scheme. Let $\Md$ be a $p$-divisible group over $S$ of height $r$ and relative dimension $d$. For $i\in\{1,2\}$ let $\Md_i$ be the pull back of $\Md$ to $S_{12}:=S\times_{{\rm Spec}({\bf F}_p)} S$ via the $i$-th projection $p_i:S_{12}\to S$. For $l\in{\bf N}$, let $I_l$ be the affine $S_{12}$-scheme that parametrizes isomorphisms between $\Md_1[p^l]$ and $\Md_2[p^l]$; it is of finite presentation. The morphism $i_l:I_l\to S_{12}$ of $S$-schemes is a ${\rm Spec}({\bf F}_p)$-groupoid that acts on $S$ in the sense of [6, Vol. I, Exp. V] and [31, Appendix A].
\medskip
{\bf 5.3.1. Basic Theorem.} {\it {\bf (a)} There exists a number $l\in{\bf N}$ effectively bounded from above only in terms of $r$ and such that for any algebraically closed field $K$ of characteristic $p$, the pulls back of $\Md$ through two $K$-valued points $y_1$ and $y_2$ of $S$ are isomorphic if and only if the $K$-valued point of $S_{12}$ defined by the pair $(y_1,y_2)$ factors through $i_l:I_l\to S_{12}$.
\smallskip
{\bf (b)} Suppose $S$ is smooth over ${\rm Spec}({\bf F}_p)$ of dimension $d(r-d)$ and $\Md$ is a versal deformation at each maximal point of $S$. Then there exists a stratification $\Ms(\Md)$ of $S$ in reduced, locally closed subschemes such that two points $y_1$ and $y_2$ as in (a) factor through the same stratum if and only if $y_1^*(\Md)$ is isomorphic to $y_2^*(\Md)$. The strata of $\Ms(\Md)$ are regular and equidimensional.
\smallskip
{\bf (c)}  The stratification $\Ms(\Md)$ of (b) satisfies the purity property. 
\smallskip
{\bf (d)} Let $q\in{\bf N}$ and let $K$ be as in (a). Then for any stratum $S_0$ of the stratification $\Ms(\Md)$ of (b) that is a subscheme of $S_K$, there exists a regular scheme $S_0[q]$ that is finite, flat over $S_0$ and such that the pull back of $\Md[p^q]$ to $S_0[q]$ is constant, i.e. is the pull back to $S_0[q]$ of a truncated Barsotti--Tate group of level $q$ over ${\rm Spec}(K)$.} 
\medskip
{\it Proof:} Part (a) follows from 1.3: as $l$ we can take any integer greater that $T(r,d)$. For the rest of the proof we will assume that $S$ is smooth over ${\rm Spec}({\bf F}_p)$ of dimension $d(r-d)$ and that $\Md$ is a versal deformation at each maximal point of $S$. We first construct the strata of $\Ms(\Md)$ that are subschemes of $S_{K}$. Let $y_1\in S(K)=S_K(K)$. Let 
$$S(y_1)^{\rm top}:=(p_2\circ i_l)_{K}((p_1\circ i_l)_{K}^{-1}(y_1))\subseteq S_K^{\rm top}.$$ 
As $I_l$ is a ${\rm Spec}({\bf F}_p)$-scheme of finite type, $S(y_1)^{\rm top}$ is a constructible subset of $S_{K}^{\rm top}$ (cf. [18, Ch. IV, (1.8.4) and (1.8.5)]). Let $\bar S(y_1)$ be the Zariski closure of $S(y_1)^{\rm top}$ in $S_{K}$; it is a reduced, closed subscheme of $S_{K}$. We identify a maximal point of $S(y_1)^{\rm top}$ with a $K$-valued point of $S_K$. We say $S(y_1)^{\rm top}$ is regular at a maximal point of it, if there exists a regular, open subscheme of $\bar S(y_1)$ which contains this point and whose topological space is contained in $S(y_1)^{\rm top}$. As $S(y_1)^{\rm top}$ is a dense, constructible subset of $\bar S(y_1)^{\rm top}$, there exists a regular, open, dense subscheme $W(y_1)$ of $\bar S(y_1)$ such that $W(y_1)^{\rm top}\subseteq S(y_1)^{\rm top}$. 
\smallskip
A point $y_3\in S(K)=S_K(K)$ belongs to $S(y_1)^{\rm top}$ if and only if $y_1^*(\Md)$ and $y_3^*(\Md)$ are isomorphic. Let now $y_2\in S(K)=S_K(K)$ be a maximal point of $S(y_1)^{\rm top}$. Let $i_{12}:y_2^*(\Md)\to y_1^*(\Md)$ be an isomorphism, cf. (a). We have $S(y_2)^{\rm top}=S(y_1)^{\rm top}$. For $i\in\{1,2\}$ let $I_{y_i}$ be the spectrum of the completion of the local ring of $S_{K}$ at $y_i$. We denote also by $y_i$ the factorization of $y_i$ through $I_{y_i}$. Let $I_{y_i}\times_{S_K} S(y_1)^{\rm top}$ be the pull back of $S(y_1)^{\rm top}$ to a constructible subset of $I_{y_i}^{\rm top}$. Due to the versal property of $\Md$ and the fact that $S$ is smooth over ${\rm Spec}({\bf F}_p)$ of dimension $d(r-d)$, the local schemes $I_{y_1}$ and $I_{y_2}$ have dimension $d(r-d)$ and moreover there exists a unique isomorphism 
$$I_{12}:I_{y_1}\tilde\to I_{y_2}$$ 
such that the following two things hold (cf. [21, 4.8]; see [9, 2.4.4] for the equivalence of the categories of $p$-divisible groups over $I_i$ and over the formal completion of $I_i$ along $y_i$):
\medskip
{\bf (i)} we have $I_{12}\circ y_1=y_2:{\rm Spec}(K)\to I_{y_2}$;
\smallskip
{\bf (ii)} there exists an isomorphism $I_{12}^*(\Md_{I_{y_2}})\tilde\to \Md_{I_{y_1}}$ that lifts $i_{12}$.
\medskip
Due to (ii) the local geometries of $S(y_1)^{\rm top}$ at $y_1$ and $y_2$ are the same. In other words, $I_{12}$ induces via restriction an isomorphism 
$J_{12}^{\rm top}:I_{y_1}\times_{S_K} S(y_1)^{\rm top}\tilde\to I_{y_2}\times_{S_K} S(y_1)^{\rm top}$
between constructible subsets. Any commutative ${\bf F}_p$-algebra of finite type is excellent, cf. [29, (34.A) and (34.B)]. So the morphism $I_{y_i}\to S_K$ is regular. From the last two sentences we get (cf. [29, (33)] for the regular part of (iv)) that:
\medskip
{\bf (iii)} the dimensions of $S(y_1)^{\rm top}$ at $y_1$ and $y_2$ are the same;
\smallskip
{\bf (iv)} if $S(y_1)^{\rm top}$ is regular at $y_2$, then $S(y_1)^{\rm top}$ is also regular at $y_1$.
\medskip
By taking $y_2$ to be a maximal point of $W(y_1)^{\rm top}$ and $y_1$ to be an arbitrary maximal point of $S(y_1)^{\rm top}$, from (iii) we get that $S(y_1)^{\rm top}$ is equidimensional and from (iv) we get that $S(y_1)^{\rm top}$ is regular at all its maximal points. So $\bar S(y_1)$ is also equidimensional. As $S(y_1)^{\rm top}$ is a constructible subset of $\bar S(y_1)^{\rm top}$, from the last two sentences we get that $S(y_1)^{\rm top}$ is the underlying topological space of an equidimensional, regular, open subscheme $S(y_1)$ of $\bar S(y_1)$. Thus $S(y_1)$ is a reduced, locally closed subscheme of $S_{K}$. Let 
$$J_{12}:I_{y_1}\times_{S_K} S(y_1)\tilde\to I_{y_2}\times_{S_K} S(y_1)$$ 
be the isomorphism of reduced schemes defined by $I_{12}$ (or $J_{12}^{\rm top}$). 
\smallskip
Let $\Ms_K(\Md)$ be the set of reduced, locally closed subschemes of $S_K$ that are of the form $S(y_1)$ for some $y_1\in S(K)=S_K(K)$. Standard Galois descent shows that there exists a set $\Ms_{{\bf F}_p}(\Md)$ of reduced, locally closed subschemes of $S$ whose pulls back to $S_{\overline{{\bf F}_p}}$ are the elements of $\Ms_{\overline{{\bf F}_p}}(\Md)$. If $L$ is an algebraically closed field that contains $K$ and if $y_1^L$ is the $L$-valued point of $S$ defined by $y_1$, then by very definitions $S(y_1^L)=S(y_1)_L$. So we have natural pull back injective maps $\Ms_K(\Md)\hookrightarrow\Ms_L(\Md)$ and $\Ms_{{\bf F}_p}(\Md)\hookrightarrow\Ms_L(\Md)$. So $\Ms_{{\bf F}_p}(\Md)$ and $\Ms_K(\Md)$'s define a stratification $\Ms(\Md)$ of $S$ in the sense of 2.1.1. So as each $S(y_1)$ is a regular and equidimensional ${\rm Spec}(K)$-scheme, (b) holds.
\smallskip
Let $v$, $n\in{\bf N}$ be as in 5.1.1 (c) for $r_1=r_2:=r$. Let $\tilde n:=\max\{q,l,n\}$. For $m\in\{\tilde n,\tilde n+v\}$ let ${\bf I}_m(y_1)$ be the $\bar S(y_1)$-scheme that parametrizes isomorphisms between $\Md_{\bar S(y_1)}$ and the pull back of $y_1^*(\Md)$ through the natural morphism $\bar S(y_1)\to {\rm Spec}(K)$. We consider the natural truncation morphism ${\bf T}_{\tilde n,v}:{\bf I}_{\tilde n+v}(y_1)\to {\bf I}_{\tilde n}(y_1)$ of $\bar S(y_1)$-schemes. 
\smallskip
Let ${\bf I}_{\tilde n,v}(y_1)$ be the minimal reduced, closed subscheme of ${\bf I}_{\tilde n}(y_1)_{\rm red}$ through which the reduced morphism defined by ${\bf T}_{\tilde n,v}$ factors. As $\tilde n\ge l$, from (a) we get that ${\bf I}_{\tilde n,v}(y_1)$ is in fact an $S(y_1)$-scheme. The resulting morphism 
$$m_{\tilde n,v}(y_1):{\bf I}_{\tilde n,v}(y_1)\to S(y_1)$$ 
is surjective, cf. the definition of $S(y_1)$. To prove (c) and (d), it suffices to show that $S(y_1)$ is an affine $\bar S(y_1)$-scheme and that ${\bf I}_{\tilde n,v}(y_1)$ is a regular scheme that is finite, flat over $S(y_1)$. It suffices to check this under the extra assumption that $S$ is affine. So the schemes $S_K$, $\bar S(y_1)$, ${\bf I}_{\tilde n}(y_1)$, and ${\bf I}_{\tilde n,v}(y_1)$ are also affine. 
\smallskip
We check that the surjective morphism $m_{\tilde n,v}(y_1):{\bf I}_{\tilde n,v}(y_1)\to S(y_1)$ 
is quasi-finite above any point $y_{\rm gen}$ of $S(y_1)$ of codimension $0$. Let $F_{\rm gen}$ and $I_{\rm gen}$ be the fibres over $y_{\rm gen}$ of ${\bf I}_{\tilde n,v}(y_1)$ and ${\bf I}_{\tilde n+v}(y_1)$ (respectively). We show that the assumption that $F_{\rm gen}$ is not of dimension 0 leads to a contradiction. This assumption implies that the image of $I_{\rm gen}$ in $F_{\rm gen}$ contains an open, dense subscheme of $F_{\rm gen}$ of positive dimension. We get the existence of an algebraically closed field $L$ that contains the residue field of $y_{\rm gen}$ and such that the number of automorphisms of $y_1^*(\Md)_L[p^{\tilde n}]$ that lift to automorphisms of $y_1^*(\Md)_L[p^{\tilde n+v}]$ is infinite. From this and 5.1.1 (c) we get that the image of ${\rm Aut}(y_1^*(\Md)_L)={\rm Aut}(y_1^*(\Md))$ in ${\rm Aut}(y_1^*(\Md)_L[p^{\tilde n}])$ is infinite. But ${\rm Aut}(y_1^*(\Md))$ is a ${\bf Z}_p$-algebra of finite rank and so this image is finite. Contradiction.
\smallskip
So $F_{\rm gen}$ has dimension $0$. Thus there exists an open, dense subscheme $U(y_1)$ of $S(y_1)$ such that the reduced ${\rm Spec}(K)$-scheme of finite type ${\bf I}_{\tilde n,v}(y_1)\times_{\bar S(y_1)} U(y_1)$ is regular as well as (cf. [18, Ch. IV, (9.6.1) and (11.1.1)] and the surjectivity of $m_{\tilde n,v}$) finite, flat over $U(y_1)$. From (ii) and constructions we get the existence of an isomorphism of $S_K$-schemes
$$K_{12}:I_{y_1}\times_{S_K} {\bf I}_{\tilde n,v}(y_1) \tilde\to I_{y_2}\times_{S_K} {\bf I}_{\tilde n,v}(y_1)$$ 
such that we have $(1_{I_{y_2}}\times_{S_K} m_{\tilde n,v}(y_1))\circ K_{12}=J_{12}\times_{S(y_1)} m_{\tilde n,v}(y_1)$.
In particular, we get:
\medskip
{\bf (v)} the morphism $m_{\tilde n,v}(y_1)$ is finite and flat above an open subscheme of $S(y_1)$ that contains $y_1$ if and only if it is so above an open subscheme of $S(y_1)$ that contains $y_2$.
\medskip
As in the above part that pertains to local geometries, the existence of such isomorphisms $K_{12}$ of $S_K$-schemes implies that ${\bf I}_{\tilde n,v}(y_1)$ is regular and equidimensional. From (v) and the existence of $U(y_1)$ we get that ${\bf I}_{\tilde n,v}(y_1)$ is a finite, flat $S(y_1)$-scheme. From this and the fact that ${\bf I}_{\tilde n,v}(y_1)$ is a regular subscheme of ${\bf I}_{\tilde n}(y_1)$, we get that (d) holds for $\tilde n$ (and so also for $q\le \tilde n$) and that ${\bf I}_{\tilde n,v}(y_1)$ is the normalization of $S(y_1)$ in the ring of fractions of ${\bf I}_{\tilde n,v}(y_1)$. As ${\bf I}_{\tilde n,v}(y_1)$ is affine and $m_{\tilde n,v}(y_1)$ is a finite, surjective morphism, from Chevalley theorem of [16, Ch. II, (6.7.1)] we get that $S(y_1)$ is affine. So (c) holds.\endproof
\medskip
{\bf 5.3.2. Ultimate stratifications.} Let $N\ge 3$ and $\Ma_{d,1,N}$ be as in 1.5. Let $(\Ma,\lambda_{\Ma})$ be the universal principally polarized abelian scheme over $\Ma_{d,1,N}$. We have:
\medskip
{\bf (a)} {\it There exists a stratification $\Ms_{d,N}$ of $\Ma_{d,1,N}$ defined by the following property: two geometric points $y_1$ and $y_2$ of $\Ma_{d,1,N}$ with values in the same algebraically closed field $K$, factor through the same stratum of $\Ms_{d,N}$ if and only if the principally quasi-polarized $p$-divisible groups of $y_1^*(\Ma,\lambda_{\Ma})$ and $y_2^*(\Ma,\lambda_{\Ma})$, are isomorphic.} 
\smallskip
{\bf (b)} {\it The stratification $\Ms_{d,N}$ of $\Ma_{d,1,N}$ satisfies the purity property and its strata are regular and equidimensional.}
\smallskip
{\bf (c)} {\it Let $q\in {\bf N}$. Let $K$ be as in (i). Then for any stratum $S_0$ of the stratification $\Ms_{d,N}$ that is a subscheme of ${\Ma_{d,1,N}}_K$, there exists a regular scheme $S_0[q]$ that is finite, flat over $S_0$ and such that the pull back to $S_0[q]$ of the principally quasi-polarized truncated Barsotti--Tate group of level $q$ of $(\Ma,\lambda_{\Ma})$ is constant, i.e. it is the pull back to $S_0[q]$ of a principally quasi-polarized truncated Barsotti--Tate group of level $q$ over ${\rm Spec}(K)$.}
\medskip
The proofs  of (a) to (c) are the same as of 5.3.1 (b) to (d), cf. 3.2.5. We only have to add that the use of [21, 4.8] in the proof of 5.3.1 (b) has to be substituted by the well known fact that the formal deformation spaces of a principally polarized abelian variety over ${\rm Spec}(K)$ and of its principally quasi-polarized $p$-divisible group, are naturally identified (cf. Serre--Tate deformation theory of [23, Ch. 1]). 
\smallskip
Let ${\rm Sch}_{\rm red}^{{\bf F}_p}$ be the category of reduced ${\rm Spec}({\bf F}_p)$-schemes endowed with the \'etale topology. Let $\Ma_{d,1}$ be the moduli stack over ${\rm Sch}_{\rm red}^{{\bf F}_p}$ of principally polarized abelian schemes of relative dimension $d$ (see [13, Ch. I, \S4, p. 17 and 4.3]). The stratification $\Ms_{d,N}$ descends to a stratification $\Ms_d$ of $\Ma_{d,1}$. As we did not formalize stratifications of stacks, we describe $\Ms_d$ directly as follows. 
\smallskip
We fix a principally quasi-polarized $p$-divisible group $\Mt:=(D,\lambda_D)$ over ${\rm Spec}(k)$ of height $r=2d$. The objects of $\Ma_{d,1}$ are principally polarized abelian schemes over reduced ${\rm Spec}({\bf F}_p)$-schemes. The substack $\Ma_{d,1}(\Mt)$ of $\Ma_{d,1}$ associated to $\Mt$ is the full subcategory of $\Ma_{d,1}$ whose objects are principally polarized abelian schemes over reduced ${\rm Spec}({\bf F}_p)$-schemes with the property that all principally quasi-polarized $p$-divisible groups obtained from them via pulls back through points with values in the same algebraically closed field $K$ that contains $k$, are isomorphic to $\Mt\times_{{\rm Spec}(k)} {\rm Spec}(K)$. 
\smallskip
We refer to $\Ms_d$ (resp. to $\Ms_{d,N}$) as the {\it ultimate stratification} of $\Ma_{d,1}$ (resp. of  $\Ma_{d,1,N}$). 
\smallskip
Let $\Ma_{d,1}(\Mt)_k:=\Ma_{d,1}(\Mt)\times_{{\rm Sch}_{\rm red}^{{\bf F}_p}} {\rm Sch}_{\rm red}^k$ and ${\Ma_{d,1}}_k:=\Ma_{d,1}\times_{{\rm Sch}_{\rm red}^{{\bf F}_p}} {\rm Sch}_{\rm red}^k$, where ${\rm Sch}_{\rm red}^k$ is the full subcategory of ${\rm Sch}_{\rm red}^{{\bf F}_p}$ formed by reduced ${\rm Spec}(k)$-schemes. The pull back of $\Ma_{d,1}(\Mt)_k$ via the 1-morphism ${\Ma_{d,1,N}}_k\to {\Ma_{d,1}}_k$ is the stratum of $\Ms_{d,N}$ that is a subscheme of ${\Ma_{d,1,N}}_k$ and that corresponds naturally to $\Mt$. Using this it can be easily checked that $\Ma_{d,1}(\Mt)_k$ is a separated, algebraic stack over ${\rm Sch}_{\rm red}^k$ in the similar sense as of [13, Ch. I, \S4, 4.6 and 4.8] but worked out using only reduced ${\rm Spec}(k)$-schemes. 
\smallskip
The following Proposition, to which we refer as the {\it integral Manin problem} for Siegel modular varieties (see [28, p. 76] and [39, p. 98] for the original Manin problem), implies that $\Ma_{d,1}(\Mt)$ is a non-empty category. 
\medskip
{\bf 5.3.3. Proposition.} {\it Let ${\rm Spec}(C)$ be a local, complete scheme whose residue field is $k=\bar k$. Then any principally quasi-polarized $p$-divisible group $\Mt_C^\prime$ over ${\rm Spec}(C)$ that lifts $\Mt$, is the one of a principally polarized abelian scheme over ${\rm Spec}(C)$.}
\medskip
{\it Proof:}
We first show that $\Mt$ is associated to a principally polarized abelian variety $(A,\lambda_{A})$ over ${\rm Spec}(k)$. Let $(M,\vph,\lambda_M)$ be the principally quasi-polarized Dieudonn\'e module of $\Mt$. Let $\tilde A$ be an abelian variety $\tilde A$ over ${\rm Spec}(k)$ whose $F$-isocrystal is $(M[{1\over p}],\vph)$, cf. [39, p. 98]. Based on [32, \S23, Cor. 1], up to an isogeny, we can choose $\tilde A$ such that it has a principal polarization $\lambda_{\tilde A}$. The principally quasi-polarized Dieudonn\'e module of $(\tilde A,\lambda_{\tilde A})$ is of the form $(\tilde M,\vph,\lambda_{\tilde M})$, where $\tilde M$ is a certain $W(k)$-lattice of $M[{1\over p}]$. 
\medskip
{\bf 5.3.3.1. Lemma.} {\it If $\lambda_1$ and $\lambda_2$ are two principal quasi-polarizations of $(M[{1\over p}],\vph)$, then the triples $(M[{1\over p}],\vph,\lambda_1)$ and $(M[{1\over p}],\vph,\lambda_2)$ are isomorphic.}
\medskip
{\it Proof:} It suffices to prove the Lemma under the assumption that there exists $\al\in [0,{1\over 2}]\cap{\bf Q}$ such that all slopes of $(M[{1\over p}],\vph)$ are $\al$ and $1-\al$, cf. Dieudonn\'e's classification of $F$-isocrystals over $k$ (see [28, \S2]) and [39, p. 98]. Let $i\in\{1,2\}$. 
\smallskip
We first consider the case when $\al\neq {1\over 2}$; so $\al\neq 1-\al$. Let $M[{1\over p}]=M_{\al}\oplus M_{1-\al}$ be the direct sum decomposition that is normalized by $\vph$ and such that for $\be\in\{\al,1-\al\}$ all slopes of $(M_{\be},\vph)$ are $\be$. We have $\lambda_i(M_{\be},M_{\be})=0$ and the bilinear form $\lambda_{i,\be}:M_{\be}\otimes_{B(k)} M_{1-\be}\to B(k)$ induced naturally by $\lambda_i$, is non-degenerate. But $\lambda_{i,1-\al}$ is determined by $\lambda_{i,\al}$. Thus $\lambda_i$ is uniquely determined by the isomorphism $j_i:(M_{\al},\vph)\tilde\to (M_{1-\al}^*,p1_{M_{1-\al}^*}\vph)$ defined naturally by $\lambda_{i,\al}$ via the rule $j_i(x)(y)=\lambda_i(x,y)=\lambda_{i,\al}(x,y)$, where $x\in M_{\al}$ and $y\in M_{1-\al}$. Let $f_{12}:=j_1^{-1}j_2:(M_{\al},\vph)\tilde\to (M_{\al},\vph)$. The automorphism $e_{12}:=f_{12}\oplus 1_{M_{1-\al}}$ of $(M[{1\over p}],\vph)=(M_{\al},\vph)\oplus (M_{1-\al},\vph)$ 
takes $\lambda_2$ into $\lambda_1$, i.e. for $x,y\in M[{1\over p}]$ we have an identity $\lambda_2(x,y)=\lambda_1(e_{12}(x),e_{12}(y))$.
\smallskip
Let now $\al$ be ${1\over 2}$. As $\al={1\over 2}$, the $F$-isocrystal $(M[{1\over p}],\vph)$ over $k$ is a direct sum of simple $F$-isocrystals  over $k$ of rank $2$. Using the standard argument that shows that any two non-degenerate, symmetric, bilinear forms on an even dimensional complex vector space are isomorphic, we get that both $(M[{1\over p}],\vph,\lambda_1)$ and $(M[{1\over p}],\vph,\lambda_2)$ are direct sums of principally quasi-polarized $F$-isocrystals  over $k$ of rank $2$. Thus we can assume $r_M=2$ (i.e. $d=1$). But in this case the Lemma is trivial (for instance, cf. [27, pp. 35--36]). \endproof
\medskip
Based on the Lemma, it suffices to prove the Proposition under the extra hypothesis that $\lambda_M=\lambda_{\tilde M}$. From the classical Dieudonn\'e theory we get directly the next property.
\medskip
{\bf 5.3.3.2. The isogeny property.} {\it There exists a unique principally polarized abelian variety $(A,\lambda_{A})$ over ${\rm Spec}(k)$ that is ${\bf Z}[{1\over p}]$-isogenous to $(\tilde A,\lambda_{\tilde A})$ and whose principally quasi-polarized Dieudonn\'e module is identifiable under this ${\bf Z}[{1\over p}]$-isogeny with $(M,\vph,\lambda_M)$.}
\medskip
So $\Mt$ is the principally quasi-polarized $p$-divisible group of $(A,\lambda_{A})$. From Serre--Tate deformation theory (see [23, Ch. 1]) and Grothendieck algebraization theorem (see [17, Ch. III, Thm. (5.4.5)]), we easily get the existence of a principally polarized abelian scheme over ${\rm Spec}(C)$ whose principally quasi-polarized $p$-divisible group is $\Mt_C^\prime$.\endproof
\medskip
{\bf 5.3.4. Remarks.} {\bf (a)} If $d\ge 2$, then the stratification $\Ms_{d,N}$ has a class which is not a set (for $d\ge 3$ this follows from 4.5.4).
\smallskip
{\bf (b)} Let $q\in{\bf N}$. Let $\Ms_{d,N,q}$ be the stratification of $\Ma_{d,1,N}$ defined by the rule: two geometric points $y_1$ and $y_2$ of $\Ma_{d,1,N}$ with values in the same algebraically closed field, factor through the same stratum of $\Ms_{d,N,q}$ if and only if the principally quasi-polarized truncated Barsotti--Tate groups of level $q$ of $y_1^*(\Ma,\lambda_{\Ma})$ and $y_2^*(\Ma,\lambda_{\Ma})$, are isomorphic. The strata of $\Ms_{d,N,q}$ are regular and equidimensional (one argues this in a way similar to the proof of 5.3.1; instead of [21, 4.8] one has to use a principal quasi-polarized version of [21, 4.7 and 4.8] and Serre--Tate deformation theory). For $q\ge T(d)$ we have $\Ms_{d,N,q}=\Ms_{d,N}$, cf. 3.2.5. The case $q=1$ was first studied by Ekedahl and Oort, cf. [34, \S1]. The strata of $\Ms_{d,N,1}$ are quasi-affine, cf. [34, 1.2]. As each stratum of $\Ms_{d,N,q}$ is a locally closed subscheme of a stratum of $\Ms_{d,N,1}$, the strata of $\Ms_{d,N,q}$ are also quasi-affine. For $1\le q< T(d)$, we do not know when the stratification $\Ms_{d,N,q}$ satisfies the purity property. 
\smallskip
{\bf (c)} The existence of the ultimate stratifications $\Ms_d$ and $\Ms_{d,N}$, though of foundation, is only a first step toward the solution of the below Main Problem. Due to the importance of Main Problem, we will state a general form of it, even if in this paper we do not formalize specializations of latticed $F$-isocrystals with a group (such specializations are standard for $p$-divisible groups; see also 3.2.7). To be short, we state Main Problem only in a context that involves tensors but no principal bilinear quasi-polarizations.  
\medskip
{\bf Main Problem.} {\it Let $(M,\vph,G,(t_{\al})_{\al\in\Mj})$ be a latticed $F$-isocrystal with a group and an emphasized family of tensors over $k$ such that the $W$-condition holds for $(M,\vph,G)$ (see 2.2.1 (b) and (d)). List using families all isomorphism classes of $(M,g\vph,G,(t_{\al})_{\al\in\Mj})$'s (where $g\in G(W(k))$) and decide which such classes specialize to which other.}   
\medskip\smallskip
{\bf 5.4. On the specialization theorem.} Let $S$ be an integral ${\rm Spec}({\bf F}_p)$-scheme. We take $k$ to be an algebraic closure of the field of fractions $k_S$ of $S$. Let ${\got C}$ be an $F$-crystal over $S$. Let $h_{\got C}\in{\bf N}\cup\{0\}$ be as in 2.1; the $h$-number of any pull back of ${\got C}$ via a geometric point of $S$ is at most $h_{\got C}$. Let $\Mn$ be the Newton polygon of ${\got C}_k$. Let 
$$U^{\rm top}:=\{x\in S^{\rm top}|x^*({\got C})\;\; {\rm has}\;\;{\rm Newton}\;\;{\rm polygon}\;\;\Mn\}.$$ 
We recall that Grothendieck proved that for any geometric point $y$ of $S$ the Newton polygon of $y^*({\got C})$ is above $\Mn$ (see [19, Appendice]) and that Katz added that moreover there exists an open subscheme $U$ of $S$ such that the notations match, i.e. $U^{\rm top}$ is the topological space underlying $U$ (see [22, 2.3.1 and 2.3.2]). 
\smallskip
We give another proof of the existence of $U$ using Grothendieck's result. This result implies that if $x\in U^{\rm top}$, then all points of the spectrum of the local ring of $x$ in $S$ belong to $U^{\rm top}$. To show the existence of $U$, it is enough to show that there exists a non-empty open subscheme $U^{\prime}$ of $S$ such that $U^{\prime {\rm top}}\subseteq U^{\rm top}$. The argument for this goes as follows. The existence of such open subschemes $U^\prime$ implies that $U^{\rm top}$ is an ind-constructible subset of $S^{\rm top}$, cf. [18, Ch. IV, (9.2.1) and (9.2.3)]. Based on this and the above part that pertains to $x\in U^{\rm top}$, from [18, Ch. IV, Thm. (1.10.1)] we get that each point of $U^{\rm top}$ is an interior point of $U^{\rm top}$. Thus $U^{\rm top}$ is an open subset of $S^{\rm top}$ and therefore $U$ exists. 
\smallskip
Let ${\got C}_0$ be a Dieudonn\'e--Fontaine $p$-divisible object over ${\bf F}_p$ of Newton polygon $\Mn$. Let $h_0$ be the $h$-number of ${\got C}_0$ and let $r_0$ be the rank of ${\got C}_0$. Let $i:{\got C}_{0k}\to {\got C}_k$ be an isogeny. Let $l\in{\bf N}$ be such that $p^l$ annihilates ${\rm Coker}(i)$. Let $v:={\rm max}\{v(r_0,r_0,h_0,b)|b\in S(0,h_{\got C})\}$ and $n:={\rm max}\{1,n(r_0,r_0,h_0,b)|b\in S(0,h_{\got C})\}$ be defined using the numbers of 5.1.1 (b). 
\smallskip
Let $i(n+v+l)$ be the reduction of $i$ mod $p^{n+v+l}$; it is a morphism of $\Mm(W_{n+v+l}(k))$ whose cokernel is annihilated by $p^l$. From 2.8.3 (see 2.8.3 (a) applied with $(V_1,V,q)=(k,k_S,n+v+l)$), we get that there exists a finite field extension $k_{\tilde S}$ of $k_S$ such that $i(n+v+l)$ is the pull back of a morphism $i(n+v+l)_{k_{\tilde S}}$ of $\Mm(W_{n+v+l}(k_{\tilde S}))$ whose cokernel is annihilated by $p^l$. Let $\tilde S$ be the normalization of $S$ in $k_{\tilde S}$ (the notations match, i.e. $k_{\tilde S}$ is the field of fractions of $\tilde S$). The continuous map $\tilde S^{\rm top}\to S^{\rm top}$ is proper, cf. the going-up theorem of [29, (5.E)]. So if there exists an open, dense subscheme $\tilde U^\prime$ of $\tilde S$ with the property that $\tilde U^{\prime {\rm top}}$ maps into $U^{\rm top}$, then we can take $U^\prime$ to be the complement in $S$ of the image of $\tilde S^{\rm top}\setminus\tilde U^{\prime{\rm top}}$ in $S^{\rm top}$. Thus it suffices to consider the case when $k_{\tilde S}=k_S$. Let $U^{\prime}$ be an open subscheme of $S$ such that we have a morphism 
$$i_{U^{\prime}}(n+v+l):{\bf E}({\got C}_0;W_{n+v+l}(U^{\prime}))\to {\bf E}({\got C};W_{n+v+l}(U^{\prime}))$$
of $\Mm(W_{n+v+l}(U^{\prime}))$ that extends $i(n+v+l)_{k_S}$ and that has a cokernel annihilated by $p^l$, cf. 2.8.3 (b). For an arbitrary geometric point $z:{\rm Spec}(K)\to U^{\prime}$, the reduction  mod $p^{n+l}$ of the pull back morphism $z^*(i_{U^{\prime}}(n+v+l))$ lifts to a morphism $i_{z}:{\got C}_{0K}\to z^*({\got C})$
of $F$-crystals over $K$ (cf. 5.1.1 (b), 5.1.2,  and the definitions of $v$ and $n$). As the cokernel of the reduction mod $p^{n+l}$ of $i_{z}$ is annihilated by $p^l$ and as $n\ge 1$, (by reasons of ranks) the morphism $i_{z}$ is injective and thus an isogeny. Thus $z^*({\got C})$ has Newton polygon $\Mn$. So $U^{\prime{\rm top}}\subseteq U^{\rm top}$. This ends the argument for the existence of $U^{\prime}$ and so also of $U$.
\bigskip\smallskip
\centerline{\bigsll {\bf \S6 Proof of Main Theorem B}}
\bigskip\smallskip
Let $S$ be a reduced ${\rm Spec}({\bf F}_p)$-scheme. Let ${\got C}$ be an $F$-crystal over $S$. Let $\Ms({\got C})$ be the Newton polygon stratification of $S$ defined by ${\got C}$, cf. [22, 2.3.1 and 2.3.2]. The stratification $\Ms({\got C})$ is of finite type and locally in the Zariski topology of $S$ has a finite number of strata. The main goal of this Section is to prove Main Theorem B stated in 1.6, i.e. to prove that $\Ms({\got C})$ satisfies the purity property (see 6.2). In 6.1 we capture the very essence of Main Theorem B for the case when $S$ is an integral, locally noetherian scheme. In 6.3 we include two remarks on the connection between 6.1 and a result of de Jong and Oort and on Newton polygon stratifications defined by certain reductions modulo powers of $p$ of $F$-crystals. We will use the notations of 2.8.2. 
\medskip\smallskip
{\bf 6.1. Theorem.} {\it Suppose $S$ is integral and locally noetherian. Let $U$ be the maximal open subscheme of $S$ with the property that the Newton polygons of pulls back of ${\got C}$ through geometric points of $U$ are all equal (see 5.4). Then $U$ is an affine $S$-scheme.}
\medskip
{\it Proof:} It suffices to prove this under the extra assumptions that $S={\rm Spec}(R)$ is affine and that the underlying $R$-module of ${\bf E}({\got C};W_1(S))$ is free. Let $R_U$ be the $R$-algebra of global functions of $U$. We have to show that $U$ is affine, i.e. the natural and functorial morphism $f_U:U\to {\rm Spec}(R_U)$ is an isomorphism. This statement is local in the faithfully flat topology of $S$ and thus we can assume that $S$ is local. Let $\hat R$ be the completion of $R$ and let $\hat S:={\rm Spec}(\hat R)$. As $\hat S$ is a faithfully flat $S$-scheme, to show that $U$ is affine (i.e. $f_U$ is an isomorphism) it suffices to show that $U\times_S \hat S$ is affine (i.e. $f_U\times_S \hat S=f_{U\times_S \hat S}$ is an isomorphism). Let $\hat S_1={\rm Spec}(\hat R_1),\ldots,\hat S_j={\rm Spec}(\hat R_j)$ be the irreducible components of the reduced scheme of $\hat S$ (here $j\in{\bf N}$); they are spectra of local, complete, integral, noetherian  ${\bf F}_p$-algebras. The scheme $U\times_S \hat S$ is affine if and only if the irreducible components $U\times_S \hat S_1,\ldots,U\times_S \hat S_j$ of the reduced scheme of $U\times_S \hat S$ are all affine (cf. [16, Ch. II, Cor. (6.7.3)]). So to prove the Theorem we can assume $R=\hat R=\hat R_1$. As $R$ is a local, complete ring, it is also excellent (cf. [29, (34.B)]). Thus the normalization $S^{\rm n}$ of $S$ is a finite $S$-scheme.  So $S^{\rm n}$ is a semilocal, complete, integral, normal  scheme. This implies that $S^{\rm n}$ is local. But $U$ is affine if and only if $U\times_S S^{\rm n}$ is affine, cf. 2.9.2. Thus to prove the Theorem, we can also assume $S$ is normal; so $S=S^{\rm n}$. 
\smallskip
We emphasize that for the rest of the proof we will only use the fact that $S$ is an integral,  normal, excellent, affine scheme (but not necessarily local and thus not necessarily complete). We group the main steps into distinct (and numbered) Subsubsections. 
\medskip
{\bf 6.1.1. Notations and two operations.} Let $k_S$, $k$, $h_{\got C}$, $\Mn$, ${\got C}_0$, $h_0$, $r_0$, $v$, and $n$ be as in 5.4. So $\Mn$ is the Newton polygon of pulls back of ${\got C}_U$ via geometric points of $U$, ${\got C}_0$ is a Dieudonn\'e--Fontaine $p$-divisible object over ${\bf F}_p$ that has Newton polygon $\Mn$, $r_0$ is the rank of ${\got C}_0$, etc. Let $q_0:=r_0!$. Below all pulls back to ${\rm Spec}(k)$ of $F$-crystals are via the natural dominant morphisms ${\rm Spec}(k)\to {\rm Spec}({\bf F}_p)$ and ${\rm Spec}(k)\to S$. 
\smallskip
We consider the following two replacement operations ($\Mr1$) and ($\Mr2$) of the triple $(S,U,{\got C})$ by a new triple $(\tilde S,\tilde U,\tilde{\got C})$. For both operations $\tilde S$ is an integral, normal, affine $S$-scheme of finite type, $\tilde{\got C}$ is ${\got C}_{\tilde S}$, and
\medskip
($\Mr1$) {\it either $(\tilde S,\tilde U)$ is the normalization of $(S,U)$ in a finite field extension of $k_S$} 
\smallskip
($\Mr2$) {\it or $U$ is an open subscheme of $\tilde S$ and $\tilde U:=U$ is $U\times_S \tilde S$.}
\medskip
The scheme $\tilde S$ is also excellent, cf. [29, (34.B)]. Moreover, $U$ is affine if and only if $\tilde U$ is affine (in connection with ($\Mr1$), cf. 2.9.2). So in what follows we will often perform one of these two operations in order to simplify the setting and to eventually end up with a situation where in fact we have $U=S$. By performing $(\Mr1)$, we can assume that $R$ is an ${\bf F}_{p^{q_0}}$-algebra (i.e. ${\bf F}_{p^{q_0}}\hookrightarrow R$).
\smallskip
Let $V$ be a local ring of $U$ that is a discrete valuation ring. Let $V_2$ be a complete discrete valuation ring  that is a faithfully flat $V$-algebra, and that has an algebraically closed residue field $k_2$. Let $V_1:=V_2^{\rm perf}$. We fix an isomorphism $V_2\tilde\to k_2[[w]]$ and we view it as an identification under which $V_2$ and $V_1$ become $k_2$-algebras. Let $\Phi_2$ be the Frobenius endomorphism of $W(k_2)[[w]]$ that takes $w$ into $w^p$ and is compatible with $\sg_{k_2}$. 
\medskip
{\bf 6.1.2. Key Lemma.} {\it There exists a number $l\in{\bf N}$ that is greater than $\max\{h_0r_0^2,h_{\got C}\}$, that depends only on ${\got C}_0$ and $h_{\got C}$ but not on $V$, and such that there exists an isogeny $i_{V_1}:{\got C}_{0V_1}\to {\got C}_{V_1}$ of $F$-crystals over ${\rm Spec}(V_1)$ whose cokernel is annihilated by $p^l$.}
\medskip
{\it Proof:} We first show by induction on $r_0\in {\bf N}$ that there exists a number $\tilde l_{C}\in{\bf N}$ that does not depend on $V$ but only on ${\got C}_0$ and $h_{\got C}$ and such that we have an isogeny $i_1:{\got C}_{CV_1}\to {\got C}_{V_1}$ whose cokernel is annihilated by $p^{\tilde l_{C}}$, where ${\got C}_C$ is an $F$-crystal over ${\rm Spec}(k_2)$ and where the role of ${\got C}_{V_1}$ is that of the pull back to ${\rm Spec}(V_1)$ of an arbitrary $F$-crystal over ${\rm Spec}(V_2)$ of constant Newton polygon which depends only on $\Mn$. 
\smallskip
Let $\al_1$ be the smallest slope of $\Mn$. Let $\tilde l_1:=-[-\al_1(r_0-1)]$. Theorem [22, 2.6.1] says that there exists an isogeny $i_1^\prime:{\got C}^\prime\to {\got C}_{V_2}$, where ${\got C}^\prime$ is an $\al_1$-divisible $F$-crystal over ${\rm Spec}(V_2)$. The $\al_1$-divisibility means that if $(M^\prime,\vph^\prime,\nabla^\prime)$ is the evaluation of ${\got C}^\prime$ at the thickening defined by the closed embedding ${\rm Spec}(V_2)\hookrightarrow {\rm Spec}(W(k_2)[[w]])$, then for all $u\in{\bf N}$ the $\Phi_2^u$-linear endomorphism $(\vph^{\prime})^u$ of $M^\prime$ is divisible by $p^{[u\al_1]}$. We can choose ${\got C}^\prime$ and $i_1^\prime$ such that ${\rm Coker}(i_1^\prime)$is annihilated by $p^{\tilde l_1}$, cf. [22, p. 153]. If $\al_1$ is the only slope of $\Mn$, then ${\got C}^\prime$ is the pull back of an $F$-crystal over ${\rm Spec}(k_2)$ (cf. [22, proof of 2.7.1]); so we can take $\tilde l_{C}$ to be $\tilde l_1$ and $i_1$ to be $i_1^\prime$. In particular, $\tilde l_{C}$ exists if $r_0=1$. 
\smallskip
We now consider the case when $\Mn$ has at least two slopes. From [22, proof of 2.6.2] we get that we have a unique short exact sequence 
$$0\to {\got C}^\prime_1\to {\got C}^\prime\to {\got C}^\prime_2\to 0$$ 
of $F$-crystals over ${\rm Spec}(V_2)$ with the property that the Newton polygons of the pulls back of ${\got C}^\prime_1$ (resp. of ${\got C}^\prime_2$) via geometric points of ${\rm Spec}(V_2)$ have all slopes equal to $\al_1$ (resp. have all slopes greater than $\al_1$). As $V_1$ is perfect, loc. cit. also proves that this short exact sequence has a unique splitting after we pull it back to ${\rm Spec}(V_1)$. Thus we have a unique direct sum decomposition ${\got C}^\prime_{V_1}={\got C}^\prime_{1V_1}\oplus {\got C}^\prime_{2V_1}$ of $F$-crystals over ${\rm Spec}(V_1)$. Using this decomposition and the fact that both $F$-crystals ${\got C}^\prime_1$ and ${\got C}^\prime_2$ over ${\rm Spec}(V_2)$ have ranks smaller than $r_0$, by induction we get the existence of $\tilde l_{C}$. This ends the induction.
\smallskip
Let $\tilde l_{f}:=\max\{d(r_0,0,c)|c\in S(0,h_{\got C}+\tilde l_{C})\}$ with $d(r_0,0,c)$'s as in 2.4.1. The $h$-number of ${\got C}_C$ is at most $h_{\got C}+\tilde l_{C}$. So from 2.4.1 (applied over $k_2$) we get that there exists a Dieudonn\'e--Fontaine $p$-divisible object ${\got C}_0^\prime$ over ${\bf F}_p$ for which we have an isogeny $i_2:{\got C}_{0V_1}^\prime\to {\got C}_{CV_1}$ whose cokernel is annihilated by $p^{\tilde l_{f}}$. As the number of isomorphism classes of Dieudonn\'e--Fontaine $p$-divisible objects over ${\bf F}_p$ that have $\Mn$ as their Newton polygons is finite, there exists a number $\tilde l_{\Mn}\in{\bf N}$ such that we have an isogeny $\break i_3:{\got C}_{0V_1}\to {\got C}_{0V_1}^\prime$ whose cokernel is annihilated by $p^{\tilde l_{\Mn}}$. 
\smallskip
As $i_{V_1}$ we can take the composite isogeny $i_1\circ i_2\circ i_3$. Thus as $l$ we can take any integer greater than $\max\{h_0r_0^2,h_{\got C},\tilde l_{C}+\tilde l_{f}+\tilde l_{\Mn}\}$.\endproof 
\medskip
{\bf 6.1.3. The open subscheme $U_0$.} With $l$ as in 6.1.2, let $m:=8l+n+v+1$. We continue the proof of Theorem 6.1 by considering (see 2.8.2) the evaluation morphism 
$${\bf E}(i_{V_1};W_{m+2v}(V_1)):{\bf E}({\got C}_0;W_{m+2v}(V_1))\to {\bf E}({\got C};W_{m+2v}(V_1))$$ 
of $\Mm(W_{m+2v}(V_1))$. We apply 2.8.3 (c) (with $q=m+2v$) to this morphism. We get that there exists a finite field extension $k_{S,V}$ of $k_S$ and an open, affine subscheme $U_{\tilde V}$ of the normalization of $U$ in $k_{S,V}$, such that $U_{\tilde V}$ has a local ring $\tilde V$ which is a discrete valuation ring and which dominates $V$ and moreover we have a morphism 
$$i_{U_{\tilde V}}(m+2v):{\bf E}({\got C}_0;W_{m+2v}(U_{\tilde V}))\to {\bf E}({\got C};W_{m+2v}(U_{\tilde V}))$$ 
of $\Mm(W_{m+2v}(U_{\tilde V}))$ whose cokernel is annihilated by $p^l$. See 2.1 for $i_{U_{\tilde V}}(m+v)$.
\smallskip
Let $\tilde m\in\{m,m+v\}$. Let $I_{\tilde m}$ be the set of morphisms
${\got C}_{0k}/p^{\tilde m}{\got C}_{0k}\to {\got C}_{k}/p^{\tilde m}{\got C}_{k}$ that are reductions mod $p^{\tilde m}$ of morphisms ${\got C}_{0k}/p^{\tilde m+v}{\got C}_{0k}\to {\got C}_{k}/p^{\tilde m+v}{\got C}_{k}$. Any morphism in $I_{\tilde m}$ lifts to a morphism ${\got C}_{0k}\to {\got C}_{k}$, cf. 5.1.1 (b), 5.1.2, and the definitions of $v$ and $n<m$; thus $I_{\tilde m}$ is a finite set. Let $J_{\tilde m}:=\{i\in I_{\tilde m}|p^l{\rm Coker}(i)=0\}$. Based on the case 2.8.3 (a) of 2.8.3, by performing $(\Mr1)$ we can assume that $J_{\tilde m}$ is the set of pulls back of a set of morphisms $L_{\tilde m}$ of $\Mm(W_{\tilde m}(k_S))$ whose cokernels are annihilated by $p^l$. The pull back of $i_{U_{\tilde V}}(m+v)$ to a morphism of $\Mm(W_{m+v}(k_{S,V}))$ is also the pull back of a morphism in $L_{m+v}$. As $V=\tilde V\cap k_S$, inside $W_{m+v}(k_{S,V})$ we have $W_{m+v}(V)=W_{m+v}(\tilde V)\cap W_{m+v}(k_S)$. This implies that the pull back of $i_{U_{\tilde V}}(m+v)$ to a morphism of $\Mm(W_{m+v}(\tilde V))$ is in fact the pull back of a morphism of $\Mm(W_{m+v}(V))$ whose cokernel is generically annihilated by $p^l$ (in the sense of 2.8.1). From the case 2.8.3 (b) of 2.8.3 (applied with $(V_1,V)$ replaced by $(V,R)$), we get the existence of an open subscheme $U_V$ of $U$ that has $V$ as a local ring and such that we have a morphism 
$$i_{U_V}(m+v):{\bf E}({\got C}_0;W_{m+v}(U_V))\to {\bf E}({\got C};W_{m+v}(U_V))$$ 
of $\Mm(W_{m+v}(U_V))$ whose cokernel is generically annihilated by $p^l$. 
\smallskip
For $i\in J_m$ let $\Mv(i)$ be the set of all those discrete valuation rings $V$ of $U$ such that the pull back of $i_{U_V}(m)$ to a morphism of $\Mm(W_m(k))$ is $i$. Let $U_i:=\cup_{V\in \Mv(i)} U_V$. 
Let
$$i_{U_i}(m):{\bf E}({\got C}_0;W_m(U_i))\to {\bf E}({\got C};W_m(U_i))$$
be the morphism of $\Mm(W_m(U_i))$ which is obtained by gluing together the morphisms $i_{U_V}(m)$'s with $V\in \Mv(i)$. We have:
\medskip
{\bf (a)} if $U_0:=\cup_{i\in J_m} U_i$, then $U^{\rm top}\setminus U_0^{\rm top}$ has codimension at least $2$ in $U^{\rm top}$;
\smallskip
{\bf (b)} for any $i\in J_m$, the cokernel of $i_{U_i}(m)$ is generically annihilated by $p^l$.
\medskip
{\bf 6.1.4. Gluing morphisms.} We now modify the morphisms $i_{U_i}(m)$'s ($i\in J_m$) so that they glue together to define a morphism 
$$i_{U_0}(m):{\bf E}({\got C}_0;W_m(U_0))\to {\bf E}({\got C};W_m(U_0))$$
of $\Mm(W_m(U_0))$ whose cokernel is generically annihilated by $p^{3l}$. If there exists $i\in J_m$ such that $U_i=U_0$, then $i_{U_0}(m)=i_{U_i}(m)$ has a cokernel generically annihilated by $p^{l}$. 
\smallskip
We now assume that for all $i\in J_m$ we have $U_0\neq U_i$. The pull back of $i_{U_i}(m)$ to a morphism of $\Mm(W_m(k))$ is $i\in J_m$. Let $f_i:{\got C}_{0k}\to {\got C}_{k}$ be a morphism such that its reduction mod $p^{m}$ is $i$, cf. 5.1.1 (b), 5.1.2, and the fact that $J_m\subseteq I_m$. As $m>l$ and $i\in J_m$, the cokernel of $f_i$ is annihilated by $p^l$ and so $f_i$ is an isogeny. 
\smallskip
Let $f_0\in\{p^lf_i|i\in J_m\}$. We have:
\medskip
{\bf (a)} the image of $f_0$ lies inside the intersection of the images of all $f_i$'s ($i\in J_m$);
\smallskip
{\bf (b)} the cokernel of $f_0$ is annihilated by $p^{2l}$.
\medskip
Let $s_i:{\got C}_{0k}\to{\got C}_{0k}$ be the isogeny such that we have $f_0=f_i\circ s_i$, cf. (a). So ${\rm Coker}(s_i)$ is annihilated by $p^{2l}$, cf. (b). We know that $s_i$ is the pull back of an isogeny ${\got C}_{0{\bf F}_{p^{q_0}}}\to{\got C}_{0{\bf F}_{p^{q_0}}}$, cf. 2.2.3 (b) applied to ${\got C}_0$. So as ${\bf F}_{p^{q_0}}\hookrightarrow R$, we get that the reduction mod $p^{m}$ of $s_i$ is the pull back of a (constant) morphism 
$$s_{U_i}(m):{\bf E}({\got C}_0;W_m(U_i))\to {\bf E}({\got C}_0;W_m(U_i))$$ 
of $\Mm(W_m(U_i))$ whose cokernel is annihilated by $p^{2l}$. If $i_1$, $i_2\in J_m$, then the pulls back of $i_{U_{i_1}}(m)\circ s_{U_{i_1}}(m)$ and $i_{U_{i_2}}(m)\circ s_{U_{i_2}}(m)$ to morphisms of $\Mm(W_m(k))$ are the reduction of $f_0$ mod $p^m$ and thus they coincide. This implies that the pulls back of $i_{U_{i_1}}(m)\circ s_{U_{i_1}}(m)$ and $i_{U_{i_2}}(m)\circ s_{U_{i_2}}(m)$ to morphisms of $\Mm(W_m(U_{i_1}\cap U_{i_2}))$ coincide. Thus the morphisms $i_{U_{i}}(m)\circ s_{U_{i}}(m)$ indexed by $i\in J_m$ glue together to define a morphism $i_{U_0}(m):{\bf E}({\got C}_0;W_m(U_0))\to {\bf E}({\got C};W_m(U_0))$ of $\Mm(W_m(U_0))$ whose cokernel is generically annihilated by $p^{3l}=p^lp^{2l}$, cf. 6.1.3 (b) and the fact that $p^{2l}$ annihilates ${\rm Coker}(s_{U_i}(m))$. 
\smallskip
It is easy to see that by performing $(\Mr1)$ we can assume $p^{3l}$ annihilates ${\rm Coker}(i_{U_0}(m))$ but this will not be used in what follows.
\medskip
{\bf 6.1.5. Lemma.} {\it By performing $(\Mr2)$, we can assume that $i_{U_0}(m)$ extends to a morphism $i_S(m):{\bf E}({\got C}_0;W_m(S))\to {\bf E}({\got C};W_m(S))$ of $\Mm(W_m(S))$.}
\medskip
{\it Proof:} Let $S^\prime$ be the affine $S$-scheme of finite type that parametrizes morphisms between the two objects ${\bf E}({\got C}_0;W_m(S))$ and ${\bf E}({\got C};W_m(S))$ of $\Mm(W_m(S))$, cf. 2.8.4.1. Let $U_0\hookrightarrow S^\prime$ be the open embedding of $S$-schemes that defines $i_{U_0}(m)$. Let $U^\prime$ be the normalization of the Zariski closure of $U_0$ in $S^\prime$. As $S$ is an excellent scheme, the $S$-scheme $U^\prime$ is integral, normal, affine, and of finite type. As $U_0$ is an open subscheme of both $U$ and $U^\prime$ and due to 6.1.3 (a), the affine morphism $U^\prime\times_S U\to U$ between integral, normal, noetherian schemes is birational and has the property that any discrete valuation ring of $U$ is also a local ring of $U^\prime\times_S U$. Thus the morphism $U^\prime\times_S U\to U$ is an isomorphism, cf. 2.9.1. So $U$ is an open subscheme of $U^\prime$. So by performing $(\Mr2)$ (with $\tilde S=U^\prime$), we can assume $U^\prime=S$. Thus we can speak about the morphism $i_S(m):{\bf E}({\got C}_0;W_m(S))\to {\bf E}({\got C};W_m(S))$ of $\Mm(W_m(S))$ that extends $i_{U_0}(m)$.\endproof
\medskip
{\bf 6.1.6. Duals.} Let ${\got C}(l)$ be the Tate twist of ${\got C}$ by $(l)$, i.e. ${\got C}$ tensored with the pull back to $S$ of the $F$-crystal $({\bf Z}_p,p^l1_{{\bf Z}_p})$ over ${\bf F}_p$. As for modules, let ${\got C}^*$ be the dual of ${\got C}$ (one could call it a latticed $F$-isocrystal over $S$). We also define the Tate twist ${\got C}^*(l)$ of ${\got C}^*$ by $(l)$; it is an $F$-crystal over $S$ (as $l>h_{\got C}$). In a similar way we define ${\got C}^*_0(l)$. As $l>h_0$, ${\got C}_{0}^*(l)$ is a Dieudonn\'e--Fontaine $p$-divisible object over ${\bf F}_p$ with non-negative slopes. 
\smallskip
We repeat the constructions we performed for ${\got C}_0$ and ${\got C}$ (like the ones through which we got $l$, $U_0$, $i_{S}(m)$, etc.) in the context of ${\got C}_0^*(l)$ and ${\got C}^*(l)$. So by enlarging $l$ and by performing $(\Mr1)$ and $(\Mr2)$, we can assume there exists a morphism 
$$i^*_{S}(m):{\bf E}({\got C}^*_0(l);W_m(S))\to {\bf E}({\got C}^*(l);W_m(S))$$
of $\Mm(W_m(S))$ whose cokernel is generically annihilated by $p^{3l}$ and which is the analogue of $i_{S}(m)$. As $l>h_{\got C}$, we think of ${\bf E}({\got C}^*(l);W_m(S))$ to be a ``twisted dual" of ${\bf E}({\got C};W_m(S))$ in the sense that there exists a morphism 
$$j_{S}(m-l):{\bf E}({\got C};W_{m-l}(S))\to {\bf E}({\got C}_0;W_{m-l}(S))$$
which at the level of $\Mo_{W_{m-l}(S)}$-modules is the dual of $i^*_{S}(m-l)$. Thus the pull back of $j_{S}(m-l)$ to an object of $\Mm(W_{m-l}(k))$ has a cokernel annihilated by $p^{3l}$. By performing $(\Mr_1)$ we can assume ${\rm Coker}(j_{S}(m-l))$ is generically annihilated by $p^{3l}$, cf. 2.8.3 (a). 
\medskip
{\bf 6.1.7. End of the proof of Theorem 6.1.} We will use the existence of the morphisms $i_S(m)$ and $i_S^*(m)$ to show that the assumption $U\neq S$ leads to a contradiction. Let $y:{\rm Spec}(k_1)\to S$ be a geometric point that does not factor through $U$. Let 
$$c_S(m-l):=j_{S}(m-l)\circ i_{S}(m-l):{\bf E}({\got C}_0;W_{m-l}(S))\to {\bf E}({\got C}_0;W_{m-l}(S)).$$
\indent
We check that ${\rm Coker}(c_S(m-l))$ is annihilated by $p^{7l}$. Let $c_S(m-l)_{\rm gen}$ and $c_S(m-2l)_{\rm gen}$ be the morphisms of $\Mm(W_{m-l}(k))$ and $\Mm(W_{m-2l}(k))$ (respectively) that are the natural pulls back of $c_S(m-l)$ and $c_S(m-2l)$ (respectively). As $c_S(m-2l)_{\rm gen}$ is the composite of two morphisms of $\Mm(W_{m-2l}(k))$ whose cokernels are annihilated by $p^{3l}$, ${\rm Coker}(c_S(m-2l)_{\rm gen})$ is annihilated by $p^{6l}$. As $c_S(m-2l)_{\rm gen}$ lifts to $c_S(m-l)_{\rm gen}$ and as $l>h_0r_0^2$, $c_S(m-2l)_{\rm gen}$ is the pull back of a morphism of $\Mm(W_{m-2l}({\bf F}_{p^{q_0}}))$ (cf. 2.2.3 (a) applied with $(K,k)$ replaced by $(k,{\bf F}_p)$). So as ${\bf F}_{p^{q_0}}\hookrightarrow R$, we get that ${\rm Coker}(c_S(m-2l))$ itself is annihilated by $p^{6l}$. Thus $p^{7l}$ annihilates ${\rm Coker}(c_S(m-l))$. 
\smallskip
As the endomorphism $y^*(c_S(m-l))=y^*(j_S(m-l))\circ y^*(i_S(m-l))$ of $\Mm(W_{m-l}(k))$ has a cokernel annihilated by $p^{7l}$, we get that $p^{7l}$ annihilates ${\rm Coker}(y^*(i_S(m-l)))$. Let $f_y:{\got C}_{0k_1}\to y^*({\got C})$ be a morphism 
that lifts $y^*(i_S(m-l-v))$, cf. 5.1.1 (b) and 5.1.2. As $m-v-l=7l+n+1\ge 7l+1$ (cf. the definition of $m$ in 6.1.3) and as $p^{7l}$ annihilates ${\rm Coker}(y^*(i_S(m-l-v)))$, (by reasons of ranks) the morphism $f_{y}$ is injective and so an isogeny. So $y^*({\got C})$ has Newton polygon $\Mn$. So $y$ factors through $U$. This contradicts the choice of $y$. Thus the existence of the morphisms $i_S(m)$ and $i_S^*(m)$ implies that $U=S$. As $U=S$, $U$ is affine. This ends the proof of Theorem 6.1.\endproof
\medskip\smallskip
{\bf 6.2. Proof of Main Theorem B.} We prove 1.6. Let $U$ be a reduced, locally closed subscheme of $S$ that is a stratum of $\Ms({\got C})$. We have to show that $U$ is an affine $S$-scheme. It suffices to check this under the extra assumptions that $S={\rm Spec}(R)$ is affine, that $U$ is an open, dense subscheme of $S$ and that the underlying $R$-module of ${\bf E}({\got C};W_1(S))$ is free. We will show that $U$ is an affine scheme. It suffices to check this under the extra assumption that $R$ is normal and perfect, cf. 2.9.2 applied with $X$ and $X^\prime$ replaced by $S$ and by the normalization of ${\rm Spec}(R^{\rm perf})$ (respectively). 
\smallskip
Let $\Mn$ be the Newton polygon of pulls back of ${\got C}_U$ via geometric points of $U$. Let ${\got C}_0$, $h_0$, $r_0$, $v$, and $n$ be associated to $\Mn$ and $h_{\got C}$ as in 5.4 (this makes sense even if $R$ is not an integral domain; see 2.1 for $h_{\got C}$). So ${\got C}_0=(M_0,\vph_0)$ is a Dieudonn\'e--Fontaine $p$-divisible object over ${\bf F}_p$ that has rank $r_0$, has $h$-number $h_0$, and has Newton polygon $\Mn$. 
\smallskip
Let $\tilde n$ be the maximum between $2+h_{\got C}+2\max\{d(r_0^2,s,c)|s\in S(0,h_{\got C}),\,c\in S(0,2h_{\got C})\}$ and $n$. If $(M,\vph)$ is an $F$-crystal over an algebraically closed field $K$ of characteristic $p$ of rank $r_0$ and $h$-number at most $h_{\got C}$, then $({\rm End}(M),\vph)$ is a latticed $F$-isocrystal over $K$ whose rank is $r_0^2$, whose $s$-number is at most $h_{\got C}$, and whose $h$-number is at most $2h_{\got C}$ (see end of 2.2.1 (e)). Thus from 3.2.8 (a) (applied with $G={\bf GL}_M$) and 2.4.1 we get:
\medskip
{\bf (i)} any $F$-crystal over $K$ whose rank is $r_0$ and whose $h$-number is at most $h_{\got C}$, is uniquely determined up to isomorphism by its reduction mod $p^{\tilde n}$. 
\medskip
We consider quadruples of the form $(\tilde k,\tilde\Phi,\tilde M,\tilde\vph)$, where:
\medskip\noindent
--  $\tilde k$ is an algebraically closed field of characteristic $p$, 
\smallskip\noindent
-- $\tilde\Phi$ is the Frobenius endomorphism of $W(\tilde k)[[w]]$ that is compatible with $\sg_{\tilde k}$ and that takes $w$ into $w^p$, 
\smallskip\noindent
-- $\tilde M$ is a free $W(\tilde k)[[w]]$-module of rank $r_0$ equipped with a $\tilde\Phi$-linear endomorphism $\tilde\vph$, 
\medskip\noindent
which have the property that the Newton polygons and the $h$-numbers of extensions of $(\tilde M,\tilde\vph)$ via $W(\tilde k)$-homomorphisms $W(\tilde k)[[w]]\to W(K)$ that are compatible with the Frobenius endomorphisms and that involve algebraically closed fields $K$ of characteristic $p$, are $\Mn$ and respectively are at most $h_{\got C}$. 
\medskip
We consider the unique $W(\tilde k)$-monomorphism $W(\tilde k)[[w]]\hookrightarrow W(\tilde k[[w]]^{\rm perf})$ that lifts the natural inclusion $\tilde k[[w]]\hookrightarrow \tilde k[[w]]^{\rm perf}$ and that is compatible with the Frobenius endomorphisms $\tilde\Phi$ and $\Phi_{\tilde k[[w]]^{\rm perf}}$, cf. [22, p. 145]; it maps $w$ into $(w,0,0,\ldots)\in W(\tilde k[[w]]^{\rm perf})$. 
\smallskip
The results [22, 2.6.1, 2.6.2, 2.7.1, and 2.7.4] hold as well in the context of pairs of the form $(\tilde M,\tilde\vph)$ that are not endowed with connections (one only has to disregard all details of loc. cit. that pertain to connections). So as in the proof of 6.1.2 we argue that there exists a number $l\in{\bf N}$ which has the properties that $l\ge {\rm max}\{h_0r_0^2,h_{\got C}\}$ and that for any quadruple $(\tilde k,\tilde\Phi,\tilde M,\tilde\vph)$ as above, there exists a monomorphism 
$$(M_0\otimes_{{\bf Z}_p} W(\tilde k[[w]]^{\rm perf}),\vph_0\otimes\Phi_{\tilde k[[w]]^{\rm perf}})\hookrightarrow (\tilde M\otimes_{W(\tilde k)[[w]]} W(\tilde k[[w]]^{\rm perf}),\tilde\vph\otimes\Phi_{\tilde k[[w]]^{\rm perf}})$$ 
whose cokernel is annihilated by $p^l$. Let $m:=8l+\tilde n+v+1$ and $\tilde m:=m+2v$. 
\smallskip
Let $(O,\vph_O):={\rm proj.}\,{\rm lim.}\,_{t\in{\bf N}} {\bf E}({\got C};W_t(S))$. As the underlying $R$-module of the object ${\bf E}({\got C};W_1(S))$ is free, $O$ is a free $W(R)$-module of rank $r_0$. Moreover $\vph_O$ is a $\Phi_R$-linear endomorphism of $O$. As $R$ is perfect, the $W_t(R)$-module of differentials $\Omega_{W_t(R)}$ is trivial. So the connection on $O$ induced by ${\got C}$ is trivial. From this and [22, p. 145] we get that the pair $(O,\vph_O)$ determines  ${\got C}$ up to isomorphism. Let $\Mb=\{e_1,\ldots,e_{r_0}\}$ be a $W(R)$-basis of $O$. Let $B\in M_{r_0\times r_0}(W(R))$ be the matrix representation of $\vph_O$ with respect to $\Mb$. Let $B_{\tilde m}\in M_{r_0\times r_0}(W_{\tilde m}(R))$ be $B$ mod $p^{\tilde m}$. Let $R^0$ be a finitely generated ${\bf F}_p$-subalgebra of $R$ which contains the components of the Witt vectors of length $\tilde m$ with coefficients in $R$ that are entries of $B_{\tilde m}$; so $B_{\tilde m}\in M_{r_0\times r_0}(W_{\tilde m}(R^0))$. Let $S^0:={\rm Spec}(R^0)$.
\smallskip
Let $\vph_O^\prime$ be a $\Phi_R$-linear endomorphism of $O$ whose matrix representation with respect to $\Mb$ is a matrix $B^\prime\in M_{r_0\times r_0}(W(R^0))\subseteq M_{r_0\times r_0}(W(R))$ that lifts $B_{\tilde m}$. Let ${\got C}_S^\prime$ be the $F$-crystal over $S$ that corresponds to the pair $(O,\vph_O^\prime)$, cf. the above part that refers to [22, p. 145]. As $B^\prime\in M_{r_0\times r_0}(W(R^0))$, loc. cit. also implies that ${\got C}_S^\prime$ is the pull back to $S$ of an $F$-crystal ${\got C}^\prime$ over ${\rm Spec}(R^{0\rm perf})$. As $B$ and $B^\prime$ are congruent mod $p^{\tilde m}$, we can identify ${\got C}/p^{\tilde m}{\got C}$ with ${\got C}_S^\prime/p^{\tilde m}{\got C}_S^\prime$. It is easy to see that due to this identification and to (i), the two Newton polygon stratifications $\Ms({\got C})$ and $\Ms({\got C}_S^\prime)$ of $S$ coincide. Thus to prove that $U$ is affine we can assume that $B=B^\prime$, $R=R^{0\rm perf}$, and ${\got C}={\got C}^\prime$. As $R=R^{0\rm perf}$, there exists a unique open subscheme $U^0$ of $S^0$ such that we have $U=U^0\times_{S^0} S$. 
\smallskip
To prove that $U$ is affine, it suffices to show that $U^0$ is affine. The scheme $U^0$ is affine if and only if its intersection with any irreducible component $C^0$ of $S^0$ is affine, cf. [16, Ch. II, Cor. (6.7.3)]. Thus by replacing $(S^0,S)$ with $(C^0,C^{0{\rm perf}})$, we can assume that both $S^0$ and $S$ are integral schemes. By replacing $S^0$ and $S$ with their normalizations (cf. 2.9.2), we can also assume that $S^0$ is a normal ${\rm Spec}({\bf F}_p)$-scheme of finite type.
\smallskip
As $B_{\tilde m}\in M_{r_0\times r_0}(W_{\tilde m}(R^0))$, for any $j\in S(1,\tilde m)$ and for every $S^0$-scheme $S_1$ we can speak about the object ${\bf E}({\got C};W_j(S_1))$ of $\Mm(W_j(S_1))$ whose underlying $\Mo_{W_j(S_1)}$-module is the free $\Mo_{W_j(S_1)}$-module that has $\Mb$ as an $\Mo_{W_j(S_1)}$-basis and whose underlying Frobenius endomorphism has a matrix representation with respect to $\Mb$ which is the natural image of $B_{\tilde m}$ in $M_{r_0\times r_0}(\Mo_{W_j(S_1)})$. If $S_1$ is an $S$-scheme, then ${\bf E}({\got C};W_j(S_1))$ is precisely the object of $\Mm(W_j(S_1))$ defined in 2.8.2.
\smallskip
Let $V^0$ be an arbitrary local ring of $U^0$ that is a discrete valuation ring. Let $w_0$ be a uniformizer of it. Let $V^0_3$ be a $V^0$-algebra that is a complete discrete valuation ring, that has $w_0$ as a uniformizer, and that has an algebraically closed residue field $k_3$. So $V^0_3$ is isomorphic to $k_3[[w_0]]$, with $w_0$ viewed as a variable. We identify $w_0^{p^{-\tilde m}}$ with a uniformizer of $V_2^0:=V_3^{0(p^{\tilde m})}$. So $k_2:=k_3^{(p^{\tilde m})}$ is the residue field of $V_2^0$. Let $V_1^0:=k_2[[w_0^{p^{-\tilde m}}]]^{\rm perf}=V_2^{0\rm perf}$. For $j\in S(1,\tilde m)$ let $W^0_{j}:=W(k_2)[[w_0^{p^{-j}}]]$ be endowed with the Frobenius endomorphism $\Phi_{W^0_j}$ that is compatible with $\sg_{k_2}$ and that takes $w_0^{p^{-j}}$ into $w_0^{p^{-j+1}}$. Let $f^0_j:W^0_j\hookrightarrow W(k_2[[w_0^{p^{-j}}]])$ be the $W(k_2)$-monomorphism that lifts the canonical identification  $k_2[[w_0^{p^{-j}}]]=W_1(k_2[[w_0^{p^{-j}}]])$ and that takes $w_0^{p^{-j}}$ into the Witt vector $(w_0^{p^{-j}},0,0,\ldots)\in W(k_2[[w_0^{p^{-j}}]])$. The following two properties hold:
\medskip
{\bf (ii)} each $f^0_j$ is compatible with Frobenius endomorphisms, and 
\smallskip
{\bf (iii)} if $j<\tilde  m$, the restriction of $f^0_{j+1}$ to the $W(k_0)$-subalgebra $W^0_j$ of $W^0_{j+1}$ is $f_j^0$. 
\medskip
We recall that if $x=(x_0,x_1,\ldots,x_{\tilde m})$ is a Witt vector of length $\tilde m$, then $px=(0,x_0^p,x_1^p,\ldots,x_{\tilde m -1}^p)$. Based on this and (ii) and (iii), by induction on $j\in S(1,\tilde m)$ we get that the image (via the natural monomorphism $R^0\hookrightarrow k_2[[w_0^{p^{-j}}]]$) of the matrix $B=B^\prime\in M_{r_0\times r_0}(W(R^0))$ in $M_{r_0\times r_0}(W_j(k_2[[w_0^{p^{-j}}]]))$, belongs to $M_{r_0\times r_0}(W^0_j/p^jW^0_j)$. Thus the image of $B=B^\prime$ in $M_{r_0\times r_0}(W_{\tilde m}(k_2[[w_0^{p^{-\tilde m}}]]))$, belongs to $M_{r_0\times r_0}(W^0_{\tilde m}/p^{\tilde m}W^0_{\tilde m})$ and so it lifts to a matrix $\tilde B\in M_{r_0\times r_0}(W^0_{\tilde m})$. 
\smallskip
Let $\tilde M:=\oplus_{i=1}^{r_0} W^0_{\tilde m}e_i$. Let $\tilde\vph$ be the $\Phi_{W^0_{\tilde m}}$-linear endomorphism of $\tilde M$ whose matrix representation with respect to $\Mb$ is $\tilde B$. The extension of $(\tilde M,\tilde\vph)$ via a $W(k_2)$-homomorphism $W^0_{\tilde m}\to W(K)$ compatible with Frobenius endomorphisms, has the $h$-number at most $h_{\got C}$ (as $\tilde m>h_{\got C}$) and has Newton polygon $\Mn$ (cf. (i) and the fact that $V^0$ is a local ring of $U^0$). So there exists a monomorphism 
$$(M_0\otimes_{{\bf Z}_p} W(V_1^0),\vph_0\otimes\Phi_{V_1^0})\hookrightarrow (\tilde M\otimes_{W^0_{\tilde m}} W(V_1^0),\tilde\vph\otimes\Phi_{V_1^0})$$ 
whose cokernel is annihilated by $p^l$, cf. the choice of $l$ (applied with $k_2[[w_0^{p^{-\tilde m}}]]$ instead of $\tilde k[[w]]$). Thus there exists a morphism ${\bf E}({\got C}_0;W_{\tilde m}(V^0_1))\to {\bf E}({\got C};W_{\tilde m}(V^0_1))$ of $\Mm(W_{\tilde m}(V_1^0))$ whose cokernel is annihilated by $p^l$. 
\smallskip
As in Subsubsections 6.1.3 to 6.1.7 we only used evaluation functors ${\bf E}$ and pulls back of $F$-crystals over $S$ via geometric points of $S$ and as for any algebraically closed field $K$ the map $S(K)\to S^0(K)$ is bijective, the rest of the proof that $U^0$ is affine is the same as Subsubsections 6.1.3 to 6.1.7 (but with the role of $(n,m)$ being replaced with the one of $(\tilde n,\tilde m)$). We will only add two extra sentences. 
\smallskip
From 2.8.3 (c) (applied with $q$ replaced by $\tilde m$) we get that there exist: 
\medskip
{\bf (iv)} an open subscheme $U^0_{\tilde V^0}$ of the normalization of $U^0$ in a finite field extension $k_{S^0,V^0}$ of the field of fractions $k_{S^0}$ of $S^0$, a local ring $\tilde V^0$ of $U^0_{\tilde V^0}$ which dominates $V^0$, and a morphism $i_{U^0_{\tilde V^0}}(\tilde m):{\bf E}({\got C}_0;W_{\tilde m}(U^0_{\tilde V^0}))\to {\bf E}({\got C};W_{\tilde m}(U^0_{\tilde V^0}))$ of $\Mm(W_{\tilde m}(U^0_{\tilde V^0}))$ whose cokernel is annihilated by $p^l$. 
\medskip
In connection with the last two paragraphs of 6.1.3 and with Subsubsections 6.1.4 to 6.1.7, we only have to add an upper right index $0$ to all schemes that (modulo the two operations of 6.1.1) are about to be introduced; thus we get open subschemes $U^0_{V^0}$, $U_i^0$, and $U^0_0$ of $U^0$, etc. This ends the proof of 1.6. 
\medskip\smallskip
{\bf 6.3. Remarks.} {\bf (a)} Let $S$ be an integral,  locally noetherian scheme. Let ${\got C}$ and $\Ms({\got C})$ be as in the beginning of \S6. Let $U$ be the unique stratum of $\Ms({\got C})$ which is an open subscheme of $S$. The open embedding $U\hookrightarrow S$ is an affine morphism, cf. 6.1. This implies that either $U=S$ or $S^{\rm top}\setminus U^{\rm top}$ is of pure codimension 1 in $S^{\rm top}$. It suffices to check this statement under the extra assumptions  that (to be compared with the first paragraph of the proof of 6.1)  $S$ is also local, complete, and normal and that $S^{\rm top}\setminus U^{\rm top}$ has pure codimension $c\in {\bf N}$ in $S^{\rm top}$. If $c>1$, then by applying 2.9.1 to the affine, birational open embedding $U\hookrightarrow  S$ we get that $U=S$. Thus $c$ must be $1$.
\smallskip
Thus 6.1 implies the following result of de Jong and Oort (see [10, 4.1]): if $S$ is a local, integral, noetherian ring and if $U$ contains the complement in $S$ of the closed point of $S$, then either the dimension of $S$ is at most $1$ or $S=U$. The converse of this implication holds, provided our scheme $S$ is locally factorial. But in general the result of de Jong and Oort does not imply 6.1. This is so as there exist integral, normal, noetherian, affine schemes $S={\rm Spec}(R)$ that have a prime Weil divisor $C$ such that the open subscheme $S\setminus C$ of $S$ is not an affine scheme. Here is one classical example. 
\smallskip
Let $R:=k[x_1,x_2,x_3,x_4]/(x_1x_4-x_2x_3)$. Let $C:={\rm Spec}(k[x_1,x_3])$ be the irreducible divisor of $S$ defined by the equations $x_4=x_2=0$. The open subscheme $S\setminus C$ of $S$ is the union of ${\rm Spec}(R[{1\over {x_2}}])$ and ${\rm Spec}(R[{1\over {x_4}}])$ and thus its $R$-algebra of global functions is $R[{1\over {x_2}}]\cap R[{1\over {x_4}}]$. But $W:={\rm Spec}(R[{1\over {x_2}}]\cap R[{1\over {x_4}}])$ is an affine $S$-scheme whose fibre over the point of $S$ defined by $x_1=x_2=x_3=x_4=0$ is non-empty. Thus the natural morphism $S\setminus C\to W$ is not an isomorphism and so the scheme $S\setminus C$ is not affine.
\smallskip
{\bf (b)} Let $(M,\vph)$ be an $F$-crystal over a perfect field $k$ of characteristic $p$. It is easy to see that [22, 1.4 and 1.5] implies the existence of a number $n_0\in{\bf N}$ such that for any $g\in {\bf GL}_M(W(k))$ the Newton polygon of $(M,g\vph)$ depends only on $g$ mod $p^{n_0}$. For instance, if $k=\bar k$ we can take $n_0$ to be the number $n_{\rm fam}$ of 3.1.5 for $G={\bf GL}_M$. One can use this (in a way similar to the first part of 6.2) to define Newton polygon stratifications for reductions modulo adequate powers of $p$ of $F$-crystals over reduced ${\rm Spec}({\bf F}_p)$-schemes. 
\smallskip
For instance, it can be easily checked starting from 1.3 and [21, 4.4 e)] that any truncated Barsotti--Tate group $\Mg_S$ of level $T(r,d)$ over a reduced ${\rm Spec}({\bf F}_p)$-scheme $S$ which has height $r$ and relative dimension $d$, defines a stratification $\Ms(\Mg_S)$ of $S$ as follows.  The association $\Mg_S\to\Ms(\Mg_S)$ is uniquely determined by the following two properties: 
\medskip
{\bf (i)} it is functorial with respect to pulls back, and 
\smallskip
{\bf (ii)} if there exists a $p$-divisible group $\Md_S$ over $S$ such that $\Md_S[p^{T(r,d)}]$ is isomorphic to $\Mg_S$, then $\Ms(\Mg_S)$ is the Newton polygon stratification of $S$ defined by the $F$-crystal over $S$ that is associated naturally to $\Md_S$.
\bigskip\smallskip
\centerline{\bigsll {\bf References}}
\medskip\smallskip
\item{[1]} A. Borel, {\it Linear algebraic groups}, Grad. Texts in Math., Vol. {\bf 126}, Springer-Verlag, 1991.
\item{[2]} S. Bosch, W. L\"utkebohmert, and M. Raynaud, {\it N\'eron models}, Ergebnisse der Math. und ihrer Grenzgebiete (3), Vol. {\bf 21}, Springer-Verlag, 1990. 
\item{[3]} J. H. Conway, R. T. Curtis, S. P. Norton, R. A. Parker, and R. A. Wilson, {\it Atlas of finite groups}, xxxiv+252 pp., Oxford Univ. Press, Eynsham, 1985.
\item{[4]} J. E. Cremona, {\it Algorithms for modular elliptic curves}, Cambridge Univ. Press, Cambridge, 1992.
\item{[5]} M. Demazure, {\it Lectures on $p$-divisible groups}, Lecture Notes in Math., Vol. {\bf 302}, Springer-Verlag, 1972.
\item{[6]} M. Demazure, A. Grothendieck, et al., {\it Sch\'emas en groupes}, Vol. I--III, Lecture Notes in Math., Vol. {\bf 151--153}, Springer-Verlag, 1970.
\item{[7]} J. Dieudonn\'e, {\it Groupes de Lie et hyperalg\`ebres de Lie sur un corps de caract\'erisque $p>0$ (IV)}, Amer. J. Math. {\bf 77} (1955), pp. 429--452.
\item{[8]} J. de Jong, {\it Finite locally free group schemes in characteristic $p$ and Dieudonn\'e modules}, Invent. Math. {\bf 114} (1993), no. 1, pp. 89--137.
\item{[9]} J. de Jong, {\it Crystalline Dieudonn\'e module theory via formal and rigid geometry}, Inst. Hautes \'Etudes Sci. Publ. Math., Vol. {\bf 82},  pp. 5--96, 1995.
\item{[10]} J. de Jong and F. Oort, {\it Purity of the stratification by Newton polygons}, J. of Amer. Math. Soc. {\bf 13} (2000), no. 1, pp. 209--241.
\item{[11]} G. Faltings, {\it Crystalline cohomology and $p$-adic Galois representations}, Algebraic analysis, geometry, and number theory (Baltimore, MD, 1988), pp. 25--80, Johns Hopkins Univ. Press, Baltimore, MD, 1989.
\item{[12]} G. Faltings, {\it Integral crystalline cohomology over very ramified valuation rings}, J. of Amer. Math. Soc. {\bf 12} (1999), no. 1, pp. 117--144.
\item{[13]} G. Faltings and C.-L. Chai, {\it Degeneration of abelian varieties},  Ergebnisse der Math. und ihrer Grenzgebiete (3), Vol. {\bf 22}, Springer-Verlag, 1990.
\item{[14]} J.-M. Fontaine, {\it Groupes $p$-divisibles sur les corps locaux}, J. Ast\'erisque {\bf 47/48}, Soc. Math. de France, Paris, 1977.
\item{[15]} J.-M. Fontaine and G. Laffaille, {\it Construction de repr\'esentations p-adiques}, Ann. Sci. \'Ecole Norm. Sup. {\bf 15} (1982), no. 4, pp. 547--608.
\item{[16]} A. Grothendieck, {\it \'El\'ements de g\'eom\'etrie alg\'ebrique. II. \'Etude globale \'el\'ementaire de quelques classes de morphisms}, Inst. Hautes \'Etudes Sci. Publ. Math., Vol. {\bf 8}, 1961.
\item{[17]} A. Grothendieck, {\it \'El\'ements de g\'eom\'etrie alg\'ebrique. III. \'Etude cohomologique des faisceaux coh\'erents, Premi\`ere partie}, Inst. Hautes \'Etudes Sci. Publ. Math., Vol. {\bf 11}, 1963.
\item{[18]} A. Grothendieck, {\it \'El\'ements de g\'eom\'etrie alg\'ebrique. IV. \'Etude locale des sch\'emas et des morphismes de sch\'ema}, Inst. Hautes \'Etudes Sci. Publ. Math., Vols. {\bf 20} (1964), {\bf 24} (1965), {\bf 28} (1966), and {\bf 32} (1967).
\item{[19]} A. Grothendieck, {\it Groupes de Barsotti--Tate et cristaux de Dieudonn\'e}, S\'em. Math. Sup. {\bf 45} (\'Et\'e, 1970), Les presses de l'Universit\'e de Montr\'eal, Montreal, Quebec, 1974.
\item{[20]} S. Helgason, {\it Differential geometry, Lie groups, and symmetric spaces}, Academic Press, New-York, 1978.
\item{[21]} L. Illusie, {\it D\'eformations des groupes de Barsotti--Tate (d'apr\`es A. Grothendieck)}, Seminar on arithmetic bundles: the Mordell conjecture (Paris, 1983/84), pp. 151--198, J. Ast\'erisque {\bf 127}, Soc. Math. de France, Paris, 1985.
\item{[22]} N. Katz, {\it Slope filtration of $F$-crystals}, Journ\'ees de G\'eom\'etrie Alg\'ebrique de Rennes (Rennes, 1978), Vol. I, pp. 113--163, J. Ast\'erisque {\bf 63}, Soc. Math. de France, Paris, 1979. 
\item{[23]} N. Katz, {\it Serre--Tate local moduli},  Algebraic surfaces (Orsay, 1976--78), pp. 138--202, Lecture Notes in Math., Vol. {\bf 868}, Springer-Verlag, 1981. 
\item{[24]} R. E. Kottwitz, {\it Isocrystals with additional structure}, Comp. Math. {\bf 56} (1985), no. 2, pp. 201--220.
\item{[25]} R. E. Kottwitz, {\it Points on some Shimura varieties over finite fields}, J. of Amer. Math. Soc. {\bf 5} (1992), no. 2, pp. 373--444.
\item{[26]} H. Kraft, {\it Kommutative algebraische p-Gruppen (mit Anwendungen auf p-divisible Gruppen und abelsche Variet\"aten)}, manuscript, 86 pages, Univ. Bonn, 1975.
\item{[27]} K.-Z. Li and F. Oort, {\it Moduli of supersingular abelian varieties}, Lecture Notes in Math., Vol. {\bf 1680}, Springer-Verlag, 1998.
\item{[28]} J. I. Manin, {\it The theory of formal commutative groups in finite characteristic}, Russian Math. Surv. {\bf 18} (1963), no. 6, pp. 1--83. 
\item{[29]} H. Matsumura, {\it Commutative algebra. Second edition}, The Benjamin/Cummings Publ. Co., Inc., Reading, Massachusetts, 1980.
\item{[30]} W. Messing, {\it The crystals associated to Barsotti--Tate groups, with applications to abelian schemes}, Lecture Notes in Math., Vol. {\bf 264}, Springer-Verlag, 1972.
\item{[31]} J. S. Milne, {\it The points on a Shimura variety modulo a prime of good reduction}, The zeta functions of Picard modular surfaces, pp. 153--255, Univ. Montr\'eal Press, Montreal, Quebec, 1992.
\item{[32]} D. Mumford, {\it Abelian varieties}, viii+242 pp, Tata Inst. of Fund. Research Studies in Math., No. {\bf 5}, Published for the Tata Institute of Fundamental Research, Bombay; Oxford Univ. Press, London, 1970 (reprinted 1988).
\item{[33]} D. Mumford, J. Fogarty, and F. Kirwan, {\it Geometric invariant theory. Third edition}, Ergebnisse der Math. und ihrer Grenzgebiete (2), Vol. {\bf 34}, Springer-Verlag, 1994.
\item{[34]} F. Oort, {\it A stratification of a moduli space of abelian varieties}, Moduli of abelian varieties (Texel Island, 1999),  pp. 345--416, Progr. Math., Vol. {\bf 195}, Birkh\"auser, Basel, 2001.
\item{[35]} F. Oort, {\it Newton polygon strata in the moduli space of abelian varieties}, Moduli of abelian varieties (Texel Island, 1999),  pp. 417--440, Progr. Math., Vol. {\bf 195}, Birkh\"auser, Basel, 2001.
\item{[36]} F. Oort, {\it Foliations in moduli spaces of abelian varieties}, J. of Amer. Math. Soc. {\bf 17} (2004), no. 2, pp. 267--296.
\item{[37]} M. Rapoport and M. Richartz, {\it On the classification and specialization of $F$-isocrystals with additional structure}, Comp. Math. {\bf 103} (1996), no. 2, pp. 153--181.
\item{[38]} J.-P. Serre, {\it Galois Cohomology}, Springer-Verlag, 1997.
\item{[39]} J. Tate, {\it Classes d'isog\'enie des vari\'et\'es sur un corps fini (d'apr\`es J. Honda)}, S\'em. Bourbaki 1968/69, Exp. 352, Lecture Notes in Math., Vol. {\bf 179}, pp. 95--110, Springer-Verlag, 1971.
\item{[40]} A. Vasiu, {\it Integral canonical models of Shimura varieties of preabelian
 type}, Asian J. Math. {\bf 3} (1999), no. 2, pp. 401--518.
\item{[41]} T. Zink, {\it On the slope filtration}, Duke Math. J. {\bf 109} (2001), no. 1, pp. 79--95.
\item{[42]} J.-P. Wintenberger, {\it Un scindage de la filtration de Hodge pour certaines
variet\'es algebriques sur les corps locaux}, Ann. of Math. (2) {\bf 119} (1984), no. 3, pp. 511--548.
\smallskip
$\;\;\;$\hbox{Adrian Vasiu}

$\;\;\;$\hbox{University of Arizona}

$\;\;\;$\hbox{Department of Mathematics}

$\;\;\;$\hbox{617 N. Santa Rita, P.O. Box 210089}

$\;\;\;$\hbox{Tucson, AZ-85721-0089, USA}

$\;\;\;$\hbox{E-mail: adrian@math.arizona.edu}
\end